\documentclass[11pt,a4paper]{article}
\usepackage{amsmath}
\usepackage{amssymb,amsfonts, amstext,amsmath}
\usepackage[latin1]{inputenc}
\usepackage{theorem}

\usepackage{graphicx}
\usepackage{xcolor}
\usepackage{fancybox}

\usepackage{hyperref}

\newcounter{exo}

\newcommand{\Q}{\mathbb{Q}}
\renewcommand{\P}{\mathcal{P}}

\newcommand{\R}{\mathbb{R}}
\newcommand{\Z}{\mathbb{Z}}
\newcommand{\N}{\mathbb{N}}
 
 \newcommand{\F}{\mathcal{F}}

\newcommand{\A}{\mathcal{A}}

\renewcommand{\S}{\mathbf{\mathcal{S}}}

\newcommand{\Irr}{\mathrm{Irr}}
\newcommand{\PF}{\mathrm{PF}}
\newcommand{\CFE}{\mathcal{C}}
\newcommand{\Gen}{\mathrm{Gen}}

\renewcommand{\phi}{\varphi}

 \newcommand{\esp}{\hspace*{0.2cm}}

\newcommand{\titre}[1]{\begin{center}  \textbf{\Large{#1}} \\ \end{center}}

\newcommand{\inter}[2]{ \{  #1,..,#2 \} }

 \theorembodyfont{\sl}  \newtheorem{montheo}{Theorem}
 \theorembodyfont{\normalfont}  
     \newtheorem{madef}{Definition} 
       \newtheorem{monalgo}{Algorithm} 
    
      \newtheorem{maprop}{Proposition}      
         \newtheorem{moncoro}{Corollary}     
                           
 \theorembodyfont{\sl} \newenvironment{mademo}{\textbf{Proof : } \small \\}{\normalsize $\blacksquare$ \\ }
 \newtheorem{monlem}{Lemma}

\begin {document}
 
 \titre{Quotients of numerical semigroups generated by two numbers }
 
 \hspace*{5.7cm} \large{Emmanuel Cabanillas \footnote{e-mail : emmanuel.cabanillas@ac-montpellier.fr}}

\vspace{0.5cm}

\normalsize 
 
\hspace{6.7cm} ABSTRACT : \\
 
\esp In this article, we study the quotients of numerical semigroups, generated by two coprime positive numbers, denoted $\frac{<a,b>}{d}$. We give  formulae for the usual invariants of these semigroups, expressed in terms of continued fraction expansions and Ostrowski-like numeration of some rationals, simply related to inputs $a,b,d$. So, we obtain  quadratic complexity algorithms to compute these invariants.  As a consequence, we will show that, for these numerical semigroups, the type is always lower than the embedding dimension and deduce Wilf's property.  We also consider the " reverse problem" : given $I$ a finite set of integers, we obtain an expression of all possible triplets $(a,b,d)$, such that $I$ is the set of minimal generators of $\frac{<a,b>}{d}$.\\

  \tableofcontents

  \newpage

 \section{Introduction}
  
     \subsection{basics about numerical semigroups}
     
     \label{subsec:intro}
 
 


  \esp We give here some definitions and elementary results without proofs. See  \textbf{[2],[8]} and \textbf{[11]} (  \nameref{biblio}) for more details. First, we will use some classical notations :  \\
  - we denote $\N$ the set of all non negative integers, $\N^*$ the set of positive integers and $\Z$ the set of integers.\\
  - for every couple of integers $p,q$, the notation $\inter{p}{q}$ means the set of integers $k$ such that $p\leqslant k\leqslant q$. We extend this notation to $\inter{p}{\infty}=\{ k\in \Z, k\geqslant p \}$. \\
   -  for an element $x$ and two subsets $I$ and $J$ of an additive group $G$, we will use the classical notations  :
    $$I+J=\{ i+j, i\in I,j\in J\} \esp ; \esp -I=\{ -i,i\in I\} \esp ; \esp  x+J=\{ x+j,j\in J\}$$ 
    - we will also denote for an integer $n$, a subset $K$ of $\Z$, an element $x$ of $G$ and a subset $J$ of $G$ :
    $$ nJ=\{ nj,j\in J \} \esp ; \esp xK=\{ kx, k\in K \}$$

$\bullet$  A \emph{numerical semigroup} $S$ is a cofinite submonoid of $(\N,+)$, where $\N$ denotes the set of non negative integers.  We can also define a numerical semigroup with a system of generators. Let $B$ be a non empty subset of $\N$  and $<B>$ the set of all non negative integer combinations of finite subsets of $B$. Then $<B>$ is a numerical semigroup ( \emph{generated by $B$}) if and only if the elements of $B$ are coprime ( not necesseraly pairwise coprime). In the case of finite $B$ ( which is the general case), say $B=\{ n_1,n_2,\cdots,n_r \}$, then :
 $$ <B> = \left\{ \sum_{k=1}^r x_kn_k \esp ; \esp  x_1,x_2,\cdots,x_r\in \N \right\} $$
 \esp \emph{Let $S$ be a numerical semigroup}.\\
 $\bullet$ An element $s$ of $S$ is \emph{irreducible} in $S$ if and only if there is no positive elements $s'$ and $s"$ of $S$ such that $s=s'+s"$. 
 The set of irreducible elements of $S$ will be denoted $\Irr(S)$. It is well known that $\Irr(S)$ is the  minimal set of generators of $S$ : that is to say $\Irr(S)$ generates $S$ and every set of generators of $S$ contains $\Irr(S)$. $\Irr(S)$ is finite and its  cardinality is called the \emph{embedding dimension of $S$}, denoted by $e(S)$.\\
 \esp $\Irr(S)$ is also the set of minimal elements of $S\backslash \{ 0 \}$, for the order $\leqslant_S $ induced by $S$ :
 $$ \forall x,y\in \Z, (x\leqslant_S y \Leftrightarrow y-x\in S )$$
 \esp  The smallest positive element of $S$ is called the \emph{multiplicity} of $S$ and is denoted $m(S)$. It is also the lowest element of $\Irr(S)$.\\
 
$\bullet$ The finite set $\N\backslash S$ is often named the \emph{set of gaps of $S$} and is denoted $G(S)$. Its cardinality, named the \emph{genus} of $S$, is denoted $g(S)$.  The greatest element of $G(S)$ is called the \emph{Frobenius number} of $S$ and is denoted $f(S)$, with the convention that $f(\N)=-1$. It is in general difficult to compute this Frobenius number, but a lot of investigations have been made ( \textbf{[11]} \nameref{biblio}).\\
 \esp The number $f(S)+1$ is sometimes called the \emph{conductor } of $S$ and denoted $c(S)$. The set of elements of $S$ lower than $f(S)$ will be denoted by $S_0$, so that we have the following partition of $\N$ : 
 $$ \N = G(S) \cup S_0 \cup  \inter{ c(S) }{ +\infty }$$
 \esp So, we have : 
    $$\#S_0=c(S)-g(S)$$
   $\bullet$  Note that if $s\in S$, then $f(S)-s\not\in S$, so that :
   $$f(S)-S_0\subset G(S)\esp \text{ so } \esp \# S_0 \leqslant g(S) \esp \text{ and } \esp \frac{g(S)}{c(S)}\geqslant \frac{1}{2}$$
   \esp A particular case if when we have an equality : then   $S$ is said \emph{symmetric}.\\
   
$\bullet$ The Frobenius number of $S$ is a particular case of the notion of \emph{pseudo-Frobenius numbers} ( PF numbers in this paper) of $S$, that can be defined as the maximal elements of $\Z\backslash S$ for the order $\leqslant_S $. $\PF(S)$ will  denote the set of these elements : it is a subset of $G(S)$ \footnote{ except when $S=\N$, where $G(S)=\emptyset$ and $\PF(S)=\{ -1 \}$}and $f(S)$ is its greatest element ( for the usual order on $\Z$). The cardinality of $\PF(S)$ is the \emph{type} of $S$, denoted $t(S)$.\\
\esp $\PF(S)$ is  also the minimal subset of $\Z$ that satisfies the following " symmetry property" :
$$ \Z\backslash S=\PF(S)-S=\bigcup_{p\in \PF(S)}(p-S)$$
\esp If we restrict to $\N$, then we obtain :
$$ G(S)\subset \PF(S)-S_0=\bigcup_{p\in \PF(S)}(p-S_0)$$

\esp The above inclusion proves that : $  g(S)\leqslant t(S) \# S_0$ and then :
$$ \frac{g(S)}{c(S)}\leqslant 1-\frac{1}{t(S)+1}$$

$\bullet$ \emph{Wilf's property} for a numerical semigroup $S$ is the following :
$$\frac{g(S)}{c(S)}\leqslant 1-\frac{1}{e(S)}$$
\esp \emph{Wilf's conjecture} claims that this inequality holds for every numerical semigroup ( \textbf{[17]}\nameref{biblio}). It has been proven in many cases ( \textbf{[5],[6],[7],[10],[15]} ), but the general case   still remains open.  According to the previous inequalities, it would be sufficient to prove that $t(S)<e(S)$. Unfortunately, it is not true for all numerical semigroups ( see example after Theorem 11 in \textbf{[7]}).\\
\esp Note that the ratio $g(S)/c(S)$ could be named : the " density of gaps of $S$".\\
 


 $\bullet$ The \emph{quotient of $S$ by a positive integer $d$} is defined by : 
 $$ S/d= \{ x \in \N, dx\in S \} $$
  \esp It is a numerical semigroup and we remark that $S/d=\N$ if and only if $d\in S$. We also have an obvious calculation rule, when we consider quotients of a quotient :
  $$ \forall d,d'\in \N^*, \esp \frac{S/d}{d'}=\frac{S}{dd'}$$
  \esp Many investigations have been made to relate invariants of $S/d$ to invariant of $S$ ( see introduction of \textbf{[1]} for details and further references).\\
  \esp In this paper, we will study the most simple case, that is, when $S$ is generated by two integers :   let $a$ and $b$ be two coprime integers greater than one. Since Sylvester in 1882 ( \textbf{[16]}), we know that $<a,b>$ is symmetric and we have simple formulae for $f$ and $g$ :
    $$ f(<a,b>)=ab-a-b \esp ;\esp g(<a,b>)= \frac{(a-1)(b-1)}{2}$$

   \esp  Curiously,  $<a,b>$ is very simple, but its quotients are not in general ! It was an open problem to express the invariants of $S= \frac{<a,b>}{d}$, in terms of $a,b,d$ ( \textbf{[4]}). A recent paper ( \textbf{[1]}\nameref{biblio}) gives an expression of $g(<a,b>/d)$ with a sum of integral parts, using Hilbert series. We will obtain another formulation for  $g(<a,b>/d)$ in Theorem \ref{theo:10}, which could easily be relied to the previous one. \\
   
   $\bullet$  Another motivation, for the study of these numerical semigroups, lies in the fact that they can be defined by two alternative ways : \\
    - let $a,b,c$ be positive integers, then the set of integers $x$ that verify $ax$ mod $ b \leqslant cx$ ( where mod denotes the remainder ) is a numerical semigroup, named  \emph{proportionally modular}.  In \textbf{[13]}, it is shown that : a numerical semigroup is  proportionally modular if and only if it exists three pairwise coprime positive integers $a,b,d$ such that $S= \frac{<a,b>}{d}$. Furthermore, it is proven  in \textbf{[12]}, that we can choose $b=a+1$. In \textbf{[13]}, they build a set of generators of these $\frac{<a,b>}{d}$, but it is not always minimal...\\
    - an equivalent definition is  : for $q$ and $q'$ two rationals such that $q<q'$, $\N\cap <[q,q']>$ defines a numerical semigroup, that is proportionally modular ( \textbf{ [9]}).\\
    
    $\bullet$ A question at the end of \textbf{[1]} is about the invariants of $S/d$, when $S$ is generated by an arithmetic sequence. We prove in our paper that, in that case, $S$ is itself a quotient of $<a,b>$ for some $a$ and $b$, which is a particular case of a more general result exposed in \textbf{[9]}( \nameref{biblio}) : a numerical semigroup $S$ is of the form $\frac{<a,b>}{d}$, for some $a,b,d$ if and only if $\Irr(S)$ can be enumerated by what we call a \emph{modular-convex} sequence of positive integers $(n_k)_{k\in\inter{0}{r}}$ :
    $$ \forall k\in \inter{1}{r-1}, \esp \frac{n_{k-1}+n_{k+1}}{n_k} \in \inter{2}{\infty}$$
    \esp We give a new proof of this result and  we also give explicit formulae for all possible $a,b$ and $d$ : Theorems \ref{theo:1} , \ref{theo:2} , \ref{theo:3} in section \ref{subsec:solutions} . When $S$ is generated by an arithmetic sequence ( which corresponds to the case when $ \frac{n_{k-1}+n_{k+1}}{n_k}=2$ for all $k$), we obtain ( see Corollary \ref{coro:1} at the end of \ref{subsec:solutions} ) : for $a,q,r$ positive integers such that $a$ and $q$ are coprime
    $$ <a,a+q,a+2q,\cdots,a+rq>= \frac{<a,a^2+q(ar+1)>}{ar+1}$$
    \esp This is a particular case of more general results and formulae, but it can be proven directly with elementary considerations. Combining this with Theorems \ref{theo:4} ,\ref{theo:5} , \ref{theo:6} , \ref{theo:7} , \ref{theo:8} , \ref{theo:9} and \ref{theo:10} ( see section \ref{sec:invariants}), that give expressions of invariants of $\frac{<a,b>}{d}$ and our remark about $(S/d)/d'=S/dd'$, we obtain the answer to the above question...

    \subsection{overview of the article}

   \esp First, in section \ref{sec:representation}, following an idea of \textbf{[13]}, we use a representation of a numerical semigroup $S$ as a subset of the cylinder $\Z^2/w\Z$, where $w=(b,-a)$ is a non null vector of $\Z^2$, via an additive map $(x,y)\to ax+by$, where  $a$ and $b$ are two coprime elements of $S$. For the numerical semigroups we are interseted in, that is $\frac{<a,b>}{d}$, this will give a " lattice-stair case shape" representation : the intersection of a lattice and a regular stair case shape subset of $\Z^2$ ( modulo $w\Z$). \\
   \esp This approach leads in \ref{subsec:relations} to simple relations between $G(S),\Irr(S),\PF(S)$ and  elements or minimal sets ( for the product order on $\Z^2$) of the above lattice in some rectangles. We remark in \ref{subsec:minimal} that these minimal sets can be parametrized by " double inductive" sequences of points, that are  convex and unimodular : we name them " modular-convex" ( see above). \\
   \esp In \ref{sec:reverse}, we relate these sequences with a kind of continued fraction expansion and use them to solve the " reverse problem" : we obtain in \ref{subsec:solutions} a complete expression of all triplets of pairwise coprime positive integers $a,b,d$ such that $1<a<b$ and $\Irr(\frac{<a,b>}{d})$ is a given finite subset of $\N$. \\
   \esp In order to solve the " direct problem" ( express invariants of $\frac{<a,b>}{d}$ in terms of $a,b,d$), we need another  parametrization of the above minimal sets by Kronecker sequences $(\{ n\alpha -\beta \},n)_n$. So, we use some general considerations about usual Diophantine approximation ( \ref{subsec:CFE} and \ref{subsec:semicv}), as well as some more unusual ( \ref{subsec:alphanum} ) Ostrowski-like numeration (  detailed in \textbf{[3]} \nameref{biblio}). \\
   \esp Finally, we give in section \ref{sec:invariants}. all formulae about the usual invariants of numerical semigroups $\frac{<a,b>}{d}$. We end this paper with the remarkable property : $e(S)-t(S)\geqslant 1$ for every numerical semigroup $S$ of this type $\frac{<a,b>}{d}$. This inequlity implies Wilf's property for these numerical semigroups, as mentioned in \ref{subsec:intro}.

 \subsection{some notations and definitions}

$\bullet$ For any real $x$, we denote respectively $\lfloor x \rfloor, \lceil x \rceil$ and $\{ x \}$ the floor, ceil and fractional part of $x$.\\

$\bullet$ For a vector $u$ of $\Z^2$, we will denote $x(u)$ and $y(u)$ its coordinates.\\  
\esp We will denote $\leqslant$ the natural partial product order on $\Z^2$ : $u\leqslant v$ if and only if $ (x(u)\leqslant x(v)$ and $ y(u)\leqslant y(v))$.\\   
\esp Then, for a subset $U$ of $\Z^2$, $\min(U)$ and $\max(U)$ are respectively the sets of \emph{minimal} and \emph{maximal} elements of $U$, for this product order on $\Z^2$.\\

$\bullet$ If $M$ is a submonoid of $(\R^n,+)$, we will denote :
$$ M^0 = \{ u\in M, -u\in M\} \esp ; \esp M^*=M\backslash M^0$$
$$ \Irr(M)=\{ u\in M^*, \forall v,v'\in M^*, u\not = v+v' \}$$

\begin{monlem} \label{lem:1}
let $\psi:\Z^p  \to \R^n$  be an additive map ( $n$ and $p$ are positive integers) and   $M$ be a submonoid of $\psi(\Z^p)$. Then :
$$ (\psi^{-1}(M))^0=\psi^{-1}(M^0) \esp ; \esp  (\psi^{-1}(M))^*=\psi^{-1}(M^*) \esp ; \esp  \Irr(\psi^{-1}(M))=\psi^{-1}(\Irr(M))$$
\end{monlem}
\begin{mademo} first two results are well known and almost obvious. Let $u\in\Irr(\psi^{-1}(M))$. We denote $m=\psi(u)$. Suppose that $m=n+n'$, with $n,n'\in M^*$. Then, we have $v,v'\in \psi^{-1}(M), w\in \ker(\psi)$, such that $u=v+v'+w$ and $\psi(v)=n,\psi(v')=n'$. But, if we denote $v"=v'+w$, then $u=v+v"$ and $\psi(v")=n'$. So, $v,v"\in (\psi^{-1}(M))^*$ : contradiction with our hypothesis on $u$. So, $m\in \Irr(M)$. \\
\esp Conversely, let $m\in \Irr(M)$ and $u\in \psi^{-1}(M^*)$, such that $\psi(u)=m$. Suppose that $u=v+v'$, with $v,v'\in (\psi^{-1}(M))^*$, then $m=n+n'$, where $n=\psi(v)$ and $n'=\psi(v')$. So, $n,n'\in M^*$ : contradiction with our hypothesis on $m$. So, $u\in \Irr(\psi^{-1}(M))$.
\end{mademo} 
 
$\bullet$ A subset $U$ of $\Z^2$ is said to be \emph{$w$-invariant} if and only if $U+w=U$.

\begin{monlem} \label{lem:2} if $U$ is a $w$-invariant subset of $\Z^2$, for $w\in\Z^2$, then $\min(U),\max(U),\Z^2\backslash U$ are $w$-invariant. ( then $\max(\Z^2\backslash U)$ is also $w$-invariant)
\end{monlem}
\begin{mademo}
  Suppose that $w+U=U$. Let $u\in\min(U)$ and $u'=u+w$, then $u'\in U$, since $U$ is $w$-invariant. If there is $u"\in U$ such that $u"\leqslant u'$, then $v=u"-w\in U$ ( $U$ is $w$-invariant) and $v\leqslant u$. So that, $v=u$ since $u$ is minimal in $U$. Then $u"=u'$ and $u'$ is minimal in $U$. This proves that $w+\min(U)\subset \min(U)$. But $U$ is also $(-w)$-invariant, then $-w+\min(U)\subset \min(U)$ and finally $\min(U)$ is $w$-invariant. The other results are proven in the same way.
   \end{mademo}



   \section{Representations in $\Z^2$}
   
   \label{sec:representation}
   
   $\bullet$ one of the simple ideas of this paper is to use Lemma \ref{lem:1} in the particular case, when $M$ is a numerical semigroup ( so $n=1$) and $p=2$ : so we have a numerical semigroup $S$ and a \emph{representation} $\psi$ of $S$ in $\Z^2$, that is a  map $\psi : \Z^2\to \R$, such that $S\subset \psi(\Z^2)$.\\
   \esp Hence, if we denote $M=\psi^{-1}(S)$, we just have to study $\Irr(M)$ to easily deduce results on $\Irr(S)$. For that aim, we will consider minimal elements of $M^*$, for the product order on $\Z^2$. So, we need to have simple relations between $\Irr(M)$ and $\min(M^*)$.\\
   
   \textbf{Question : } for what kind of submonoids $M$ of $(\Z^2,+,\leqslant)$, do we have : $\Irr(M)=\min(M^*)$ ? \\
   
   \esp We can give an elementary exemple : $\Irr(\N^2)=\min((\N^2)^*)=\{ (1,0);(0,1)\}$. \\
   \esp We will not answer the question, but only give a partial result, that will be sufficient for our purpose.

   \begin{monlem} \label{lem:3} let $L$ be a sublattice of $\Z^2$ and $w\in L$, such that $w,-w\not\in \N^2$.\\
   \esp Then, for $M=(L\cap \N^2)+w\Z$, we have :
   $$ M^0=w\Z \esp ; \esp \Irr(M)=\min(M^*)$$
   \end{monlem}
  
  \begin{mademo} 
  \esp First, $w\Z\subset M^0$. On the other hand, if $m\in M^0$, then $m=m'+k'w$ and $-m=m"+k"w$, where $m',m"\in L\cap\N^2$ and $k',k"\in\Z$. Then, $m'+m"=jw$, where $j\in \Z$. So, $jw \geqslant 0$ : contradiction with our hypothesis on $w$, if $j\not = 0$. Hence, $j=0$ and $m'+m"=0$. We deduce that $m'=m"=0$, for $m',m"\in\N^2$. So, $m\in w\Z$.\\
  \esp We remark that $M\subset L$, since $w\in L$.\\
  \esp  Let $u\in  \Irr(M)$, then $u\in M^*$. We suppose that there exists  $v\in M^*$ such that $u>v$. Then $v'=u-v >0$ and $v'\in L$, for $u,v\in L$. So, $v'\in M^*$, for $M^0=w\Z$ and $w,-w\not\in \N^2$ : no elements of $M^0$ is positive. So $u=v+v'$, with $v,v'\in M^*$ : contradiction with our hypothesis on $u$. So, $u\in \min(M^*)$.\\
  \esp Conversely : if $u\in \min(M^*)$. Suppose that $u=v+v'$ and $v,v'\in M^*$. Then, we can find $m,m'\in  L\cap (\N^2)^*$ and $k,k'\in \Z$, such that : $v=m+kw$ and $v'=m'+k'w$. Now, we denote $v"=v'+kw$ and we have $u=m+v"$, where $m>0$ and $v"\in M^*$, for $v'=v" $ ( mod $M^0$). So, $u>v"$ : contradiction with our hypothesis on $u$. So, $u\in \Irr(M)$. 
  \end{mademo}

  $\bullet$ For $a$ and $b$ two coprime positive integers, we define two maps  by :
 $$\phi : \begin{cases} \Z^2\to \Z \\ (x,y)\to ax+by \end{cases} \esp \esp ; \esp \esp   \psi :\begin{cases} \Z^2 \to (1/d) \Z \\ (x,y)\to (ax+by)/d\end{cases}$$
 
 \esp  These maps are additive group morphisms. They are  strictly increasing, surjective  and their kernel is $w\Z$, where $w=(-b,a)$.
 
\esp  Let $S=<a,b>$ be the numerical semigroup generated by $a$ and $b$ and $S'=\frac{S}{d}$. We denote $L$ the lattice defined by $L=\{ (x,y)\in\Z^2, ax+by=0 $ ( mod $d$)$\}$. We will prove in Lemma \ref{lem:5} :
$$ \phi^{-1}(S)=\N^2+w\Z \esp ; \esp \psi^{-1}(S')=L\cap (\N^2+w\Z )= (L\cap \N^2)+w\Z $$
\esp So, this representation of $\frac{<a,b>}{d}$ in $\Z^2$ is of the above type and will allow us ( see Lemma \ref{lem:1} and \ref{lem:3}) to claim that :
$$ \psi^{-1}\left(\Irr\left(\frac{<a,b>}{d}\right)\right)=\min(L\cap ((\N^2)^*+w\Z))$$ 
\esp This is detailed in this section and in \ref{subsec:irreducible}.

 \newpage

   \subsection{a lattice-stair case shape representation}

   $\bullet$ Let $S$ be a numerical semigroup and $d$ a positive integer. We have then  :
   $$ d ( S/d)=\{ s\in S, s=0 ( \text{ mod } d )\}=S\cap d\Z$$
   \esp  Let $a$ and $b$ be two coprime integers. With the two maps defined just above ( previous page), we have :
      $$ \psi^{-1}(S/d)=\phi^{-1}( d ( S/d))=\phi^{-1}(S) \cap L$$
   \esp where $L$ is the lattice :
   $$ L=\{ (x,y)\in\Z^2, ax+by =0 ( \text{ mod } d )\}=\phi^{-1}(d\Z)=\psi^{-1}(\Z)$$
   \esp and : 
  
       $$ \psi^{-1}(\N)=L\cap H_+ \esp ; \esp   \psi^{-1}(G(S/d))=(H_+\backslash \phi^{-1}(S))\cap L$$
   \esp  where $H_+$ denotes the halfplane $\{ (x,y)\in\Z^2, ax+by\geqslant 0\}$. \\

   $\bullet$ Now we suppose that $d$ is coprime with $a $ and $b$.
    With elementary arithmetic, we find that, in every square $[x,x+d[\times [y,y+d[$, $L$ has exactly one point on every line and one point on every column of $\Z^2$.\\
   \esp  $L$ is generated by $w=(b,-a), v_2=(0,d)$ and $v_1=(d,0)$ ( but it is not a basis)  : indeed these  three vectors are in $L$ and since $a$ and $d$ are coprime, for every $y$ in $\Z$, we can find $x\in\Z$ such that $(x,y)\in w\Z+v_2\Z$. So with the former remark, we obtain, adding  multiples of $v_1$, all the points of $L$.\\ 

   $\bullet$ We will study the simplest case, namely when $S$ is generated by two coprime integers $a,b$. We begin with some remarks.
   
   \begin{monlem}  \label{lem:4}.\\
   (i) the general case for $\frac{<a,b>}{d} $ is when $a,b,d$ are pairwise coprime. Indeed, for any integer $d'$ coprime with $b$ : 
   $$ \frac{<d'a,b>}{dd'}=\frac{<a,b>}{d} $$
   (ii) if $d$ divides $f(<a,b>)=ab-a-b$, then $\frac{<a,b>}{d}$ is symmetric and $f(\frac{<a,b>}{d})=\frac{f(<a,b>)}{d}$.
   \end{monlem}
   \begin{mademo} (i) let $k$ be an integer  :
   $$ kdd' \in <d'a,b> \Leftrightarrow \exists x,y\in\N, kdd'=xd'a+yb \Leftrightarrow  \exists x,z\in\N, kd=xa+zb \Leftrightarrow kd\in <a,b>  $$
   \esp as if $d'$ divides $yb$, then it divides $y$, because $d'$ and $b$ are coprime. \\
   (ii) we denote $f=\frac{f(<a,b>}{d}$ and suppose that $f$ is an integer. Then $df\not\in <a,b>$, so $f$ is a gap of $\frac{<a,b>}{d}$. In addition, $dn\in <a,b>$ for all integer $n$ such that $n>f$, because $dn> f(<a,b>)$. So, $f=f(\frac{<a,b>}{d})$. \\
   \esp If $n$ is a gap of $\frac{<a,b>}{d}$, then $dn$ is a gap of $<a,b>$, so $df-dn\in <a,b>$, for $<a,b>$ is symmetric. Hence, $f-n\in \frac{<a,b>}{d}$. Conversely, if $s\in \frac{<a,b>}{d}$ and $s<f$, then $f-s$ is a gap of $\frac{<a,b>}{d}$. So, $\frac{<a,b>}{d}$ is symmetric.  
   \end{mademo} 
   
  \esp Thus, we will suppose that $d$ is a positive integer, coprime with $a$ and $b$. We will also suppose that $d\not = a $ and $d\not = b$. Indeed : if, more generally, $d\in <a,b>$ then $kd\in <a,b>$ for all non negative integer $k$ and then $\frac{<a,b>}{d} =\N$.\\

 \newpage
   
  $\bullet$ The representation of $S=<a,b>$ and $S'=S/d$ are simple :  $S$ is represented in $\Z^2$ by a " full stair-case shape" and $S'$ by its intersection with the lattice $L$. \\
  \esp This is proven by the following Lemma :

  \begin{monlem} \label{lem:5} let $a,b,d$ be three positive pairwise coprime integers.\\
  \esp  Let $S=<a,b>$ and $S'=\frac{<a,b>}{d}$. We denote $w=(b,-a)$. \\
  Let $T$ be the triangle, " upper-half" of the rectangle $\inter{1}{b-1}\times \inter{-(a-1)}{-1}$ : 
   $$T=\{ (x,y)\in\Z^2, x\in \inter{1}{b-1}, y\in \inter{-a+1}{-1}, ax+by>0 \} $$
  
  (i) $(\N^2+w\Z ,  T+w\Z)$ is a partition of $H_+=\{ (x,y)\in\Z^2, ax+by\geqslant 0\}$.\\
  
  (ii)
   $$ \phi^{-1}(S)= \N^2+w\Z \esp ; \esp \phi^{-1}(G(S))= T+w\Z $$
 
   (iii)
    $$ \psi^{-1}(S')=(L\cap \N^2)+w\Z \esp ; \esp \psi^{-1}(G(S'))= ( L\cap T )+w\Z$$
   \esp  $\psi$ is bijective from $L\cap T$ to $G(S')$.
    \end{monlem}
   \begin{mademo}
   (i) if $u\in (\N^2+w\Z) \cap  (T+w\Z)$, then we have $u=u_1+k_1w=u_2+k_2w$, where $u_1\in T,u_2\in\N^2,k_1,k_2\in\Z$. So $u_1-u_2=kw$, where $k\in \Z$. Note $u_1=(x_1,y_1)$ and $u_2=(x_2,y_2)$. We have $x_1\in  \inter{1}{b-1}, x_2\in\N, y_1\in  \inter{-a+1}{-1}, y_2\in\N,k\in\Z$ such that : $x_1-x_2=kb$ and $y_2-y_1=ka$. So $ka>0$, thus $k\geqslant 1$ and $x_2<0$ : contradiction... This proves that $ (\N^2+w\Z) \cap  (T+w\Z)=\emptyset$.\\
   \esp On the other hand, if $(x,y)\in H_+$, then $\frac{x}{b}\geqslant -\frac{y}{a}$. We consider two cases : \\
   \esp Case 1 : if $\lfloor x/b \rfloor = \lfloor -y/a \rfloor$ and $x/b,y/a$ are not integers, then if we denote $k$ this common integer, we have : $x=kb+x'$ and $y=-ka+y'$ , with $(x',y')\in T$. , so $(x,y)\in T+w\Z$.\\
   \esp Case 2 : if $x/b=-y/a$ and are integers, then $(x,y)\in w\Z\subset(\N^2+w\Z)  $. \\
   \esp Case 3 :  if $\lfloor x/b \rfloor > \lfloor -y/a \rfloor$, then if we denote $k=\lfloor x/b \rfloor$, we have $x=kb+x'$ and $y=-ka+y'$,  with $(x',y')\in \N^2$, so $(x,y)\in \N^2+w\Z$. \\
   (ii)  We have $\phi(\N^2)=S$ and use (i) for $G(S)$.\\
   (iii) with the general remarks at the beginning of the section and the fact that $L$ is $w$-invariant :
   $$ \psi^{-1}(S')=L\cap (\N^2+w\Z)=(L\cap \N^2)+w\Z$$
   \esp the same arguments rules for $G(S')$. \\
   \esp Now, $\psi$ is bijective from $T\cap L$ to $G(S')$, since for every $u\in T, k\in \Z^*, u+kw\not\in T$.
   \end{mademo}

 $\bullet$  For now in this section, we adopt the following hypothesis and notations : \\
 - we suppose that  $a,b,d$ are pairwise coprime. \\
- we denote $L'=L\backslash \{ (0,0)\}$ and  $S'=\frac{<a,b>}{d}$.\\ 
 - the symbol $A\overset{\psi} {\simeq}B$ means that $A$ and $B$ are in bijective correspondence via $\psi$.\\
 - we use the usual product order on $\Z^2$ ( partial order).
 
 \begin{monlem} \label{lem:6}.\\
  (i)   
 $$S'\backslash \{ 0 \} \overset{\psi} {\simeq} (\N\times \inter{0}{a-1})\cap L'$$
 (ii) 
 $$ \Irr(S') \subset \psi(\min((\inter{0}{d}\times \inter{0}{a-1})\cap L')$$
 (iii) 
 $$ d<b \Rightarrow  \Irr(S')\overset{\psi} {\simeq} \min((\inter{0}{d}\times \inter{0}{a-1})\cap L')$$
 \end{monlem}
 
 \textbf{Remark : } $\Irr(S')$ is the minimal set of $S'\backslash \{ 0 \} $ for the order $\leqslant_S$ on $\Z$. The result (i) would suggest that, since $\psi$ is increasing, then $\Irr(S')$ is via $\psi$ in one to one correspondence with $\min(E)$, where $E= (\N\times \inter{0}{a-1})\cap L'$..., but this is not true in general, because $\psi$ is not an order isomorphism on these sets.  We can have $u,u'\in E$, with $\psi(u)\leqslant_{S'} \psi(u')$ and $u\not \leqslant u'$ : it is the case if $u'\geqslant u+w$ and $y(u')<y(u)$. \\
 
 \begin{mademo}
--- (i) $\psi$ is injective on  $(\N\times \inter{0}{a-1})\cap L'$, because for two different elements $u,u'$ of this set, $u'-u \not\in w\Z$ ( as $| x(u'-u) |<a$). Moreover,  let $u\in (\N^2\cap L')+w\Z$, then  : 
  $$u=kw+u' \text{ with } k=-\left\lfloor \frac{y(u)}{a}\right\rfloor \text{ and } y(u')\in \inter{0}{a-1}$$
  \esp in addition, $u=jw+u"$ for some integer $j$ and $u"\in \N^2$, so : $u'=u"+(j-k)w$. We deduce $y(u')\geqslant (k-j)a$, so $k-j\leqslant 0$, and $x(u')\geqslant (j-k)b\geqslant 0$. Then, $u'\in (\N\times \inter{0}{a-1})\cap L$. In addition, $u\not\in w\Z$, so $u"\not\in w\Z$, hence $u'\not = (0,0)$. \\
 --- (ii)   Let $s\in \Irr(S')$ and $u\in (\N\times \inter{0}{a-1})\cap L'$ such that $s=\psi(u)$ ( see (i)). If $x(u)>d$, then $u=(d,0)+u'$ with $u'\in (\N\times \inter{0}{a-1})\cap L'$, so $s=a+s'$, with $s'=\psi(u')\in S'$ and $s'\not = 0$, which is in contradiction with $s \in \Irr(S')$, since $a\in S'$. We have proven that $x(u)\leqslant d$, and so that $u\in (\inter{0}{d}\times \inter{0}{a-1})\cap L'$.\\
   \esp Now, we want to prove that $u$ is minimal in $M=(\inter{0}{d}\times \inter{0}{a-1})\cap L'$. If it is not, then we have $u'\in M$, such that $u'<u$. Let $u"=u-u'$, then $u"\in L'$ and $0\leqslant x(u')\leqslant x(u) \leqslant d$, so $x(u")\in \inter{0}{d}$, and with similar arguments, $y(u")\in\inter{0}{a-1}$. So $u"\in M$ and $s=s'+s"$, with $s'=\psi(u'),s"=\psi(u")\in S'\backslash \{ 0 \}$. It is impossible since $s\in \Irr(S')$. \\
 ---  (iii) Reversely : let $u\in \min(M)$ and $s=\psi(u)$. Then, $s\in S'$ and $s\not = 0$. We would like to prove that $s$ is irreducible in $S'$ : if it was not, then we would have $s_1$ and $s_2$ in $S'$ such that $s=s_1+s_2$ and $s_1,s_2$ non null. Let $u_1$ and $u_2$ be the unique elements in $(\N\times \inter{0}{a-1})\cap L'$ ( see (i)), such that $s_1=\psi(u_1)$ and $s_2=\psi(u_2)$. Then, $\psi(u)=\psi(u_1+u_2)$ and $u_1+u_2\in (\N\times \inter{0}{2a-2})\cap L'$. So, $u=u_1+u_2$ or $u=u_1+u_2+w$. In the last case, we have $x(u)\geqslant b$, so that, if $d<b$, we obtain a contradiction. In that case $d<b$, we can conclude that $u=u_1+u_2$ and then $u_1,u_2\in M$, since they have non negative $x$. Then, $u$ is not minimal in this set, as $u_1>0$ and $u_2>0$ : contradiction.   So $\psi(S')=\min(M)$ and we conclude with (i). 
 \end{mademo}
    

   \subsection{relations between invariants of $\frac{<a,b>}{d}$ and a lattice in a rectangle}

   \label{subsec:relations}

   $\bullet$ We begin with the genus $g(S')$. At the beginning of the previous subsection, we have seen in Lemma 5 that $T\cap L$ and $G(S')$ ( see notations in this Lemma) are in   bijective correspondence via $\psi$. We will use another property with the following Proposition :

 \begin{maprop} \label{prop:1} let $R_0$ be the rectangle : $\inter{1}{b-1}\times \inter{-(a-1)}{-1}$ and $L$ the lattice $L=\{ (x,y)\in \Z^2, ax+by=0[d] \}$.
   $$ g\left(\frac{<a,b>}{d}\right)= \frac{1}{2}\#(L\cap R_0) $$
  
   \end{maprop} 
 
 \begin{mademo}
  The above remark proves that $g(S')= \# (T\cap L)$. In addition, the symetry of $\Z^2$ , $\sigma: (x,y)\to (b-x,-a-y)$ carries $T$ onto $T'$ such that $(T,T')$ is a partition of $R_0$.\\
  \esp Indeed, $T=\{ (x,y)\in R_0, ax+by>0 \}$, so $T'=\{ (x,y)\in R_0, ax+by<0 \}$. Moreover, there is no point $(x,y)$ in $R_0$ such that $ax+by=0$, since $a$ and $b$ are coprime. So, the result is proven. 
  \end{mademo}

  \newpage

  $\bullet$ Next proposition is our main result to express minimal generators of $\frac{<a,b>}{d}$ : 
  
   \begin{maprop} \label{prop:2} let $L$ be the lattice $L=\{ (x,y)\in \Z^2, ax+by=0[d] \}$.\\ 
   $\blacktriangleright$ Case 1 : if $d<a<b$, then $a,b\in \Irr\left(\frac{<a,b>}{d}\right) $ and :
   $$  \Irr\left(\frac{<a,b>}{d}\right)\backslash\{a,b\}\overset{\psi} {\simeq}\min(L \cap \inter{1}{d-1}^2)$$
    $\blacktriangleright$ Case 2 : if $a<d<b$, then $a\in \Irr\left(\frac{<a,b>}{d}\right),b\not\in \Irr\left(\frac{<a,b>}{d}\right)$ and :
    $$ \Irr\left(\frac{<a,b>}{d}\right)\backslash\{a\}\overset{\psi} {\simeq} \min(L \cap (\inter{1}{d-1}\times\inter{1}{a-1}))$$
     $\blacktriangleright$ Case 3 : if $a<b<d$, then :
    $$ \Irr\left(\frac{<a,b>}{d}\right)\overset{\psi} {\simeq} \min(L\cap \inter{1}{x_1}\times \inter{0}{a-1})$$
    where $x_0=\min\{ x(u), u\in (\N\times \inter{1}{a-1})\cap L \}$ and $x_1=\min(d,x_0+b-1)$.\\
    so $a\in\Irr\left(\frac{<a,b>}{d}\right) \Leftrightarrow x_0>d-b$. 
     \end{maprop}

\textbf{Remark  : } if $d>\max(a,b)$ and $x_0+b>d$, then the result of Case 3 is the same as for Case 2.\\
\esp We can summarize all these cases as follows : there exists two positive integers $c$ and $c'$ such that :  $\Irr\left(\frac{<a,b>}{d}\right)\overset{\psi} {\simeq} \min(L'\cap \inter{0}{c}\times \inter{0}{c'} )$, where $L'= L\backslash \{ (0,0)\}$.\\
\esp For case 1, $c=c'=d$, for case 2, $c=d$ and $c'=a-1$, for case 3, $c=x_1$ and $c'=a-1$.\\

  \begin{mademo} We use Lemma \ref{lem:6} (iii). \\
 $\blacktriangleright$ Case 1 : if $d<a<b$ : then, $(d,0)$ and $(0,d)\in  \min(L' \cap (\inter{0}{d}\times\inter{0}{a-1}))$ ( remind that $L$ has exactly one point on every $\Z$-line or row of length $d$), so $a,b\in \Irr(S')$ and we have also  $ \Irr(S')\backslash\{a,b\}\overset{\psi} {\simeq}\min(L \cap \inter{1}{d-1}^2)$, since $(0,d)$ and $(d,0)$ are not comparable with all elements of this set. Moreover, if $u\in \inter{0}{d}\times \inter{d}{a-1}$, then $(0,d)<u$. \\
   $\blacktriangleright$ Case 2 : if $a<d<b$ : then $(d,0)\in \min(L' \cap (\inter{0}{d}\times\inter{0}{a-1}))$, but $(0,d)$ is not in this set ! We conclude in the same way as above. \\
  $\blacktriangleright$ Case 3 : if $a<b<d$ : this case is more intricate.\\
--- we have seen in  Lemma \ref{lem:6} (ii), that : if $s\in \Irr(S')$, then $s=\psi(u)$, with $u\in  \min((\inter{0}{d}\times \inter{0}{a-1})\cap L')$.\\
 \esp  We have $x(u)\geqslant x_0$, with the definition of $x_0$, for if $y(u)=0$, then $u=(d,0)$ and $x_0<d$ ( see property of $L$ again). Let us prove that "$x(u)\leqslant x_0+b-1$ ". If $x(u)\geqslant x_0+b$, we denote $u'=u-w$, so $u_0\leqslant u'$, where $u_0=(x_0,y_0)\in L' \cap (\inter{0}{d}\times\inter{0}{a-1})$. So, if we denote $s_0=\psi(u_0)$, we obtain $s_0\leqslant_{S'} s=\psi(u')$. But, $s_0\in S'$ and $s_0\not = 0$ : contradiction with $s\in \Irr(S')$ ! \\
 \esp So, $u\in (\inter{x_0}{x_1}\times \inter{0}{a-1})\cap L)$ and is minimal in it, since it is minimal in a bigger set. \\
 --- on the other hand, suppose that $u\in \min(    (\inter{x_0}{x_1}\times \inter{0}{a-1})\cap L)$ and $s=\psi(u)$, then $s\in S'$ and $s\not=0$. Suppose that  $s=s_1+s_2$, with $s_1$ and $s_2$ in $S'$ and $s_1,s_2$ non null. Let $u_1$ and $u_2$ be the unique elements in $(\N\times \inter{0}{a-1})\cap L'$, such that $s_1=\psi(u_1)$ and $s_2=\psi(u_2)$. Then, $\psi(u)=\psi(u_1+u_2)$ and $u_1+u_2\in (\N\times \inter{x_0}{2a-2})\cap L'$. So, $u=u_1+u_2$ or $u=u_1+u_2+w$. In the last case, we have $x(u)\geqslant b+x_0>x_1$ : contradiction. So, $u=u_1+u_2$ and $u$ is not minimal in $(\inter{x_0}{x_1}\times \inter{0}{a-1})\cap L)$, since $u_1,u_2$ are in this set and are not null ( $x_0\leqslant x(u_1),x(u_2)$ and $x(u_1)+x(u_2)=x(u)\leqslant x_1$). \\
 \esp We have proved that : if $u\in \min(    (\inter{x_0}{x_1}\times \inter{0}{a-1})\cap L)$, then $\psi(u)\in \Irr(S')$, so the result is obtained with Lemma \ref{lem:6} (i), since $\psi$ is bijective over this set. \\
 \esp To finish : $(0,d)\not\in \inter{x_0}{x_1}\times \inter{0}{a-1})\cap L)$ and $(d,0)\in \inter{x_0}{x_1}\times \inter{0}{a-1})\cap L$ if and only if $x_1\geqslant d$. In that case, it is minimal in this set !    
  \end{mademo}

\newpage

  $\bullet$ Now, in the same way as for $\Irr(S')$, the set $\PF(S')$ can be related with maximal points of a certain rectangle of $L$ : 
  
  \begin{monlem}\label{lem:7}  we denote $R$ the rectangle of $\Z^2$ defined by : $ R= \inter{b-d}{b-1}\times \inter{-d}{-1} $. We denote the numerical semigroup $S=\frac{<a,b>}{d}$, the lattice $L=\{ (x,y)\in \Z^2, ax+by=0[d] \}$ and the triangle $T=\{ (x,y)\in\Z^2, x\in \inter{1}{b-1}, y\in \inter{-a+1}{-1}, ax+by>0 \} $. Then : 
  $$  \PF(S)\overset{\psi} {\simeq}\max(L\cap R)\cap T  $$
  We can precise if we distinguish several cases : \\
   $\blacktriangleright$ Case 1 : if $d<a<b$ :
  $$ \PF(S) \overset{\psi} {\simeq} \max(L\cap R) $$
   $\blacktriangleright$ Case 2 : if $a<d<b$ :
  $$   \PF(S) \overset{\psi} {\simeq} \max(L\cap \inter{b-d}{b-1}\times \inter{-(a-1)}{-1})\cap T $$
   $\blacktriangleright$ Case 3 : if $a<b<d$ :
  $$   \PF(S) \overset{\psi} {\simeq} \max(L\cap \inter{1}{b-1}\times \inter{-(a-1)}{-1})\cap T $$
  
\end{monlem}

\textbf{Remark : } this result will be improved in next subsection, after we have proved that these maximal points are in $T$ ( see Proposition 4). But now, we can use the following argument :  let $C$ be a  square of $\Z^2$, say $C=\inter{x_0}{x_0+c-1}\times \inter{y_0,y_0+c-1}$, and $E$ a    subset of $C$,  that contains exactly one point on every row and  one point on every line. Then, minimal points of $E$ are all under ( or on) the diagonal
$x+y=x_0+y_0+c-1$ : these points $(x,y)$ verify : $x+y\leqslant x_0+y_0+c-1$. The proof is obvious, with reducio ad absurdum for instance. \\

 \begin{mademo}
 \esp We will use similar arguments as in proof of the previous result. As a reminder : $G(S')\overset{\psi} {\simeq}L\cap T$ ( see Lemma \ref{lem:5} (iii)). \\
 --- Let $s\in \PF(S')$, then $s\in G(S')$, so we have $u\in L\cap T$ such that $s=\psi(u)$.\\
 \esp  We claim :  $u\in L\cap R$ : indeed, if $u\not\in L\cap R$ ( only possible if $d<\max(a,b)$), then $u+(d,0)$ or $u+(0,d)\in L\cap T$ and so $s+a$ or $s+b\in G(S')$, which contradicts our hypothesis on $s$. Suppose that $u$ is not maximal in $L\cap R$, then we have $u'\in L\cap R$ such that $u<u'$. Let $u"=u'-u$, then $u"\in \N^2\cap L'$ so $s"=\psi(u")\in S'$. But, $s"\not = 0$, since $u">0$. Now $s+s"=\psi(u')\not\in S'$ :  contradiction with our hypothesis on $s$.\\
--- Conversely, let $u\in \max(L\cap R)\cap T$ and $s=\psi(u)$. Then $u\in L\cap T$, so $s\in G(S')$. Suppose that $s\not\in \PF(S')$, then we have $s'\in S'$, such that $s'\not = 0$ and $s"=s+s'\not\in S'$. Let $u'\in L'\cap \N^2$, such that $\psi(u')=s'$ and $u"\in L\cap T$, such that $\psi(u")=s"$ ( indeed $s"\in G(S')$, for $s">0$). Then, $\psi(u+u')=\psi(u")$, so $u+u'=u"+kw$, for some integer $k$. 
With first coordinate in $\Z^2$, we obtain :
 $$1\leqslant x(u+u')=x(u"+kw)<(k+1)b$$
 \esp  so $k\geqslant 0$. With second coordinate, we obtain : 
 $$-a<y(u+u')=y(u"+kw)< -ak$$
 \esp so $k\leqslant 0$. We deduce : $k=0$ and so $u<u"$, since $u'>0$ : it is in contradiction with our hypothesis.\\
 \esp We have proved the first assertion. \\

  $\blacktriangleright$ Case 1 : if $d<a<b$ : with the remark above, $\max(R\cap L)$ is a subset of the upper-half triangle of the square $R\cap L$. Then  this triangle is a subset of $T$ and the result is proven.\\
  $\blacktriangleright$ Case 2 : if $a<d<b$ : let denote $R'=\inter{b-d}{b-1}\times \inter{-(a-1)}{-1})$. Then, $\max(R\cap L)\cap T\subset \max(R'\cap L)\cap T $, for $R\cap T\subset R'\subset R$. On the other hand, no element of $R\backslash R'$ is greater than any element of $R'$ ( look at $y$), so the result is proven.\\
 $\blacktriangleright$ Case 3 : if $a<b<d$ :  same arguments as in Case 2 with $R"=  \inter{1}{b-1}\times \inter{-(a-1)}{-1})$.  \end{mademo}



\newpage

    \subsection{minimal points of a lattice in a rectangle}

    \label{subsec:minimal}

    $\bullet$ Proposition \ref{prop:2} and Lemma \ref{lem:5} are about minimal points of a lattice in a rectangle ( in $\R^2$ with the product order). Indeed, if we use an "half turn  symmetry" in $\R^2$, such as $u\to v-u$, then maximal points are turned into minimal points and vice versa...

   \begin{monlem}   \label{lem:8} Let $L$ be a lattice of $\R^2$ and $R=\inter{x_0}{x_0+c}\times \inter{y_0}{y_0+c'}\backslash \{ (x_0,y_0)\}$, where $x_0,y_0$ are two reals and $c,c'$ are two positive reals. \\
   \esp if $(u_k)_{k\in \inter{1}{e}}$ is  $x$-monotonous or $y$-monotonous and describes $\min(L\cap R)$, then : \\
   (i) both $(x(u_k))_k$ and $(y(u_k))_k$ are monotonous : one is increasing and the other is decreasing.\\
   (ii)
   $$ \forall i\in\inter{2}{e-1}, u_{i+1}+u_{i-1} \geqslant 2u_i$$
   
   that is to say : $(u_i-u_{i-1})_i$ is an increasing sequence in $\Z^2$ ( with the product order). \\
   (iii) we have equality in (ii) for a fixed $i\in \inter{2}{e-1}$ if and only if $2u_i-u_{i-1} \in R$.
   \end{monlem}

   
   \textbf{Remark 1 : } we could name \emph{ order-convex} the property (ii).\\
    \textbf{Remark 2 : }  the curve connecting the  $(u_k)_k$ is the " under border" of the convex hull of $L\cap R$ in $\R^2$.  \\ 
    \textbf{Remark 3 : } (iii) means that $(u_i-u_{i-1})_i$ is constant as long as possible :  the curve connecting the  $(u_k)_k$ remains straight as long as it does not get out of the rectangle.\\
    
   \begin{mademo}
   (i) if only one of $(x(u_k))_k$ and $(y(u_k))_k$ is monotonous, then we can find $i$ such that $u_i$ and $u_{i+1}$ are comparable, which is not possible ! \\
   (ii) we remark that the conclusion is unchanged if we reverse the indexation of the sequence, so that we can suppose that $(x(u_k))_k$ is decreasing. Thus, $(y(u_k))_k$ is increasing. \\
  \esp Let $i\in\inter{2}{e-1}$ and we denote $v_i=u_{i+1}+u_{i-1}-u_i$. Then, $v_i\in L\cap R$, for :
  $$x(u_{i+1})< x(v_i)< x(u_{i-1}) \esp ; \esp y(u_{i-1})<y(v_i)<y(u_{i+1})$$
  \esp So, there exists $k\in\inter{1}{e}$, such that $u_k\leqslant v_i$. But, the previous inequalities show that $k$ can not be different from $i$, so $u_i\leqslant v_i$...q.e.d. \\
  (iii) if we have equality, then $2u_i-u_{i-1}=u_{i+1}\in R$. Conversely : if  $2u_i-u_{i-1}\in R$, then if we denote $t_i$ that element, we have $u_{i+1}\geqslant t_i$ ( see proof of (ii)) and $t_i\in L\cap R$. But $u_{i+1}$ is minimal in $L\cap R$, so $u_{i+1}=t_i$. 
   \end{mademo}

  $\bullet$  Now, we are able to propose a statement that improves Lemma \ref{lem:7} :
   \begin{maprop} \label{prop:3} $a,b,d$ are three pairwise coprime positive integers.\\
   \esp   We denote $\psi: \Z^2\to (1/d)\Z ; (x,y)\to (ax+by)/d$ and $L=\{ (x,y), ax+by=0[d] \}=\psi^{-1}(\Z)$. \\
   \esp The symbol $\overset{\psi} {\simeq}$ means that both sets are in bijective correspondance via $\psi$. 
    $$ \PF\left(\frac{<a,b>}{d}\right) \overset{\psi} {\simeq} \max(L\cap \inter{b-\delta}{b-1}\times \inter{-\delta'}{-1}) $$
   \esp where $\delta=\min(d,b-1)$ and $\delta'=\min(d,a-1)$.
  
\end{maprop}
 \begin{mademo} we denote $R_1= \inter{b-\delta}{b-1}\times \inter{-\delta'}{-1})$. If we compare this statement with Lemma \ref{lem:7}, we have just to show that $\max(L\cap R_1) \subset T$, with notations of Lemma \ref{lem:7}. This has been done in the case $d<\min(a,b)$. \\
 \esp If $a<d<b$, then $R_1=\inter{b-d}{b-1}\times \inter{-(a-1)}{-1}$. We denote $R_0=\inter{b-d}{b}\times \inter{-a}{-1}$ and we claim that $\max(L\cap R_1)=\max(L\cap R_0)\backslash\{ (b,-a)\}$. Indeed, $(b,-a)$ is a point of $L$ and it is the only one in $R_0$ on the row $y=b$. So, $(b,-a)\in \max(L\cap R_0)$ and $\max(L\cap R_1)\subset\max(L\cap R_0)$, since there is no point in $L\cap (R_0\backslash R_1)$ that is greater than a point of $R_1$. Our equality is proven. By now, our former Lemma proves ( via an half-turn ) that $\max(L\cap R_0)$ can be enumerated by an "order-concave" sequence of points : $\max(L\cap R_0)= \{ u_k , k\in\inter{0}{t}\}$, such that $u_0=(b,-a)$ and $(u_k-u_{k-1})_k$ is decreasing. We deduce that : $\max(L\cap R_1)=\{ u_k , k\in\inter{1}{t}\}\subset T$, for $u_t=(x_1,-1)$ the unique point of $L$ that lies on the $\Z^2$-segment of length $d$ : $x\in \inter{b-d}{b-1}, y=-1$. \\
 \esp If $d>\max(a,b)$, then $R_1=\inter{1}{b-1}\times \inter{-(a-1)}{-1}$. We use similar arguments as in previous case with $R_0=\inter{0}{b-1}\times \inter{-a}{0}$, as we remark that $(0,0)$ and $(b,-a)$ are points of $L$, that are maximal points of $L\cap R_0$ ( no point of $L$ on the $\Z^2$-segment $x\in \inter{1}{b},y=0$). Moreover, no point of $R_1$ is greater than one of these two points.
\end{mademo}

   $\bullet$ In order to study the minimal generators of $\frac{<a,b>}{d}$, we consider a particular case of Lemma \ref{lem:8} :

   \begin{monlem} \label{lem:9}.\\
    let $L$ be a lattice of $\R^2$ and $R=\inter{0}{c}\times \inter{0}{c'}\backslash \{ (0,0) \}$, where $c,c'$ are positive reals.\\
   \esp We can describe $\min(L\cap R)$ by a sequence $(u_i)_{i\in\inter{1}{e}}$, such that  \esp $(u_1,u_2)$ is a $\Z$-basis of $L$ ( if $e\geqslant 2$),  $(x(u_i))_i,(y(u_i))_i$ are monotonous and :
   $$ \forall i\in\inter{2}{e-1},\exists k_i\in\inter{2}{+\infty}, u_{i+1}+u_{i-1}=k_iu_i \esp (1) $$
  \esp moreover :
   $$ \forall i\in\inter{2}{e-1},k_i= \min(n\in\N, nu_i\geqslant u_{i-1})  \esp (2) $$
   
   \end{monlem}   
    \begin{mademo}
      Let us consider $(u_k)_{k\in \inter{1}{e}}=\min(L\cap R)$, such that $(x(u_k))_k$ is decreasing and $e=\#(\min(L\cap R))$. Then, $(y(u_k))_k$ is increasing ( see arguments of the proof of the previous Lemma). \\
     \esp If $e\geqslant 2$, then $(u_1,u_2)$ is a $\Z$-basis of $L$, because there is no point of $L$ in the interior  ( in $\R^2$, or $\Q^2$) of the triangle based on $(0,u_1,u_2)$, hence by symmetry $u\to u_1+u_2-u$, there are no points in the interior of the parallelogram based on $(0,u_1,u_2,u_1+u_2)$.\\
   \esp  Now, let us prove the induction formula. Let $i\in\inter{2}{e-1}$. Let $j$ be the grestest integer such that $u_{i+1}+u_{i-1}\geqslant ju_i$ and $\rho=u_{i+1}+u_{i-1}- ju_i$. Then, $\rho\geqslant (0,0)$ and $j\geqslant 2$, according to the previous Lemma. If $\rho=(0,0)$, our result is proven. Else, we have $x(\rho)<x(u_i)$ or $y(\rho)<y(u_i)$. Suppose that $x(\rho)<x(u_i)$. We also have :
    $$ y(\rho)\leqslant y(u_{i+1})+y(u_{i-1})-2y(u_i)< y(u_{i+1})-y(u_i)<y(u_{i+1})$$
    \esp Hence, $\rho\in L\cap R$ and it exists $k\in\inter{1}{e}$, such that $u_k\leqslant \rho$. But $x(\rho)<x(u_k)$ if $k\leqslant i$ and $y(\rho)<y(u_k)$ if $k>i$, so this is impossible. We use similar arguments if $y(\rho)<y(u_i)$, proving that $x(\rho)<x(u_{i-1})$. \\
    \esp Conclusion : $\rho=(0,0)$ and our formula is proven.\\
    \esp We end this proof with the expression of $k_i$ as a minimum. Let $i\in \inter{1}{e-1}$. As $u_{i+1}\geqslant 0$, then $k_iu_i\geqslant u_{i-1}$. If there exists a positive integer $n<k_i$ such that $nu_i\geqslant u_{i-1}$, then : we denote $w=nu_i-u_{i-1}$ , so $u_{i+1}=w+(k_i-n)u_i>w$. So, $0\leqslant w < u_{i+1}$ and $w\in L$. In addition, $w\not =0$, since $u_{i-1}\not = nu_i$ ( this would imply $u_i\leqslant u_{i-1}$). Hence, $w\in R\cap L$ and $w< u_{i+1}$ : contradiction.
  \end{mademo}
  
  \textbf{Remark : } with hypothesis of previous Lemma, we can also conclude that $(u_i,u_{i+1})$ is a $\Z$-basis of $L$ if $i\in\inter{1}{e-1}$ ( obvious induction with formula (1)).\\
    
\esp This has a direct consequence on $\Irr\left(\frac{<a,b>}{d}\right)$, when we use Proposition \ref{prop:2} : \\
\esp Let $a,b,d$ be three pairwise coprime positive integers, then  we can find a finite sequence of positive integers $(m_k)_{k\in \inter{1}{e}}$, such that $m_1$ and $m_2$ are coprime, $\Irr\left(\frac{<a,b>}{d}\right)=\{ m_k, k\in\inter{1}{e}\}$ and :
  $$ \forall i\in\inter{2}{e-1}, \exists k_i\in\inter{2}{+\infty},m_{i+1}+m_{i-1}=k_im_i$$
    

 \esp As proven in \textbf{[9]}, the converse is also true, if in addition we suppose that $m_i\not\in <m_k,k\not = i >$ for all $i\in\inter{1}{e}$ ! We will explore that in next section.

  \newpage

  \section{The reverse problem}
  
  \label{sec:reverse}

 $\bullet$ Let us consider the set $A$ of triplets of positive integers $(a,b,d)$ such that $1<a<b$ and  $a,b,d$ are pairwise coprime. We also consider the following map : 
 $$\mathcal{I} :  \begin{cases} A \to \P(\N) \\
 (a,b,d) \to \Irr\left(\frac{<a,b>}{d}\right) \end{cases}$$

  \esp In \textbf{[9]}, we can find a precise description of the range of $\mathcal{I} $ : it is the set $B$ of finite subset of $\N$,  that  can be enumerated by a sequence $(n_k)_{k\in\inter{1}{e}}$ , such that : 
  $$ \begin{cases}  (i) \esp \forall k\in\inter{1}{e},n_k\not\in < n_i,i\not = k > \\
  (ii) \esp \gcd(n_1,n_2)=1  \\
  (iii)\esp \forall i\in\inter{2}{e-1}, \exists c_i\in\N, n_{i-1}+n_{i+1}-2n_i=c_in_i \end{cases} $$
  \esp We have proved half of this result at the end of previous section. We will prove the other half, giving the preimage of each element of $B$. \\ 
  
  
  \textbf{Remark 1 : } we mention that $< P >$ is the $\N$-span of $P$ in a monoid, so (i) could be named " $\N$-independance of the $(n_i)_i$" and is obvious for a minimal set of generators of a semigroup. This implies that the $(n_k)_k$ are pairwise distincts.\\
  \esp  Property (ii) could be generalized ( as a consequence of (iii)) : 
  $$ \forall i\in\inter{1}{e-1}, \gcd(n_i,n_{i+1})=1 $$
  \esp Property (iii) implies that $(n_i)_i$ is convex and that, in a sense, "the second derivative of $(n_i)_i$" is a ( integral) multiple of $(n_i)_i$ : if we denote $\underline{\delta}(n)=(n_i-n_{i-1})_i$ and $\overline{\delta}(n)=(n_{i+1}-n_i)_i$ then
  $$\underline{\delta}\overline{\delta}(n) =\overline{\delta}\underline{\delta}(n)=c\times n $$
  \esp where $c$ is a $\N$-valued sequence. \\
 
 \textbf{Remark 2 : } If the sequence $(n_k)_{k\in\inter{1}{e}}$ satisfies conditions (i),(ii) and (iii), then, for every $k_1,k_2\in\inter{1}{e}$, such that $k_1\leqslant k_2$, the sequence $(n_k)_{k\in\inter{k_1}{k_2}}$ satisfies (i),(ii) and (iii).\\

\textbf{Question : } given such a finite subset $I$ of $\N$, how many indexations of its elements do exist such that the above conditions (i),(ii),(iii) are satisfied ? \\
\esp If $I$ contains more that one element, then it exists as least two indexations : $(n_k)_{k\in\inter{1}{e}}$ and $(n_{e+1-k})_{k\in\inter{1}{e}}$.  In some cases, there exist other indexations : for example, if $I=\{ 7,11,59 \}$, we have 4 indexations of $I$ that satisfy (i),(ii) and (iii) :
$$  (59,11,7) \esp ; \esp (7,11,59)  \esp ; \esp (11,7,59)  \esp ; \esp (59,7,11)$$  
 \esp As a consequence, we can find 2 triplets $(a,b,d)$ of pairwise coprime integers such that $d<a<b$ and $\Irr(\frac{<a,b>}{d})= I$ ( see Thoerem \ref{theo:1} in 3.4) :
 $$ \Irr\left(\frac{<7,59>}{6}\right)= \{ 7,11,59 \} = \Irr\left(\frac{<11,59>}{10}\right)$$

\newpage

     \subsection{modular-convex sequences}
    
      \label{subsec:modular}
      
    $\bullet$ Let us denote $\Delta$ the operator of discrete second derivative of $\R^n$-valued sequences ( $n\in\N^*$) :
      $$ \Delta : u \to v \text{ such that } \forall i, \esp v_i=u_{i-1}+u_{i+1}-2u_i $$
      \esp Then, we have $\Delta=  \overline{\delta}\underline{\delta}=\underline{\delta} \overline{\delta}$, with notations of remark 1 above.  \\
  \esp In 2.4 and Lemma \ref{lem:9}, we have proved that the minimal set of a lattice of $\R^2$ in a rectangle $R=(\inter{0}{c}\times \inter{0}{c'})\backslash \{ (0,0)\}$ can be enumerated by a sequence $u$, satisfying the relation : 
     
    $$ \Delta(u)=c\times u $$ 
    \esp where  $c$ is a sequence valued in $\N$. In other words :  
     $$ \forall i,\exists c_i\in\N, \esp  u_{i+1} +u_{i-1}=(c_i+2)u_i $$
    
   \esp We also have a $\mathrm{SL}_2(\Z)$-matricial version : 
    $$ \forall i, \exists d_i \in \inter{2}{+\infty},\esp  \begin{pmatrix} u_{i+1} \\ u_i \end{pmatrix} = \begin{pmatrix} d_i & -1 \\ 1 & 0 \end{pmatrix} \begin{pmatrix} u_{i} \\ u_{i-1} \end{pmatrix} $$

   \esp For convenience, we will name such sequences with the very unperfect term " modular-convex", for this mix of unimodular matrices and convexity-like property : 
    \begin{madef}  \label{def:1} a $\R^n$ valued sequence $(u_i)_{i\in J}$ is  \textbf{modular-convex} if and only if it exists a $\N$-valued sequence $(c_i)_{i\in\overset{\sim}{J}}$, such that $\Delta(u)=c\times u$. ( $J$ is an interval of integers and $\overset{\sim}{J}$ is $J$ without its  finite extremities )
    \end{madef}
    
\textbf{Remark 1 : } a numerical modular-convex sequence is not always convex. But, we have a simple result : if $(x_i)_{i\in\inter{0}{r}}$ is a numerical modular-convex sequence such that $x_1\geqslant \max(0,x_0)$, then it is convex. Indeed, by obvious induction, we prove that $x_i$ is non negative for $i>0$ and $(x_i-x_{i-1})_{i\geqslant 1}$ is non decreasing and non negative. \\      
    
    $\bullet$ Now, we take interest in these sequences and will show that they are fully given by three terms : first, second and last one. In addition, the coefficients $(d_i)_i$ are related with some continued fraction development :

    \begin{maprop} \label{prop:4} let $u,v,w \in \R^n$ such that $u,v$ are linearly independant ( so $n\geqslant 2$)
      $$ \exists r\in\N^*, \exists (u_i)_{i\in\inter{0}{r}} \text{ modular-convex,  such that } \begin{cases} u_0=u \\ u_1=v \\ u_r=w \end{cases} \Leftrightarrow \exists p,q \in\N, \begin{cases} 0<p<q \\ \gcd(p,q)=1 \\  w=-pu+qv \end{cases} $$ 
 \esp in that case, the sequence $(u_i)_i$ is unique and the integers $c=(c_i)_i$ such that $\Delta(u)=c\times u$ are given by the following algorithm, where we denote $d_i=c_i+2$ : 
   $$t_0=\dfrac{p}{q }, \esp  \forall i\in\inter{1}{r-1}, \esp d_i=\left\lceil \frac{1}{t_i} \right\rceil \esp ; \esp t_{i+1}=\left\{ - \frac{1}{t_i}\right\} \text{ and } t_r=0  \text{ ( stop condition )}$$
     \end{maprop} 
     
     \textbf{Remark 2 : } the above coefficients $(d_i)_{i\in\inter{1}{r-1}}$ are the coefficients of the \emph{ ceiling continued fraction expansion} of $\frac{q}{p}$. That is to say :
    $$ \frac{q}{p}= d_1-\dfrac{1}{d_2-\dfrac{1}{d_3-\cdots }}$$ 
    \esp Unlike the usual continued fraction algorithm, that is of quadratic complexity in terms of the binary size of inputs, this one is exponential in worst cases : indeed, if we take $q=n$ and $p=n-1$, then we obtain a sequence of $n-1$ digits all equal to 2, for the  ceiling continued fraction expansion of $1/n$.\\
    \esp There exists a relation between this two types of continued fraction expansion, but it is not obvious and will be explicited in \ref{subsec:CFE}.\\
    \esp However, we can easily give some relation for \emph{  ceiling convergents} : we denote $\lceil d_1,d_2,\cdots,d_r \rceil$, the  ceiling continued fraction expansion of $\frac{q}{p}$, as above. We also define, for $i\in\inter{1}{r}$,  the \emph{ ceiling convergent} $\frac{q_i}{p_i}=\lceil d_1,d_2,\cdots,d_i\rceil$. Hence :
    $$ \forall i\in\inter{1}{r},\esp p_i=d_ip_{i-1}-p_{i-2}, \esp q_i=d_iq_{i-1}-q_{i-2}$$
    \esp with $p_{-1}=-1, \esp p_0=0, \esp q_{-1}=0, \esp q_0=1$. The proof of this result is similar to the proof for classical partial quotients...\\

   \begin{mademo}
     (i) We suppose that we have a modular-convex sequence $(u_i)_{i\in\inter{0}{r}}$ in $\R^n$ such that $u_0=u,u_1=v$ and $u_r=w$. Let us denote $(-x_i,y_i)$ the coordinates of $u_i$ in $(u,v)$, for it is clear that all $u_i$ are integer combinations of $(u,v)$. We also denote $(d_i)_{i\in\inter{1}{r-1}}$, such that :
     $$ \forall  i\in\inter{1}{r-1}, \esp u_{i+1}=d_iu_i-u_{i-1}$$
     \esp Then, all $d_i$ are integers greater than one. We have $x_0=-1,x_1=0,x_2=1$ and $y_0=0,y_1=1,y_2=d_1$. Yet, $(x_i)_i$ and $(y_i)_i$ are modular-convex and $x_1\geqslant \max(x_0,0)$ and $y_1\geqslant \max(y_0,0)$, so ( see Remark 1) they are convex. In addition, $x_1-x_0>0$ and $y_1-y_0>0$, so  $(x_i)_{i\geqslant 0}$ and $(y_i)_{i\geqslant 0}$ are increasing sequences. Likewise, if we denote $\epsilon_i=y_i-x_i$ for $i\in\inter{0}{r}$, then $(\epsilon_i)_i$ is modular-convex, so convex and  non decreasing, since $\epsilon_0=\epsilon_1 > 0$. Thus, $\epsilon_i>0$ for all $i$. So, $0<p<q$ and half of (i) is proven.\\
     \esp With same arguments, we prove that $(d_1x_i-y_i)_{i\geqslant 2}$ is a non negative and increasing sequence, for $d_1x_2-y_2=0$ and $d_1x_3-y_3=1$. Similarly, we prove that $(y_i-(d_1-1)x_i)_{i\geqslant 1}$ is a positive and increasing sequence, for $y_1-(d_1-1)x_1=1y_2-(d_1-1)x_2$. So :
     $$  \forall i\in\inter{2}{r}, \esp d_1=\lceil y_i/x_i  \rceil $$ 
     \esp Thus, in particular, $d_1= \lceil y_r/x_r \rceil = \lceil q/p \rceil$. Now, we can do the same with $(u_1,u_2)$ instead of $(u_0,u_1)$, for $u_1$ and $u_2$ are linearly independant and all $u_i, i\in\inter{1}{r}$ are in the lattice generated by these two vectors, and so on... \\
     \esp If we denote, for $k\in\inter{1}{r},(-X_k,Y_k)$ the coordinates of $u_r$ in the basis $(u_{k-1},u_k)$ of the lattice $L$ generated by $(u_0,u_1)$, then : $X_1=x_r=p,Y_1=y_r=q$ and $X_r=0,Y_r=1$. We also have the following induction relations :
       $$ \forall k\in\inter{1}{r-1}, \esp X_k=Y_{k+1} \esp ; \esp Y_k=d_kY_{k+1}-X_{k+1} $$ 
    \esp We deduce, by obvious decreasing induction on $k$, that $\gcd(X_k,Y_k)=1$ for all $k$, so  $\gcd(p,q)=1$. \\
    \esp Now,  we have proved that $d_1=\lceil Y_1/X_1 \rceil$ and more generally that $d_k=\lceil Y_k/X_k\rceil$ for all $k\in\inter{1}{r-1}$. So, if we denote $t_k= \frac{X_k}{Y_k}$, then $t_r=0$ and :
    $$\forall k\in\inter{1}{r-1}, \esp d_k=\left\lceil \frac{1}{t_k}\right \rceil \esp ; \esp t_{k+1}=d_k-\frac{1}{t_k}=\left\{ -\frac{1}{t_k}\right\}$$

    (ii) Now, we suppose that $w=qv-pu$, with $p,q$ two coprime positive integers, such that $p<q$.  We define a finite sequence $(d_i)_{i\in\inter{1}{r-1}}$, as we use the algorithm detailed in Proposition 4. Indeed, we obtain $t_r=0$ for a certain positive integer $r$, because if we denote $t_k=\frac{X_k}{Y_k}$ the reduced fraction (  for $k\geqslant 1$), then : 
     $$\forall k\geqslant 1, \esp 0 \leqslant X_{k+1}<Y_{k+1}=X_k$$
     \esp for $t_{k+1}\in [0,1[$. Since, for all $k$, we have $t_k\in [0,1[$, then $1/t_k>1$, so $d_k\geqslant 2$. We define the sequence $(u_i)_{i\in\inter{1}{r}}$ in $L$ by :
     $$ u_0=u, \esp u_1=v, \esp \forall i\in \inter{1}{r-1}, \esp u_{i+1}=d_iu_i-u_{i-1}$$
     \esp This sequence is modular-convex and $u_r=w$, according to (i) above.   
   \end{mademo}

\newpage

   \subsection{irreducible elements of certain submonoids of $\Z$ or $\Z^2$} 
   
   \label{subsec:irreducible}
 
   $\bullet$ The following result is a complement to Lemma \ref{lem:9} and is important for next subsection : as usual, we consider the product order on $\Z^2$ and denote $C^*=\N^2\backslash \{ (0,0)\}$ the \emph{positive cone}.
   
   \begin{monlem} \label{lem:10} let $d$ and $m$ be two coprime integers such that $0<m<d$. Let $u_0=(d,0),u_1=(m,1)$ and $L=\mathrm{Span}_{\Z}(u_0,u_1)$. We denote $(d_i)_{i\in\inter{1}{r-1}}$ the coefficients of the ceiling continued development of $m/d$ and $(u_i)_{i\in\inter{0}{r}}$ the modular-convex sequence defined by $u_0,u_1$ and $u_r=(0,d)$ ( see Proposition \ref{prop:4}, for $(0,d)=-mu_0+du_1$). Reminder : $(x(u_i))_i$ is decreasing and $(y(u_i))_i$ is increasing.\\
   \esp  Let $a,b$ be two integers such that $0<a<b$ and $w=(b,-a)\in L$. We denote $\S^*=C^*+ w\Z$.\\
   
   (i) 
   $$ \min(L\cap C^*)= \{ u_i, i\in\inter{0}{r}\}$$
   (ii)
   $$ \min(L\cap \S^*)\subset \min(L\cap C^*)+w\Z$$
   (iii) \\
   $\blacktriangleright $ Case 1 : if $d<a<b$ :
   $$ \min(L\cap \S^*)= \{ u_i, i\in\inter{0}{r}\}+w\Z$$
    $\blacktriangleright $ Case 2 : if $a<d<b$ :
   $$ \min(L\cap \S^*)= \{ u_i, i\in\inter{0}{s}\}+w\Z$$
   \esp where $y(u_s)<a\leqslant y(u_{s+1})$.\\
   $\blacktriangleright $ Case 3 : if $d<a<b$ :
   $$ \min(L\cap \S^*)= \{ u_i, i\in\inter{\sigma}{s}\}+w\Z$$
   \esp where $y(u_s)<a\leqslant y(u_{s+1})$ and $x(u_{\sigma})<b+x(u_s)\leqslant x(u_{\sigma-1})$.
  
  \end{monlem}
  
\vspace*{1cm}  
  
  \textbf{Remark 1 : } as we have seen in Lemma \ref{lem:2}, $ \min(L\cap \S^*)$ is $w$-invariant, since $\S^*$ and $L$ are $w$-invariant. So, we can choose representatives of "classes modulo $w$" of minimal points in the strip :
  $$\A=\{ (x,y), x\in\N, y\in \inter{0}{a-1}\}$$
  \esp We remark that the $u_i$ that appear above are in $\A$.\\
  
  \textbf{Remark 2 : } we could summarize the above three cases in one formulation : 
   $$ \min(L\cap \S^*)= \{ u_i, i\in\inter{\sigma}{s}\}+w\Z$$
   \esp where $y(u_s)<a\leqslant y(u_{s+1})$ and $x(u_{\sigma})<b+x(u_s)\leqslant x(u_{\sigma-1})$. With additional border conditions : $x_{-1}=y_{r+1}=+\infty$. \\

    
  \begin{mademo}
  (i) we remark that for all $x,y\in\Z,xu_0+yu_1=(xd+ym,y)$, so : on every " horizontal segment"  of $\Z^2$ of length $d$, there is a unique point of $L$. On the other hand, since $m$ and $d$ are coprime, $xd+ym=0$ only if $y$ is a multiple of $d$ : so on every " vertical segment"  of $\Z^2$ of length $d$, there is a unique point of $L$. We deduce that  $u_0,u_1,u_r\in  \min(L\cap C^*)$ and so $\min(L\cap C^*)=\min(L\cap R)$, with $R=(\inter{0}{d}^2)^*$.\\
  \esp  Now, according to Lemma \ref{lem:9}, $\min(L\cap R)$ can be enumerated by a modular-convex sequence $(v_i)_{i\in\inter{0}{r'}}$ such that $(x(v_i))_i$ is decreasing and $(y(v_i)_i$ is increasing. So, $v_0=u_0,v_1=u_1$ and $v_{r'}=u_r$. Now, using Proposition \ref{prop:4}, we deduce that $r'=r$ and $v_i=u_i$ for all $i\in\inter{0}{r}$. \\
  
  (ii) let $u\in \min(L\cap \S^*)$ and $u\in \A$. Then, $u\in L\cap C^*$, since $\A\subset C^*$. If there exists $v\in L\cap C^*$ such that $v<u$, then $v\in L\cap \S^*$, which contradicts the fact that $u$ is minimal in $L\cap \S^*$. So, $u\in \min(L\cap C^*)$. We have proved that :
  $$  \min(L\cap \S^*)\cap \A \subset \min(L\cap C^*)$$
  \esp We deduce (ii) with Remark 1 above. \\
  
  (iii) \\
  $\blacktriangleright $ Case 1 : if $d<a<b$. Let $i\in\inter{0}{r}$, then $u_i\in \A$. Moreover : 
  $$0\leqslant x(u_i)\leqslant d<b \esp \text{ and } \esp 0\leqslant  y(u_i) \leqslant d <a$$
  \esp  Now, we suppose that we have $v\in L\cap \S^*$ such that $v<u_i$. Then, there exists an integer $k$, such that $v'=v-kw\in C^*$. But, $k$ can not be positive, for we would have $x(v')=x(v)-kb<0$ and $k$ can not be negative,  for we would have $y(v')=y(v)+ka<0$. Hence $k=0$ and  $v\in C^*$, which contradicts (i)... \\
  \esp So, $u\in  \min(L\cap \S^*)$. We conclude with (ii).\\
   $\blacktriangleright $ Case 2 : if $a<d<b$. We use the same arguments as in Case 1, if $i\in\inter{0}{s}$. If, $i\in\inter{s+1}{r}$, then $u_i\not\in\A$. In addition, $x(u_i)\geqslant 0$ and $y(u_i)\geqslant a$, so $u_i>(-b+d,a)=u_0-w$. But, $u_0-w\in L\cap \S^*$, so $u_i\not\in   \min(L\cap \S^*)$.\\
    $\blacktriangleright $ Case 3 : if $a<b<d$. If $i\in \inter{s+1}{r}$, then $u_i\not\in \A$. If $i\in \inter{0}{\sigma -1}$, let $v=u_s+w\in L\cap \S^*$. Then $v<u_i$, for $x(v)=b+x(u_s)\leqslant x(u_i)$ and $y(v)=y(u_s)-a<0\leqslant y(u_i)$. So  $u_i\not\in   \min(L\cap \S^*)$.\\
    \esp Now, let $i\in \inter{\sigma}{s}$. We suppose that there exists $v\in L\cap \S^*$ such that $v<u_i$. Then, we claim that $v\in C^*$. Indeed : we have an integer $k$ such that $v'=v-kw\in \A$. \\
    ----- if $k<0$, then $y(v)=y(v')-ka\geqslant a$, which contradicts $v<u_i$, for $y(u_i)<a$. \\
    ----- if $k>0$, then 
    $$0\leqslant x(v')\leqslant x(v)-kb\leqslant x(u_i)-b<x(u_s)$$
    \esp But, $v'\in L\cap C^*$, so there exists $j\in\inter{0}{r}$, such that $v'\geqslant u_j$ ( see (i)). The above inequality proves that $j>s$ and so $y(v')\geqslant y(u_j)\geqslant a$, which contradicts $v'\in \A$. \\
    \esp So, $v\in L\cap C^*$ and $v<u_i$ : it is impossible according to (i). 
    \end{mademo}
  

 $\bullet$ the previous Lemma is important, because these minimal sets are related with minimal generators and irreducible elements ( see beginning of section 2.) in the following meaning : \\
 \esp Let $M$ be a submonoid of $(\R^n,+)$ such that $0$ is the unique invertible element of $M$. $\mathrm{Gen}(M)$ denotes  the intersection of all sets of generators ( as a submonoid) of $M$. The elements of $\mathrm{Gen}(M)$ are the \emph{minimal generators } of $M$.

  
  \begin{monlem} \label{lem:11} let $L$ be a lattice of $\Z^2$.
  $$ \min(L\cap (\N^2)^*)=\Gen(L\cap \N^2)=\Irr(L\cap \N^2)$$
   \end{monlem}
   \begin{mademo} we denote $C^*=(\N^2)^*$.\\
  - " $\Irr  \subset \min $" : if $u\in L\cap C^*$ is not minimal in $L\cap C^*$, then there exists $v\in L\cap C^*$ such that $v<u$. Note $v'=u-v$, then $v'\in L\cap C^*$ and $u=v+v'$, so $u\not \in \Irr(L\cap \N^2)$. \\
  -  " $\min \subset \Gen$" : if $m\in \min(L\cap C^*)$ and $(u_k)_{k\in I}$ is a system of generators of $L\cap \N^2$, say $L\cap \N^2=<u_k,k\in I>$. Then , $m=\sum_{k\in I}\alpha_ku_k$, where $\alpha_k\in\N$ for all $k$. Yet, $m\not = (0,0)$, thus there exists a $k\in I$, such that $\alpha_k\geqslant 1$ and so $m\geqslant u_k$, since all $u_i\geqslant (0,0)$. Yet $u_k\in L\cap C^*$, so $m=u_k$ and $m$ belongs to every set of generators of $L\cap \N^2$.\\
  - " $\Gen \subset \Irr $" : suppose that $u\in L'\cap \N^2$ and $u\not\in \Irr(L\cap \N^2)$. Then, we can find $u_1$ and $u_2$ in $L'\cap \N^2$, such that $u=u_1+u_2$. So that we can replace $u$ by $(u_1,u_2)$ in any system of generators of $L\cap \N^2$ that would contain $u$ : $u\not\in \Gen(L\cap \N^2)$.   \end{mademo}

  \newpage

   $\bullet$ We now state a result that will we improved later ( see  Theorem \ref{theo:5}). It gives an expression of $\Irr(\frac{<a,b>}{d})$ in terms of $a,b,d$ and is a  consequence of Lemma \ref{lem:8} and Proposition \ref{prop:2}.
   
   \begin{maprop} \label{prop:5} let $a,b,d$ be three pairwise coprime positive integers such that $d\not \in <a,b>$ and $d>1$. We denote $m\in \inter{1}{d-1}$, such that $am+b=0[d]$. We denote $\lceil d_1,\cdots,d_{r-1}\rceil$ the  ceiling continued fraction expansion of $\frac{d}{m}$ ( see the algorithm of Proposition \ref{prop:4}) and define  sequences of integers, by double induction :
     $$ y_0=0, y_1=1, \esp \forall i\in\inter{1}{r-1}, \esp y_{i+1}=d_iy_i-y_{i-1} $$
   $$ x_0=d, x_1=m, \esp \forall i\in\inter{1}{r-1}, \esp x_{i+1}=d_ix_i-x_{i-1} $$
   \esp in addition, we denote $y_{r+1}=x_{-1}=+\infty$.
   $$ n_0=a, n_1=\frac{am+b}{d}, \esp \forall i\in\inter{1}{r-1}, \esp  n_{i+1}=d_in_i-n_{i-1} $$
   \esp Then $x_r=0,\esp y_r=d,\esp n_r=b,\esp (x_i)_i$ is decreasing and $(y_i)_i$ is increasing. In addition : \\
   
   $\blacktriangleright$ Case 1 : if $d<a<b$
   $$ \Irr\left(\frac{<a,b>}{d}\right)=\{ n_i, i\in \inter{0}{r} \}$$

   $\blacktriangleright$ Case 2 : if $a<d<b$
   $$\Irr\left(\frac{<a,b>}{d}\right)=\{ n_i, i\in \inter{0}{s} \}$$
 \esp   where $s$ is such that $y_s<a\leqslant y_{s+1}$. \\
 
  $\blacktriangleright$ Case 3 : if $a<b<d$
   $$\Irr\left(\frac{<a,b>}{d}\right)=\{ n_i, i\in \inter{\sigma}{s} \}$$
 \esp   where $s$ is such that $y_s<a\leqslant y_{s+1}$ and $\sigma$ is such that $x_{\sigma}<x_s+b\leqslant x_{\sigma-1}$.

   \end{maprop}

   \begin{mademo}
   If we denote $q_i/p_i$ the reduced fraction of $\lceil d_1,d_2,\cdots,d_i\rceil$ ( see the remark below Proposition \ref{prop:4}), then, with obvious induction : 
   $$ \forall i \in\inter{0}{r}, \esp x_i=-dp_{i-1}+mq_{i-1}, \esp y_i=q_{i-1}, \esp n_i=-ap_{i-1}+\tau q_{i-1} $$
   \esp where $\tau=\frac{am+b}{d}$. So $x_r=-dp_{r-1}+mq_{r-1}=-dm+md=0, \esp y_r=q_{r-1}=d$ and $n_r=-ap_{r-1}+\tau q_{r-1}=-am+\tau d =b$. Now, by obvious induction, $(y_i)_i$ is increasing and $(x_i)_i$ is decreasing. \\
   \esp Let us consider now ( with the product order in $\Z^2$) :
   $$L=\{ (x,y)\in\Z^2, ax+by=0[d] \} \esp ; \esp M=\min(L\cap \N^2\backslash \{ (0,0)\} )$$
   \esp  But $(d,0)$ and $(0,d)$ are in $M$, so $M=\min (L\cap ((\inter{0}{d})^2)^*)$. So, according to Lemma \ref{lem:7}, $M$ can be enumerated by a finite modular-convex sequence $(u_i)_{i\in\inter{0}{r}}$ such that $(x(u_i))_i$ is decreasing and $(y(u_i))_i$ is increasing. We have $u_0=(d,0)$ and $u_r=(0,d)$. Moreover, $(m,1)\in L$, so $(m,1)\in M$ and $u_1=(m,1)$. Now, $u_r=-mu_0+du_1$ and $0<m<d$. In addition, $\gcd(m,d)=1$, since $\gcd(b,d)=1$. With Proposition \ref{prop:4}, $(u_i)_i$ satisfy the induction relation :
   $$ \forall i\in\inter{1}{r-1}, \esp u_{i+1}=d_iu_i-u_{i-1}$$
   \esp So, $x_i=x(u_i)$ and $y_i=y(u_i)$ for $i\in\inter{0}{r}$.  We deduce : $n_i=\psi(u_i)$ for $i\in\inter{0}{r}$, where $\psi : (x,y)\to \frac{ax+by}{d}$. \\
   \esp At the end, results for $ \Irr\left(\frac{<a,b>}{d}\right)$ are direct consequences of Proposition \ref{prop:2} ...   
   \end{mademo}

   \newpage

   \subsection{$\N$-independance of sequences of integers} 
   
   $\bullet$ Now, we are almost ready to give the  solutions of the " reverse problem", that is to say : to provide all triplets of coprime positive integers $a,b,d$ such that $0<a<b$ and $\Irr(\frac{<a,b>}{d})$ is a given finite subset of $\N$. At the beginning of this section, we have reminded the necessary conditions for a finite subset $I$ of $\N$ to be the set of minimal generators of some $\frac{<a,b>}{d}$. One of them is of course, that elements of $I$ are $\N$-independant. Let us precise this notion and give a useful criterion :   \\
   \esp First, we remind a notation : for a subset $A$ of a monoid $M$, $<A>$ denotes the submonoid of $M$ generated by $A$.
   
   \begin{madef} \label{def:2}
       elements of a subset $A$ of a monoid $M$ are \emph{$\N$-independant} if and only if :
   $$ \forall a\in A, a\not\in <A\backslash \{ a \}>$$
   \end{madef}
   \textbf{Remark : } we will name also " $A$ is $\N$-independant" such a definition. \\
   
   \esp For our purpose, we deal with the particular case of the monoid  $\N$. Let $A$ be a subset of  $\N$. We claim that :
   $$ A \text{ is $\N$-independant } \Leftrightarrow \Irr(<A>)=A$$
   
  \textbf{Question : } given a finite set $A$ of positive integers, do we have an efficient algorithm to know if $A$ is $\N$-independant or not ? \\
   \esp It is difficult to answer this question in general, but in the particular case when $A$ can be enumerated by a modular-convex sequence, we have a simple criterion, stated in Proposition \ref{prop:6}. But, first, we need an arithmetic Lemma, that will be also useful later ( see Theorem \ref{theo:2} and \ref{theo:3}) :
   
   \begin{monlem} \label{lem:12} let $a,b,c$ be three non null integers such that $b$ and $c$ are coprime.\\
    There exists an infinity of positive integer $k$, such that $a$ and $kc-b$ are coprime.
   \end{monlem}
   \begin{mademo} for any prime divisor $p$ of $a$, we consider :
   $$ K_p=\{ k\in\Z, kc=b[p]\} $$
   if $c=0[p]$, then $K_p$ is empty, since $b\not =0[p]$ ( $b$ and $c$ are coprime). Else, $K_p$ is a class of congruence modulo $p$, say : $K_p=\{ k\in\Z,k=r_p[p]\}$. So, if we denote $\P(a)$ the finite set of prime divisor  of $a$ that does not divide $c$  :
   $$ \forall k\in\Z, ( \gcd(a,kc-b)=1 \Leftrightarrow \forall p\in \P(a) , k\not = r_p [p] )$$
   \esp With chinese remainder Lemma, we obtain that $k$ is such that $kc-b$ and $a$ are coprime if and only if $k$ take certain values modulo $q=\prod_{p\in\P(a)}p$. 
     \end{mademo}

   \esp Now, we propose a result that gives a criterion for a modular-convex sequence of positive integers to be $\N$-independant, but also the description of $\Irr(<A>)$ for $A$ enumerated by such a sequence. \\
  \esp We will use what we have named the " ceiling continued fraction expansion " of a rational $x$, denoted $\lceil  d_1,d_2,\cdots,d_r \rceil$, as well as $\left(\frac{q_i}{p_i}\right)_{i>0}$ the finite sequence of " ceiling partial quotients" of $x$, defined by ( see \ref{subsec:modular}):
  $$ p_{-1}=-1,p_0=0, \esp \forall i\in\inter{1}{r}, p_i=d_ip_{i-1}-p_{i-2}$$
   $$ q_{-1}=0,q_0=1, \esp \forall i\in\inter{1}{r}, q_i=d_iq_{i-1}-q_{i-2}$$

   \begin{maprop} \label{prop:6} let $(n_i)_{i\in\inter{0}{r}}$ be a modular-convex sequence of positive integers such that $n_0$ and $n_1$ are coprime. We denote $(d_i)_i$ the sequence of integers such that :
   $$ \forall i\in\inter{1}{r-1}, \esp d_i\in\N, \esp d_i\geqslant 2, \esp n_{i+1}=d_in_i-n_{i-1}$$
   \esp We also denote $\frac{q_i}{p_i}$ the reduced fraction of $\lceil d_1,d_2,\cdots,d_i \rceil$ for all $i\in\inter{1}{r-1}$ ( see \ref{subsec:modular}).  More explicitly, $(q_i)_{i\in\inter{-1}{r-1}}$ is the following sequence :\\
   $$ q_{-1}=0, q_0=1, \esp \forall i\in\inter{0}{r-1}, \esp  q_i=d_iq_{i-1}-q_{i-2} $$
   (i) 
   $$ \{ n_i,i\in\inter{0}{r}\}   \text{ is $\N$-independant } \Leftrightarrow  q_{r-1}<\min(n_0,n_r) $$
   (ii)
   $$ \Irr( < n_i,i\in\inter{0}{r}>)= \{ n_i, i\in\inter{\sigma}{s} \} $$
    \esp   where $s$ is such that $y_s<n_0\leqslant y_{s+1}$ and $\sigma$ is such that $x_{\sigma}<x_s+n_r\leqslant x_{\sigma-1}$, with  $ x_{-1}= y_{r+1}=+\infty$ and : 
      $$  \forall i\in\inter{0}{r}, \esp y_i=q_{i-1} \esp ; \esp x_i=p_{r-1}q_{i-1}-q_{r-1}p_{i-1}$$

   \end{maprop}

   \vspace*{0.5cm}

   \textbf{Remark 1 : } as we have already seen in \ref{subsec:modular},  the algorithm that computes the " ceiling continued fraction expansion" of a fraction is of exponential complexity in worst cases. But here, the entries are the $(d_i)_i$, so computing $q_{r-1}$ is very fast...\\

   \textbf{Remark 2 : }  again, we have already mentioned that, for a modular-convex sequence $(n_j)_{j\in\inter{0}{r}} $, we have $\gcd(n_{i-1},n_{i})=\gcd(n_0,n_1)$ for all $i\in\inter{1}{r}$, so $\gcd((n_j)_{j\in\inter{0}{r}})=\gcd(n_0,n_1)$.\\
   \esp A bit more general statement of Proposition 6 would propose the criterion :
    $$ \{ n_i,i\in\inter{0}{r}\}   \text{ is $\N$-independant } \Leftrightarrow  q_{r-1}<\dfrac{\min(n_0,n_r)}{\gcd(n_0,n_1)}$$
   
    \textbf{Remark 3 : } The proof will appear to be an application of Lemma \ref{lem:10}. It is quite direct, if the first and last term of the sequence are coprime, but more intricate when they are not. We will use an argument of limit in that last case...\\

   \begin{mademo}
   (i) is a direct consequence of (ii), for $y_r=q_{r-1}=x_0$ and $x_r=0$. ( see below) \\
   (ii) let $u_0=(q_{r-1},0),u_1=(p_{r-1},1)$ and $(u_i)_{i\in\inter{0}{r}}$ defined by :
   $$ \forall i\in\inter{1}{r-1},u_{i+1}=d_iu_i-u_{i-1}$$
   \esp Then $x(u_i)=x_i$ and $y(u_i)=y_i$ for all $i\in\inter{0}{r}$, by obvious induction. So, $u_r=(0,q_{r-1})$. We also remark that :
   $$ \forall i\in\inter{0}{r}, \esp x_i=p_{r-1}q_{i-1}-q_{r-1}p_{i-1}= \begin{vmatrix} p_{r-1} & p_{i-1} \\ q_{r-1} & q_{i-1} \end{vmatrix}$$
\esp We denote $a=n_0$ and $b=n_r$ and define a map $\psi$ from $\Z^2$ to $(1/q_{r-1})\Z$ :
   $$ \forall x,y\in\Z,\esp \psi(x,y)=\frac{ax+by}{q_{r-1}}$$
   \esp We have $\psi(u_0)=n_0$ and $\psi(u_r)=n_r$. We claim that : $\psi(u_i)=n_i$ for all $i$. Indeed, if we denote $\epsilon_i=\psi(u_i)-n_i$ for all $i$, then $\epsilon_0=0$ and : 
   $$ \forall i\in\inter{1}{r-1},\epsilon_{i+1}=d_i\epsilon_i-\epsilon_{i-1}$$
   \esp Then, by obvious induction, $\epsilon_i=\epsilon_1y_i$ for all $i\in\inter{0}{r}$. But, $\epsilon_r=0$ and $y_r=q_{r-1}$, so $\epsilon_1=0$. \\
   \esp We deduce :
   $$ \forall i\in\inter{0}{r}, \esp \psi(u_i)=n_i$$
    \esp Now, we denote $L=\mathrm{Span}_{\Z}(u_0,u_1)$. With Lemma \ref{lem:10} and \ref{lem:11}, we obtain :
    $$\Gen(L\cap\N^2)=\min(L\cap (\N^2)^*)=\{ u_i,i\in \inter{0}{r}\}$$
    \esp Let $S=<n_i, ,i\in \inter{0}{r}>$. Then, $S=\psi(L\cap \N^2)$. \\
    \esp We also remark that, if we denote $w=(b,-a)$, then  $w=n_1u_0-au_1 \in L$, for $n_1=\frac{ap_{r-1}+b}{q_{r-1}}$.\\
    
   $\blacktriangleright$ Case 1 : if $\gcd(n_0,n_r)=1$. Then, $\ker(\psi)=w\Z$, so :
    $$ \psi^{-1}(S)=  ( L\cap \N^2)+w\Z= L\cap (\N^2+w\Z) $$ 
   \esp So with Lemma \ref{lem:1} and \ref{lem:3}, $\psi^{-1}(\Irr(S))=\min(L\cap \S^*)$, where $\S^*=(\N^2)^*+w\Z$. We conclude with Lemma \ref{lem:10}.  \\
  $\blacktriangleright$ Case 2 : if $\gcd(n_0,n_r)>1$. According to Lemma \ref{lem:12}, there exist an infinite number of integer $d_r$, such that $d_r\geqslant 2$ and $n_0,n_{r+1}$ are coprime, if we denote $n_{r+1}=d_rn_r-n_{r-1}$. For these values of $n_{r+1}$, the sequence $(n_j)_{j\in\inter{0}{r+1}}$ satisfy the assumptions of Case 1. We denote :
  $$ q_r=d_rq_{r-1}-q_{r-2} \esp ; \esp p_r=d_rp_{r-1}-p_{r-2}$$
  \esp So, we obtain :
   $$      \Irr( < n_i,i\in\inter{0}{r+1}>)= \{ n_i, i\in\inter{\sigma}{s'} \} $$ 
  \esp   where $s'$ is such that $y'_{s'}<n_0<y'_{s'+1}$ and $\sigma$ is such that $x'_{\sigma}<x'_{s'}+n_{r+1}<x'_{\sigma-1}$, with $ x'_{-1}= y'_{r+2}=+\infty$ and $(x'_i)_i, ( y'_i)_i$ defined by :
  $$ \forall i\in\inter{-1}{r+1}, \esp y'_i=q_{i-1} \esp ; \esp  x'_i=\begin{vmatrix} p_r & p_{i-1}\\ q_r & q_{i-1} \end{vmatrix}$$ 
  \esp So, $y'_i=y_i$ for all $i\in\inter{0}{r}$, while $y_{r+1}=+\infty$ and $y'_{r+1}=q_r$. So, we just have to choose $d_r$ large enough so that $q_r >n_0$, to obtain $s'=s$ ( $s$ being defined as in Proposition \ref{prop:5}). Again, we can choose $d_r$ large enough such that $n_{r+1}\in <n_0,n_1>$. Then,  $ < n_i,i\in\inter{0}{r+1}>= < n_i,i\in\inter{0}{r}>$. \\
  \esp Now, we denote $(X_i)_i$ the sequence defined by :
  $$ \forall i\in \inter{0}{r-1}, \esp X_i=\begin{vmatrix} p_{r-2} & p_{i-1}\\ q_{r-2} & q_{i-1} \end{vmatrix}$$   
  \esp Then, we can write :
  $$ \forall i\in\inter{0}{r-1},\esp  x'_i=tx_i-X_i \esp \text{ where } t=d_r $$
  \esp We also have : $n_{r+1}=tn_r-n_{r-1}$. So, the conditions on $(x'_i)_i$, mentioned above, become :
  $$\begin{cases} n_{r-1} + X_s-X_{\sigma}< t(n_r+x_s-x_{\sigma})  \hspace{2.7cm} (1) \\
  t(n_r+x_s-x_{\sigma -1} ) \leqslant n_{r-1}+X_s-X_{\sigma -1} \hspace{2cm} (2)\end{cases} $$
  \esp Since $t$ can be as large as we want,  condition (1) is equivalent to :
   $$ n_r+x_s-x_{\sigma}>0 \text{ or } \begin{cases} n_r+x_s-x_{\sigma}=0 \\ n_{r-1} + X_s-X_{\sigma}< 0 \end{cases} $$
  \esp Suppose now that $x_{\sigma}-x_s=n_r=b$. We have $u_{\sigma},u_s\in L$ and $w\in L $ ( see remark just before Case 1), so $v=u_{\sigma}-u_s-w \in L$. But $x(v)=0$, so $y(v)=kq_{r-1}$, for some integer $k$ ( there is only one point of $L$ on each vertical line of $\Z^2$ of length $q_{r-1}$). So, $n_{\sigma}-n_s=kn_r$ : contradiction with the fact that $n_s,n_{\sigma}\in  \Irr( < n_i,i\in\inter{0}{r}>)$. We have proved that $n_r+x_s-x_{\sigma}\not = 0$, so condition (1) and (2) ( that we treat with similar arguments) are equivalent ( for $t$ large enough) to :
   $$ x_{\sigma} < x_s+n_r < x_{\sigma-1} \Leftrightarrow  x_{\sigma} < x_s+n_r \leqslant  x_{\sigma-1}  $$ 
   \end{mademo}

   \newpage

   \subsection{solutions}
   
   \label{subsec:solutions}

  $\bullet$ We remind our " reverse problem" :  if $r$ is a positive integer and $(n_i)_{i\in \inter{0}{r}}$ is a finite sequence of positive integers satisfying conditions (i),(ii),(iii) detailed at the beginning of section \ref{sec:reverse}, what are the triplets $(a,b,d)$ of coprime positive integers, such that $a<b$ and $\Irr(\frac{<a,b>}{d})=\{ n_i,i \in\inter{0}{r} \}$ ? \\
  
  \esp We also recall some notations : we have defined  by induction in \ref{subsec:modular} the notation $\lceil d_0,d_1,d_2,\cdots,d_{i}\rceil$, as follows :
  $$\lceil d_0,d_1,d_2,\cdots,d_{i}\rceil=d_0-\frac{1}{\lceil d_1,d_2,\cdots,d_{i}\rceil}$$
 \esp This is a variant of the usual  convergents of continued fraction and has been named : " ceiling convergents". The $(d_i)_i$ are integers such that $d_i\geqslant 2$ for all positive $i$.\\
 \esp To every real $\alpha$, we can associate such a sequence of $(d_i)_i$ which is infinite if $\alpha$ is irrational and finite else ( plus an infinite $(+\infty)$-tail).  Then, $\alpha$ is the limit of the sequence of above " ceiling convergents". The algorithm to obtain these coefficients has been detailed for rationals in \ref{subsec:modular} and is generalized in \ref{subsec:CFE} to the case of irrationals.\\
 
   \textbf{Remark : } in the proof of next three Theorems, we will have integers $a,b,d,m,n_1$ such that $b=dn_1-am, \gcd(d,m)=1$ and $\gcd(a,n_1)=1$. We will deduce then :
  $$\gcd(a,b)=\gcd(a,d)=\gcd(b,d)$$

 $\bullet$ First, we look for solutions $(a,b,d)$ such that $d<a<b$. In that case, the situation is very constraint and we obtain a finite number of solution ( one in general) :
  
  \begin{montheo} [ Reverse-problem, solutions $(a,b,d)$ such that  $d<a<b$].  \label{theo:1}\\
  let  $I$ be a finite subset of $\N\backslash\{ 0,1\}$. \\
  $\blacklozenge$ (i) if $I$ is $\N$-independant and  can be enumerated by a modular-convex  sequence $(n_k)_{k\in\inter{0}{r}}$ , such that : $ \gcd(n_0,n_r)=1 \text{ and } n_0<n_r $, then $I=\Irr(\frac{<a,b>}{d})$, where $a=n_0,b=n_r$ and $d\in\inter{1}{a-1}$ such that $dn_1=b[a]$.\\
 $\blacklozenge$ (ii) there are as many triplets of coprime integers $(a,b,d)$ such that $1<d<a<b$ and $I=\Irr(\frac{<a,b>}{d})$, as there are sequences of integers $(n_k)_k$ that verify  conditions of (i).

  \end{montheo}

 \begin{mademo}
 (i) we denote $a=n_0,b=n_r$ and $(d_i)_{i\in\inter{1}{r-1}}$ the sequence of integers such that :
 $$ \forall i \in\inter{1}{r-1}, \esp d_i\geqslant 2, \esp n_{i+1}=d_in_i-n_{i-1}$$
 \esp Let $m,d$ be coprime positive integers such that :
  $$1-\frac{m}{d}=\lceil 1,d_1,d_2,\cdots,d_{r-1}\rceil \esp \text{ that is : } \esp \frac{d}{m}=\lceil d_1,d_2,\cdots,d_{r-1}\rceil$$
  \esp Since $I$ is $\N$-independant, we can use Proposition \ref{prop:6} and claim that : $d<\min(a,b)$. \\
  \esp In addition, by obvious induction ( see the proof of Proposition \ref{prop:5}, for example) :
  $$\forall i\in\inter{0}{r},\esp n_i=q_{i-1}n_1-p_{i-1}n_0 \esp ; \esp \text{  where } \esp \frac{q_i}{p_i}=\lceil d_1,d_2,\cdots,d_{i}\rceil$$
 \esp so $b=n_r=dn_1-mn_0$ and $am+b=dn_1=0[d]$. In addition, $a$ and $b$ are coprime, so $a,b,d$ are pairwise coprime ( see remark above). Yet, $m\in\inter{1}{d-1}$ and $d\not\in <a,b>$, so we can apply Proposition \ref{prop:5} ( Case 1) and we obtain our result.\\ 
 
 (ii) direct consequence of (i) and Proposition \ref{prop:5} ( Case 1). 
 \end{mademo}

  \textbf{Examples : } at the beginning of this section, we have seen an example of a $\N$-independant set of integers, namely $I=\{ 7,11,59\}$, that can be enumerated by two different modular-convex sequences verifying conditions of Theorem \ref{theo:1} (i) : $(7,11,59)$ ( $d_1=6$) and $(11,7,59)$ ( $d_1=10$), so we have : 
  $$  \Irr(\frac{<7,59>}{6})=\Irr(\frac{<11,59>}{10})=\{ 7,11,59 \}$$
  \esp What about $I=\{ 10,17,24 \}$ ? $I$ is $\N$-independant ( esay to check) and is enumerated by the modular-convex sequence $(10,17,24)$ ( $d_1=2$), but $10$ and $24$ are not coprime. We also remark that $(17,10,24)$ is not modular-convex. So, there is no enumeration of $I$ by a sequence verifying conditions of (i) : for all triplets $(a,b,d)$ of pairwise coprime positive integers such that $d<a<b$, we can claim that $I\not = \Irr(\frac{<a,b>}{d})$.  \\

  $\bullet$ Secondly, we consider solutions such that $a<d<b$. In that case, we can extend in one direction the sequence $(n_i)_i$ in a modular-convex sequence. We then obtain an infinite number of solutions :

  \begin{montheo} [ Reverse-problem, solutions $(a,b,d)$ such that $a<d<b$]. \label{theo:2} \\
  let  $I$ be a finite subset of $\N\backslash\{ 0,1\}$ that is $\N$-independant and  can be enumerated by a modular-convex  sequence $(n_k)_{k\in\inter{0}{r}}$, such that $n_0$ and $n_1$ are coprime.\\
  \esp  We denote $(d_i)_{i\in\inter{1}{r-1}}$ the sequence of integers such that :
   $$ \forall i \in\inter{1}{r-1}, \esp d_i\geqslant 2, \esp n_{i+1}=d_in_i-n_{i-1}$$
   $\blacklozenge$ (i) Let $\rho$ be an integer such that $\rho>r$ and $(d_i)_{i\in \inter{r}{\rho-1}}$ a sequence of integers bigger than $1$.\\
   \esp  We denote $d$ and $m$ the coprime positive integers such that $  \frac{d}{m}=\lceil d_1,d_2,\cdots,d_{\rho-1}\rceil$ and $\frac{q_i}{p_i}$ the reduced fraction of $\lceil d_1,d_2,\cdots,d_i\rceil$ for all $i\in\inter{1}{\rho-1}$. \\
    \esp We denote $a=n_0$ and $b=dn_1-mn_0$. Then :

  $$ \begin{cases} \gcd(a,b)=1 \\ a\leqslant q_r \\ d<b \end{cases} \Rightarrow \begin{cases} a<d<b \\ \Irr\left(\frac{<a,b>}{d}\right)=I  \end{cases}$$ 
   $\blacklozenge$ (ii) all triplets of pairwise coprime $(a,b,d)$ such that  $\Irr\left(\frac{<a,b>}{d}\right)= I$ and $a<d<b$ can be obtained as in (i).\\   
    $\blacklozenge$ (iii) there is an infinity  of pairwise coprime $(a,b,d)$ such that  $\Irr\left(\frac{<a,b>}{d}\right)=I,a<d<b$ and  $\rho=r+1$ in (i).
  
  \end{montheo}

   \begin{mademo}
   (i) we suppose that $a$ and $b$ are coprime, $a\leqslant q_r$ and $d<b$. We have $d=q_{\rho-1}$ and $(q_i)_i$ is an increasing sequence, so : $ a\leqslant q_r \leqslant d<b $.\\
   \esp Yet, $a\not =d$, for $\gcd(a,b)=1$. So $a<d<b$. We use Proposition \ref{prop:5} ( case 2) and obtain :
   $$ \Irr\left(\frac{<a,b>}{d}\right)= \{ n_i, i\in\inter{0}{s}\}$$
   \esp where $s$ is the index such that $y_s<a\leqslant y_{s+1}$ and $y_i=q_{i-1}$ for all $i$ ( obvious induction on $i$). But, $a\leqslant q_r$ by hypothesis and $q_{r-1}<a$, because $I$ is $\N$-independant ( see  Proposition \ref{prop:6}). So $s=r$ and we obtain our result (i). \\ 
  (ii)  a direct consequence of Proposition \ref{prop:5} ( case 2).   \\
  (iii) using (i), the question is to find what values of $d_r$, integer greater than $1$, are such that : 
  $$\gcd(a,d_rn_r-n_{r-1})=1 \esp (1) \esp ; \esp  a< d_rq_{r-1}-q_{r-2}<d_rn_r-n_{r-1} \esp (2) $$
  \esp Indeed, if we choose $\rho=r+1$, then $d=q_r=d_rq_{r-1}-q_{r-2}$ and $b=n_{r+1}=d_rn_r-n_{r-1}$. But, $q_{r-1}<n_r$, for $I$ is $\N$-independant ( see  Proposition \ref{prop:6}). So, condition (2) is valid for $d_r$ large enough. Now, $n_{r-1}$ and $n_r$ are coprime, since $n_0,n_1$ are coprime ( see remark   1 at the beginning of \ref{sec:reverse}). So, with Lemma \ref{lem:12}, condition (1) is valid for an infinite number of positive integers $d_r$... that proves our result.
   \end{mademo}

   \textbf{Examples : } we go back to our previous examples at the end of Theorem \ref{theo:1}.\\
   \esp We take $I=\{ 7,11,59 \}$, we denote $a=7$ and we search integers $d_2$ such that $d_2\geqslant 2,b=59d_2-11, d=6d_2-1$ and $\gcd(a,b)=1, a<d<b$. We remark that the condition $a<d<b$ will always be verified here and that the condition $\gcd(a,b)=1$ is equivalent to $d_2\not = -1[7]$. So :
   $$ \forall k\in \inter{2}{+\infty}, k\not = -1[7] \Rightarrow \Irr\left(\frac{<7,59k-11>}{6k-1}\right)= \{ 7,11,59 \} $$  
   \esp Same arguments with the enumeration of $I$ by $(11,7,59)$, that is also modular-convex, gives :
    $$ \forall k\in \inter{2}{+\infty}, k\not = -1[11] \Rightarrow \Irr\left(\frac{<11,59k-7>}{10k-1}\right)= \{ 7,11,59 \} $$  
  
  \esp We consider now $I=\{ 10,17,24 \}$ and obtain with similar arguments :
  $$ \forall k\in \inter{6}{+\infty}, k\not = 3 [5] \Rightarrow \Irr\left( \frac{<10,24k-17>}{2k-1}\right)=\{ 10,17,24 \}$$
 \esp What do we obtain for $k=2,4,5$, that is for values of $k$ such that $b>a>d$ ? Using Proposition 5, we obtain ( $I=\{ 10,17,24 \}$) :
 $$ \Irr\left( \frac{<10,31>}{3}\right)=I\cup \{ 31 \} \esp ; \esp \Irr\left( \frac{<10,79>}{7}\right)=I\cup \{79 \} \esp ; \esp \Irr\left( \frac{<10,103>}{9}\right)=I\cup \{ 103\}$$

 
 $\bullet$ Finally, we search  the solutions such that $a<b<d$. In that case,   we can extend in both directions the sequence $(n_i)_i$ in a modular-convex sequence. We then  obtain an infinite number of solutions :

 \begin{montheo}[ Reverse-problem : solutions $(a,b,d)$ such that $a<b<d$]. \label{theo:3}\\
  let  $I$ be a finite subset of $\N\backslash\{ 0,1\}$ that is $\N$-independant and  can be enumerated by a modular-convex  sequence $(n_k)_{k\in\inter{0}{r}}$, such that $n_0$ and $n_1$ are coprime.\\
  \esp  We denote $(d_i)_{i\in\inter{1}{r-1}}$ the sequence of integers such that :
   $$ \forall i \in\inter{1}{r-1}, \esp d_i\geqslant 2, \esp n_{i+1}=d_in_i-n_{i-1}$$
   $\blacklozenge$ (i) Let $\rho$ and $\nu$ be integers such that $\rho > r$ and $\nu < 0$. Let  $(d_i)_{i\in \inter{\nu+1}{\rho-1}}$ be an extension of $(d_i)_{i\in\inter{1}{r-1}}$ in a sequence of integers bigger than $1$, that defines an extension of $(n_i)_{i\in\inter{0}{r}}$ in a modular-convex sequence $(n_i)_{i\in\inter{\nu}{\rho}}$, such that $n_i>0$ for all $i$.\\
   \esp  We denote $d$ and $m$ the coprime positive integers such that $  \frac{d}{m}=\lceil d_{\nu+1},\cdots,d_{\rho-1}\rceil$ and $\frac{q_i}{p_i}$ the reduced fraction of $\lceil d_{\nu+1},d_{\nu+2},\cdots,d_i\rceil$ for all $i\in\inter{\nu+1}{\rho-1}$. Moreover $q_{\nu-1}=0,q_{\nu}=1$ and $p_{\nu-1}=-1,p_{\nu}=0$.\\
    \esp We denote $a=n_{\nu}$ and $b=dn_{\nu+1}-mn_{\nu}$ and $x_i=mq_{i-1}-dp_{i-1}$ for all $i\in \inter{\nu}{\rho}$.\\
    \esp Then, $n_{\rho}=b$ and :

  $$ \begin{cases} \gcd(a,b)=1 \\ a\leqslant q_r \\ a<b<d \\ x_0<x_r+b \leqslant x_{-1} \end{cases} \Rightarrow \Irr\left(\frac{<a,b>}{d}\right)=I  $$ 
   $\blacklozenge$ (ii) all triplets of pairwise coprime $(a,b,d)$ such that  $\Irr\left(\frac{<a,b>}{d}\right)= I$ and $a<b<d$ can be obtained as in (i).\\   
    $\blacklozenge$ (iii) there is an infinity  of pairwise coprime $(a,b,d)$ such that  $\Irr\left(\frac{<a,b>}{d}\right)=I$, $a<b<d$ such that $\nu=-1,\rho=r+1$ in (i).
  
  \end{montheo}
  
  \begin{mademo} 
  (i) by obvious induction, we obtain :
  $$ \forall i\in \inter{\nu}{\rho}, \esp n_i=\frac{ax_i+by_i}{d} \esp \text{ where } \esp y_i=q_{i-1}$$
  \esp so $n_{\rho}=b$, for $x_{\rho}=dm-md=0$ and $y_{\rho}=q_{\rho-1}=d$.\\
   \esp Now, we suppose that $a$ and $b$ are coprime, $a\leqslant q_r$ and $a<b<d$. Then, $a,b,d$ are pairwise coprime and using Proposition \ref{prop:5} Case 3, with indices in range $\inter{\nu}{\rho}$ instead of $\inter{0}{r}$, we obtain :
   $$ \Irr\left(\frac{<a,b>}{d}\right)= \{ n_i, i\in\inter{\sigma}{s}\}$$
   \esp where $s,\sigma\in \inter{\nu}{\rho}$ such that $y_s<a\leqslant y_{s+1}$ and $x_{\sigma}< x_s+b \leqslant x_{\sigma-1}$.  \\
   \esp But, $a\leqslant q_r$ by hypothesis and $q_{r-1}<a$, because $I$ is $\N$-independant ( see  Proposition \ref{prop:6}). So $s=r$. Our hypothesis $x_0<x_r+b \leqslant x_{-1}$ gives $\sigma=0$, so we obtain result (i).\\
   
   (ii) a direct consequence of Proposition \ref{prop:5} ( case 3), with a shift on indices of $(n_i)_i$  so that $\sigma=0$.\\
   
   (iii) so, we restrict ourselves to the extension of $(d_i)_i$ and $(n_i)_i$ by one term in both directions : that is to say, we choose $\rho=r+1$ and $\nu=-1$. Then, we search $d_0$ and $d_r$ two integers such that  $d_0,d_r\geqslant 2$ and :
      $$ n_{r+1},n_{-1}\in\N^* \esp \mathbf{(1)} \esp  \esp ; \esp \gcd(n_0,n_{r+1})=1 \esp  \mathbf{(2)} \esp \esp ; \esp d>\max(n_{-1},n_{r+1}) \esp  \mathbf{(3)} \esp \esp ; \esp x_0<x_r+n_{r+1} \leqslant x_{-1} \esp  \mathbf{(4)} \esp$$
      \esp where : 
      $$n_{r+1}=d_rn_r-n_{r-1}, \esp n_{-1}=d_0n_0-n_1 \esp \text{ and } \esp x_{-1}=d \esp , \esp x_0=m$$
      \esp  $d/m$ being the reduced fraction of $\lceil d_0,d_1,\cdots,d_r\rceil$ and $q_i/p_i$ being the reduced fraction of $\lceil d_0,d_1,\cdots,d_i\rceil$. In addition, $p_{-1}=0$ and $q_{-1}=1$. So, $p_r=m$ and $q_r=d$.\\
      \esp We can simplify condition \textbf{(4)}. Indeed, we have here $x_r=1$, for $x_r=D_r$ with following notation and result :
      $$ \forall i\in\inter{\nu}{\rho}, \esp D_i=\begin{vmatrix} q_{i-1} & q_i\\ p_{i-1} & p_i \end{vmatrix} \esp \text{ then } D_{i+1}=D_i$$
      \esp yet $D_{\nu}=1$, so $x_r=1$.  Condition \textbf{(4)} can be rewritten : $m\leqslant n_{r+1}< d$.\\
      \esp Now, $d$ and $m$ depends on $d_0$ and $d_r$ : we denote $Q_i/P_i$ the reduced fraction of $\lceil d_1,\cdots,d_i\rceil$, for all $i\in\inter{1}{r}$. We can also set $P_{-1}=-1, P_0=0$ and $Q_{-1}=0, Q_0=1$. By obvious induction, we deduce :
      $$ \forall i\in\inter{-1}{r}, \esp p_i=Q_i \esp ; \esp q_i=d_0Q_i-P_i $$
      \esp So :
      $$ m= d_rQ_{r-1}-Q_{r-2} \esp ; \esp d=d_0( d_rQ_{r-1}-Q_{r-2})-(d_rP_{r-1}-P_{r-2})$$
      \esp Let us summarize : \\
      - condition \textbf{(1)} is satisfied if $d_0$ and $d_r$ are large enough.\\
      - condition \textbf{(2)} is satisfied for an infinite number of $d_0$, when $d_r$ is fixed ( for $\gcd(n_{1},n_0)=1$) and for an infinite number of $d_r$, when $d_0$ is fixed ( for $\gcd(n_{r-1},n_r)=1$), according to Lemma \ref{lem:12}.\\ 
      - condition \textbf{(3)} reduces to $d>n_{-1}$ if  condition \textbf{(4)} is satisfied. But, if $d_r$ is large enough so that : $ d_rQ_{r-1}-Q_{r-2}>n_0$, then $d>n_{-1}$, if $d_0$ is large enough ( see expressions of $d$ and $n_{-1}$ above) \\
      - condition \textbf{(4)} is equivalent to $m\leqslant  n_{r+1}<d$. But, since $(n_0,\cdots,n_r)$ is $\N$-independant, we have $Q_{r-1}<n_r$ ( see Proposition \ref{prop:6} (i)), so $m<n_{r+1}$ for $d_r$ large enough.  Now, if $d_r$ is fixed, then $n_{r+1}, P_r$ and $Q_r$ are fixed and for $d_0$ large enough, we will have $d>n_{r+1}$, because :
      $$ d= d_0Q_r-P_r $$
      
      \esp Conclusion : we first choose ( and fix) $d_r$ large enough so that $n_{r+1}>0, \esp d_rQ_{r-1}-Q_{r-2}>n_0$ and $n_{r+1}\geqslant m$. Then, we can choose an infinity of integers $d_0\geqslant 2$ such that $n_{-1}>0,\gcd(n_0,n_{r+1})=1, \\  d>n_{-1}$ and $d>n_{r+1}$. That proves (iii). 
  \end{mademo}

   \newpage

  \textbf{Examples : } again, we take our two examples : \\
   \esp First, $I=\{ 7,11,59 \}$, enumerated by the modular-convex sequence $(7,11,59)$ ( $d_1=6$). \\
  \esp We denote $k=d_2$ and $j=d_0$. Then $d=6kj-k-j, \esp  m=6k-1, \esp b=59k-11$ and $a=7j-11$.  If $k$ and $j$ are two integers bigger than $1$ and such that $d > \max(a,b) \geqslant m>0$ and $\gcd(a,b)=1$, then $\Irr(\frac{<a,b>}{d})=I$. So :
  $$ \forall k,j\in\N, \esp \begin{cases} j \geqslant 10 , \esp k\geqslant 2 \\ \gcd(7j-11,59k-11)=1 \end{cases} \Rightarrow \Irr\left( \frac{<7j-11,59k-11>}{6kj-k-j}\right)=\{ 7,11,59 \}$$ 

  \esp On the other hand,  $I=\{ 7,11,59 \}$ is enumerated by the modular-convex sequence $(11,7,59)$ ( $d_1=10$).\\
   \esp We denote $k=d_2$ and $j=d_0$. Then $d=10kj-k-j, \esp m=10k-1, \esp  b=59k-7$ and $a=11j-7$.  If $k$ and $j$ are two integers bigger than $1$ and such that $d > \max(a,b) \geqslant m>0$  and $\gcd(a,b)=1$, then $\Irr(\frac{<a,b>}{d})=I$. So :
    $$ \forall k,j\in\N, \esp \begin{cases} j \geqslant 6 , \esp k\geqslant 2 \\ \gcd(11j-7,59k-7)=1 \end{cases} \Rightarrow \Irr\left( \frac{<11j-7,59k-7>}{10kj-k-j}\right)=\{ 7,11,59 \}$$

   \esp Now, we consider $I=\{ 10,17,24 \}$, enumerated by the modular-convex sequence $(10,17,24)$ ( $d_1=2$). \\ 
      \esp We denote $k=d_2$ and $j=d_0$. Then $d=2kj-k-j, \esp m=2k-1, \esp  b=24k-17$ and $a=10j-17$.  If $k$ and $j$ are two integers bigger than $1$ and such that $d > \max(a,b) \geqslant m>0$  and $\gcd(a,b)=1$, then $\Irr(\frac{<a,b>}{d})=I$. So :
   $$ \forall k,j\in\N, \esp \begin{cases} j \geqslant 13, \esp k\geqslant 6  \\ \gcd(10j-17,24k-17)=1 \end{cases} \Rightarrow \Irr\left( \frac{<10j-17,24k-17>}{2kj-k-j}\right)=\{ 10,17,24  \}$$

\vspace*{0.5cm}   
 
 $\bullet$ Now, we end this section with a very particular case : when $I$ is an arithmetic sequence.
 
 \begin{moncoro} \label{coro:1}
    let $a$ and $k$ be two coprime integers greater than $1$ and $r\in \inter{1}{a-1}$. Then :
    $$ < a+jk, j\in \inter{0}{r}>=\frac{<a,a^2+dk>}{ar+1}$$
    \esp If $a$ and $r$ are coprime, wa have a simpler result :
    $$  < a+jk, j\in \inter{0}{r}>=\frac{<a,a+rk>}{r}$$
    
    \end{moncoro}
 \begin{mademo} first, the sequence $( a+jk)_{   j\in\inter{0}{r} }$ is $\N$-independant. Indeed, with Proposition \ref{prop:6} and its notations, we have $d_i=2$ for all $i\in\inter{0}{r-1}$, so $q_i=i+1$, by obvious induction. Hence, we obtain $q_{r-1}=r<a<a+rk$. \\
 \esp If $a $ and $r$ are coprime, then $a$ and $b$ are coprime, when we set $b=a+rk$. Moreover, if we denote $d=r$ and $n_1=a+k$, then $d\in\inter{1}{a-1}$ and :
 $$ dn_1=ar+rk=b \text{ mod } a$$
 \esp So, we can apply Theorem \ref{theo:1} and obtain the result in that case. \\
 \esp In the general case, we use Theorem \ref{theo:2} (i) and its notations, with $\rho=r+1$ and $d_r=a+1$. We obtain $d=ar+1$ and $b=a^2+dk$. In addition, $a$ and $b$ are coprime and $a\leqslant q_r=d<b$.
 \end{mademo}

 \textbf{Remark : } as mentioned in the introduction, this result can be proven directly...

 \newpage

  \section{Diophantine approximation}
  
     \subsection{two kinds of continued fraction expansions}
    
       \label{subsec:CFE}
            
        \esp All  results given in this subsection are well known and we just want to recall some notations and simple facts. For all reals $x$, $\lfloor x \rfloor$ denotes its floor ,$\lceil x \rceil $ its ceiling and $\{ x \}$ its fractional part. We have $x=\lfloor x \rfloor + \{ x \}$. \\
     
     $\bullet$ We have met in \ref{subsec:modular} the notion of " ceiling continued fraction expansion of a rational". We will precise the definition and remind some results about usual continued fraction expansion. \\
     \esp One usually uses the Gauss map $T$, defined on $]0,1[$, by $T(x)=\{ 1/x \}$ to define the usual continued fraction expansion of an irrational. We will rather use here a variant : \\
      \esp Let $F$ and $C$ be the following maps defined on $\R$ :
     $$ F : \begin{cases} \R\cup \{ +\infty\} \to ]1,+\infty] \\ x\to +\infty \text{ if } x\in\Z\cup  \{ +\infty\}  \\ x \to \frac{1}{\{ x \}} \text{ else }\end{cases} \esp ; \esp \esp  C : \begin{cases} \R\cup \{ +\infty\} \to ]1,+\infty] \\ x\to +\infty \text{ if } x\in\Z\cup  \{ +\infty\}  \\ x \to \frac{1}{1-\{ x \}} \text{ else }\end{cases}$$  
     \esp We use the natural convention : $\lfloor +\infty \rfloor = \lceil +\infty \rceil=+\infty$. By obvious induction, we can prove that : if $\alpha$ is not rational, then $F^k(\alpha)\not = +\infty$ and $C^k(\alpha)\not = +\infty$ for every $k\in\N$. It is a consequence of euclidean algorithm on integers that : if $\alpha$ is rational, then $F^k(\alpha)=C^k(\alpha)=+\infty$ for some positive integer $k$.  \\
     \esp Now, given a real $\alpha$, we define two sequences ( named \emph{partial quotients} ) by : 
     $$ \forall k\in\N, \esp a_k=\lfloor F^k(\alpha) \rfloor \esp ; \esp d_k= \lceil C^k(\alpha) \rceil $$
     \esp If, $\alpha$ is an irrational, then $a_k$ and $d_k$ are positive integers for all $k\in\N$ and $d_k\geqslant 2$ if $k>0$. \\
     \esp If, $\alpha$ is a rational, then $a_k=d_k=+\infty$ for $k$ large enough. \\
     \esp $[a_k]_{k\in\N}$ is the \emph{usual continued fraction expansion} of $x$ ( we could also name it the " floor continued fraction expansion" of $x$) and $\lceil d_k \rceil _{k\in\N}$ is the " ceiling continued fraction expansion" of $x$.  \\
     \esp We usually define the \emph{( usual or ceiling) convergents } of $\alpha$ relative to these continued fraction expansions by : 
     $$ \forall i\in \N, \esp \frac{p_i}{q_i}=[a_0,a_1,\cdots,a_i] \esp ; \esp \frac{p'_i}{q'_i}=\lceil d_0,d_1,\cdots,d_i\rceil $$
     \esp if necessary, we will precise $ \frac{p_i(\alpha)}{q_i(\alpha)}$. The above brackets are inductively defined by  ( with the usual convention : $1/\infty =0$) :
     $$ [a_0,a_1,\cdots,a_i]=a_0+\frac{1}{[a_1,\cdots,a_i]} \esp ;\esp \lceil d_0,d_1,\cdots,d_i\rceil =d_0-\frac{1}{\lceil d_1,\cdots,d_i\rceil }$$
     \esp Thus, we can express these  with finite fraction superpositions : 
    $$  [a_0,a_1,\cdots,a_i]=a_0+\dfrac{1}{a_1+\dfrac{1}{a_2+\dfrac{1}{a_3+\cdots }}} \esp ; \esp \lceil d_0,d_1,\cdots,d_i\rceil = d_0-\dfrac{1}{d_1-\dfrac{1}{d_2-\dfrac{1}{d_3-\cdots }}}$$  
    
     \esp The sequences of ( usual or ceiling) convergents of $\alpha$ converge towards $\alpha$. \\
     \esp For example : if $\alpha= \frac{1+\sqrt{5}}{2}$ ( the so-called " golden ratio"), then : $a_k=1$ for all $k\in\N$, $d_0=2$ and $d_k=3$ for all $k\in\N^*$. Indeed : for these expansions we have :
     $$ x=1+\frac{1}{x} \text{ and } x=2-\frac{1}{1+x} $$

   \esp   With these parallel definitions, we could expect similar properties...but, this is not the case :

   \newpage
   
   \esp First, a useful Lemma that relates iterates of $F$ and $C$ : 
   \begin{monlem} \label{lem:13} let $\alpha$ be a real, but not an integer. We denote $a=\lfloor F(\alpha) \rfloor$.
   $$ \forall i\in\inter{1}{a}, \esp C^i(\alpha)= 1+ \frac{1}{F(\alpha)-i}  $$
   \esp So :
   $$ C^{a}(\alpha)=1+F^2(\alpha)$$
    \end{monlem}
   \begin{mademo} Second result is a direct consequence of the first result, which is obtained by finite induction on $i$ : \\
   - it is true for $i=1$, since, if we denote $\epsilon=\{ \alpha \}$, then : 
   $$ C(\alpha)= \frac{1}{1-\epsilon} \esp \text{ and } \esp
   1+\frac{1}{F(\alpha)-1}=1+ \frac{1}{\frac{1}{\epsilon}-1}=1+ \frac{\epsilon}{1-\epsilon}=  \frac{1}{1-\epsilon} $$
   - suppose it is true for $i-1$, with $i\in\inter{2}{a}$, then we remark that $C^{i-1}(\alpha)\in ]1,2]$ for $F(\alpha)\geqslant i$. The special case when $i=a=F(\alpha)$ leads to $C^{i-1}(\alpha)=2$ and so $C^i(\alpha)=+\infty=1+\frac{1}{0}$ ( with usual convention). Now, we suppose that $F(\alpha)>a$, so that $C^{i-1}(\alpha)\in ]1,2[$. Then :
   $$ C^i(\alpha)=\frac{1}{1-\{ C^{i-1}(\alpha)\}}=\frac{1}{2-C^{i-1}(\alpha)}=\frac{F(\alpha)-(i-1)}{F(\alpha)-i}=1+ \frac{1}{F(\alpha)-i}  $$

   \end{mademo}

   \esp The following Lemma states the relations between the two kinds of partial quotients seen above, for the case of  irrationals. The case of rationals will be treated later, for we will choose an alternate version of usual " last partial quotient" for a rational.
     
     \begin{monlem} \label{lem:14} let $\alpha$ be an irrational and $[a_k]_{k\in\N},\lceil d_k \rceil_{k\in\N}$ its usual and ceiling continued fraction expansions.\\
   (i) $\lceil d_k \rceil_{k\in\N^*}$  is obtained from $[a_k]_{k\geqslant 1}$ by the following substitution of patterns :
   $$ \forall i\in\N^*, \esp (a_{2i-1},a_{2i}) \to ((2,)^{a_{2i-1}-1},a_{2i}+2)$$   
  (ii) for the converse : we denote $(k_i)_{i\in\N}$ an increasing sequence of non negative integers such that $k_0=0$ and $(d_{k_i})_{i\in\N^*}$ is the subsequence of values of $(d_k)_{k\in\N^*}$ such that  $d_k\not = 2$. Then : 
      $$ \forall i\in \N^*,\esp  \begin{cases} a_{2i-1}=k_i-k_{i-1} \\  a_{2i}=d_{k_i}-2 \end{cases} $$
      
    \end{monlem}

    \begin{mademo} 
    (ii) is a direct consequence of (i). For (i), we use previous Lemma : we will explain the process for the first step, namely $i=1$. First, we recall that $C^i(\alpha)$ and $F^j(\alpha)$ are never integer since $\alpha$ is not rational.\\
    \esp  Then, as mentioned in the proof of Lemma \ref{lem:13}, $C^i(\alpha)\in ]1,2[$ for all $i\in \inter{1}{a_1-1}$, so $d_i=2$ for these values of $i$. \\
    \esp Later, $C^{a_1}(\alpha)=1+F^2(\alpha)$ and $F^2(\alpha)$ is not an integer, so :
    $$ \lceil C^{a_1}(\alpha)\rceil = 2+\lfloor F^2(\alpha) \rfloor $$
    \esp which proves that $d_{a_1}=a_2+2$. Finally :
    $$ \{ C^{a_1}(\alpha) \}=\{ F^2(\alpha) \}$$
    \esp so, the process continues at the next step...   
    \end{mademo}

   
   
  $\bullet$   There are two ways for usual continued fraction expansion  of rationals : those ending with an integer greater than $1$ ( and an infinite sequence of $\infty$) or those ending with $1$ ( and an infinite sequence of $\infty$).\\
   \esp For example, $5/7$ has two usual continued fraction expansions  : $[0,1,2,2]$ and $[0,1,2,1,1]$ ( the ending $\infty$ sequence has been omitted). We will choose the second way and denote CFE this kind of usual continued fraction expansion.

   
\esp  The set of sequences representing CFE of this type, for reals,  is then : 
   
    $$ \CFE= \{ (t_k)\in \Z\times \N^* \times (\overline{\N^*})^{\N}, \forall k\geqslant 2, ( t_k= \infty \Rightarrow ( t_{k+1}=\infty \text{ and } t_{k-1}\in \{ \infty  ,1 \} ) \} $$

 \esp Since we will be interested by natural order on reals and especially on rationals, we mention the corresponding order on CFE sequences :   if we define $\phi$ by 
   $$ \phi : \begin{cases} \CFE \to \R \\ (t_k)_{k\in\N} \to [t_k]_{k\in\N} \end{cases}$$
   \esp this map is bijective and increasing, with the Alternate Lexicographic Order ( ALO) on $\CFE$ defined by : 
   $$ (t_k)_{k\in\N}\leqslant_A (t'_k)_{k\in\N} \Leftrightarrow  (\forall k\in\N, t_k=t'_k ) \text{ or }  \exists j\in\N, \begin{cases} \forall k<j, t_k=t'_k  \\ (-1)^jt_j< (-1)^j t'_j \end{cases}$$

$\bullet$ Now is the time to give a " rational version" of Lemma \ref{lem:14} : 
 \begin{monlem} \label{lem:15} let $\alpha$ be a rational and $[a_0,\cdots,a_r,1], \lceil d_0,d_1,\cdots,d_{\rho}\rceil$ its usual and ceiling continued fraction expansions. ( $\infty$-tails  are not mentioned) \\
 (i) $r$ and $\rho$ are related by the following formula :
 $$  \rho = \sum_{i=1}^{\lfloor (r+1)/2 \rfloor }  a_{2i-1}$$
 (ii)    $\lceil d_k \rceil_{k\in\inter{1}{\rho}}$  is obtained from $[a_k]_{k\in \inter{1}{r}}$ by the following substitution of patterns :
   $$ \forall i\in\inter{1}{\lfloor r/2 \rfloor}, \esp (a_{2i-1},a_{2i}) \to ((2,)^{a_{2i-1}-1},a_{2i}+2)$$ 
   \esp if $r$ is odd, then we add $d_{\rho}=2$.\\
   
   (iii) for the converse : let $(k_i)_{i\in\inter{0}{s}}$ be an increasing sequence of non negative integers such that $k_0=0, k_{s}=\rho$ and $(d_{k_i})_{i\in\inter{1}{s -1}}$ is the subsequence of values of $(d_k)_{k\in\inter{1}{\rho-1}}$ such that $d_k\not = 2$. Then $s=\lfloor (r+1)/2 \rfloor$ and :
     $$ \forall i\in \inter{1}{s }, \esp \begin{cases} a_{2i-1}=k_i-k_{i-1} \\  a_{2i}=d_{k_i}-2 \esp ( \text{ if } 2i\leqslant r )\end{cases} $$
  
\end{monlem} 
 \begin{mademo} (i) and (ii) are direct consequences of (iii), that we prove as in the proof of Lemma \ref{lem:14}, adding the following fact. If $F(\alpha)$ is an integer, then $r=1,a_1=F(\alpha)-1$ and $C^{a_1}(\alpha)=2$ ( see Lemma \ref{lem:13}, with $a=a_1+1$). So, $d_{a_1}=2,\rho=a_1,s=1$ and $a_1=k_1-k_0$. If $F^2(\alpha)$ is an integer, then $a_1=\lfloor F(\alpha)\rfloor,a_2=F^2(\alpha)-1$ and $C^{a_1}(\alpha)=a_2+2$ ( see Lemma \ref{lem:13}). So, $r=2,d_{a_1}=a_2+2, \rho=a_1$ and $s=1$.\\
 \esp As in the proof of Lemma \ref{lem:14}, this can be generalized to the case : $F^j(\alpha)$ is an integer for a positive integer $j$.
  \end{mademo}
 \newpage

   \subsection{semi-convergents and best rational in an interval}  
   
  \label{subsec:semicv}

     $\bullet$ Let $\alpha$ be a real with CFE $[a_k]_{k\in\N}$ and $(p_k/q_k)_k$ its convergents sequence.\\
     \esp  A \emph{semi-convergent} of $\alpha$ is any rational of the form $\frac{mp_k+p_{k-1}}{mq_k+q_{k-1}}$, with $m\in \inter{0}{a_k}$ and $k\in \N$ such that $a_k<\infty$ ( we take $m>0$ if $k=0$ to avoid $1/0$ !). So, convergents are particular semi-convergents.

   \begin{monlem} \label{lem:16} Let $\alpha$ be a real with CFE $[a_k]_{k\in\N}$. Semi-convergents of $\alpha$ are exactly the rationals with CFE $[a_0,\cdots,a_{s-1},b_s,1]$, such that $s\in\N,b_s\in\inter{1}{a_s}$ and $a_{s+1}<\infty$.           
    \end{monlem}  
    
     \textbf{Remark : } it would not be difficult to deduce that the  ceiling convergents of $\alpha$ are exactly the  semi-convergents of $\alpha$ that are greater or equal than $\alpha$. ( see Lemma \ref{lem:15} and order on CFE) \\

     \begin{mademo}
    Consequence of the definition and the well known fact : 
    $ \forall m\in \N^*, [a_0,\cdots,a_{s-1},m]=\frac{mp_{s-1}+p_{s-2}}{mq_{s-1}+q_{s-2}}$.    
     \end{mademo}

       
       
             
  $\bullet$ Now, we would like to precise the CFE of reals in $\overset{\longleftrightarrow}{[\theta, \theta']}$ ( denotes the set of reals that are between $\theta$ and $\theta'$, even if $\theta > \theta'$), where $\theta$ and $\theta'$ are two different reals and find the rationals in this interval with the lowest reduced denominator.   : 
 
   
   \begin{madef} \label{def:3} let $x$ be a real. We name \emph{ CFE-depth} of $x$ the non negative integer, denoted $\mu(x)$ and defined by :   $\mu(x)=+ \infty$ if $x$ is  irrational  and $\mu(x)=s$, if $x=[a_0,a_1,\cdots,a_s,1]$ is the CFE of $x$.\end{madef} 
   
   \esp We  remark that : 
   $$ \mu(x)=0 \Leftrightarrow x\in\Z \esp ; \esp \forall n\in\Z, \mu(x+n)=\mu(x) \esp ; \esp \forall x\not\in \Z, \mu(T(x))=\mu(x)-1$$
   
    \esp We denote $\theta = [t_k]_{k\in\N}$ and       $\theta' = [t'_k]_{k\in\N}$ and  will abreviate $t$ and $t'$ these CFE-sequences. We denote $r$ the smallest integer $k$ such that $t_k\not = t'_k$. If we suppose that $t_r<t'_r$, then we have $r\leqslant \min(\mu(\theta),\mu(\theta'))+2$, when $\theta$ or $\theta'$ is rational ( if they are both irrationals, $r$ is finite ! ). Indeed, the extremal case when $r=\mu(\theta)+2$ for example corresponds to $\theta=[t_0,\cdots,t_{r-2},1]$ and $\theta'=[t_0,\cdots,t_{r-2},1,t'_r,...]$, with $t'_r<\infty$.

  
   \esp We remark that, all integers in $\overset{\longleftrightarrow}{[\theta, \theta']}$ minimize the denominator of their reduced fraction : it is $1$ !! So, we can suppose that $\lfloor \theta \rfloor = \lfloor \theta' \rfloor$ and even that $\theta, \theta'\in [0,1[$.\\
   \esp The following result claims that, in that case, there is only one rational in $\overset{\longleftrightarrow}{[\theta, \theta']}$, that minimizes the value of its denominator : it is usually named  the " \emph{best rational}" in $\overset{\longleftrightarrow}{[\theta, \theta']}$ 
   
   \begin{maprop} \label{prop:7} let $\theta$ and $\theta'$ be two different reals in $[0,1[$ and       $\theta = [t_k]_{k\in\N},\theta' = [t'_k]_{k\in\N}$  their respective CFE. We denote $r$ the lowest integer $k$ such that $t_k\not = t'_k$. \\
   \textbf{(i)} there is a unique rational in  $\overset{\longleftrightarrow}{[\theta, \theta']}$ that minimizes the denominator. We denote it $\gamma$.\\
   - if $r\leqslant \min(\mu(\theta),\mu(\theta'))$, then $\gamma= [t_0,\cdots,t_{r-1},\min(t_r,t'_r),1]$.\\
   - else, $\mu(\theta)<\mu(\theta')$ ( up to swap) and $\gamma=\theta$.\\
  \textbf{(ii)} in both cases, $\mu(\gamma)\leqslant \min(\mu(\theta),\mu(\theta'))$ and $\gamma = [t_0,\cdots,t_{s-1},\min(t_s,t'_s),1]$, where $s=\mu(\gamma)\leqslant r$, so $\forall k\in\inter{0}{s-1},t_k=t'_k$.\\
 \textbf{(iii)} the best rational in  $\overset{\longleftrightarrow}{[\theta, \theta']}$ is the common semi-convergent of $\theta$ and $\theta'$ with the greatest denominator.
   
  
   \end{maprop}    

   \begin{mademo} see \textbf{[3]}( \nameref{biblio}) 1.4 Proposition 1.  
   \end{mademo}
  
  \textbf{Remark : } as a direct consequence of (iii) : $\theta$ is the best rational in $\overset{\longleftrightarrow}{[\theta, \theta']} $ if and only if $\theta$ is a semi-convergent of $\theta'$. \\

  $\bullet$  Let $\alpha$ be a real, $[a_k]_{k\in \N^*}$ its CFE and $r=\mu(\alpha)$, the CFE-depth of $\alpha$. So, we denote $[a_0,a_1,\cdots,a_r,1]$ the CFE of $\alpha$ if $\alpha$ is rational. We also denote $(p_n/q_n)_n$ the usual sequence of convergents of $\alpha$.\\
 \esp We consider the usual notion of best rational approximation of a real $\alpha$ : for $p,q$ two integers, $p/q$ is said a \emph{best rational approximation} of $\alpha$ if and only if :
 $$ \forall q'\in \inter{1}{q}\esp , \esp \forall p'\in\Z \esp , \esp  \left|\dfrac{p'}{q'}-\alpha \right| \geqslant \left|\dfrac{p}{q}-\alpha \right|$$
 \esp It is well known that best rational approximation of a real are exactly its reduced convergents. \\
 \esp Now, we can consider two sided similar definitions :  for $p,q$ two integers, $p/q$ is said a \emph{best left rational approximation} of $\alpha$ if and only if :
 $$ \forall q'\in \inter{1}{q}\esp , \esp \forall p'\in\Z \esp , \esp  \dfrac{p'}{q'} \leqslant  \dfrac{p}{q} \leqslant \alpha \text{ or }  \dfrac{p'}{q'}>\alpha$$
  \esp $p/q$ is said a \emph{best right rational approximation} of $\alpha$ if and only if :
 $$ \forall q'\in \inter{1}{q}\esp , \esp \forall p'\in\Z \esp , \esp  \dfrac{p'}{q'} \geqslant  \dfrac{p}{q} \geqslant \alpha \text{ or }  \dfrac{p'}{q'}<\alpha$$
 
 \esp Here is a corollary of Proposition \ref{prop:7} :
       
  \begin{moncoro} \label{coro:2} .\\
  (i) best left rational approximations of $\alpha$ are the semi-convergents of $\alpha$, that are lower than $\alpha$.\\
   (ii) best right rational approximations of $\alpha$ are the semi-convergents of $\alpha$, that are greater than $\alpha$. 
   \end{moncoro}
   \begin{mademo}(i) we remark that $p/q$ is a best left rational approximation of $\alpha$ if and only if $p/q$ is the best rational in $[p/q, \alpha]$ and use the remark below Proposition \ref{prop:7}. Same arguments for (ii).   
   \end{mademo}    
       
    \esp If  we denote  $(p_k/q_k)_k$ the reduced convergents of $\alpha$, then : \\
    - its best left rational approximations are  : 
    $$\frac{p_{2i}+mp_{2i+1}}{q_{2i}+mq_{2i+1}} \esp ; \esp i\in \inter{0}{(\mu(\alpha)-1)/2} \esp ; \esp m\in \inter{0}{a_{2i+2}}$$
     - its best right rational approximations are  : 
    $$\frac{p_{2i-1}+mp_{2i}}{q_{2i-1}+mq_{2i}} \esp ; \esp i\in \inter{1}{\mu(\alpha)/2} \esp ; \esp m\in \inter{0}{a_{2i+1}}$$

   \begin{maprop}  \label{prop:8} .\\
   (i) let $\alpha$ and $\alpha'$ be two reals such that $0<\alpha<\alpha'<1$. We denote $\gamma$ the best rational in $]\alpha, \alpha']$ and $q$ the denominator of its reduced fraction. Then    $$ q = \max \{ N\in\N, \forall n\in \inter{0}{N-1}, \lfloor n\alpha \rfloor = \lfloor n\alpha' \rfloor \}$$  
  (ii) let $\alpha$ be a real in $[0,1)$ and $p/q$ a reduced fraction, with $q\in\N^*$, such that $\alpha$ is not the best left  strict convergent of $p/q$. 
   $$ p/q \text{ is a semi-convergent of } \alpha \Leftrightarrow  \forall k\in \inter{0}{q-1}, \lfloor k\alpha \rfloor = \lfloor kp/q \rfloor$$ 
   
     \end{maprop}
    \begin{mademo} see proof of Proposition 8 in 4.2 of \textbf{[3]} ( \nameref{biblio}) \end{mademo}  
    

\newpage

      \subsection{$\alpha$-numeration for a rational $\alpha$}

  \label{subsec:alphanum}
     
   $\bullet$ Let $\alpha$ be a rational in $[0,1[$ and $\alpha=[0,a_1,\cdots,a_r,1]$ its CFE. We will denote $(p_k/q_k)_{0\leqslant k \leqslant r }$ its convergents, so that $\alpha= p_{r+1}/q_{r+1}$.
   
   \begin{madef} \label{def:4} a sequence $d$ in $\N^r$ is said \emph{$\alpha$-admissible} if and only if :
   $$\forall j\in\inter{1}{r},\begin{cases}  d_j\in\inter{0}{a_j}\\  d_j=0 \Rightarrow ( \forall i\geqslant j, d_i=0 )\text{ or } d_{j-1}=a_{j-1} \end{cases}$$
   \end{madef}
   
    \esp We will denote $E_{\alpha}$ the set of $\alpha$-admissible sequences. \\
   
   \textbf{Remark : } for $j=1$, the second condition reduces to $d_1=0\Rightarrow \forall i\geqslant 1,d_i=0$. So to say, $d=(0,\cdots,0)$ is the only element of $E_{\alpha}$ , whose first coordinate is $0$. \\

    \esp We consider  the \textbf{reversed lexicographic order ( RLO) } denoted $ \underset{R}{\leqslant}$ on $\N^r$ : 
     $$ d \underset{R}{\leqslant} d' \Leftrightarrow d=d' \text{ or } \exists j\in \inter{1}{r}, \begin{cases} d_j<d'_j \\
     \forall i\in \inter{j+1}{r}, d_i=d'_i \end{cases} $$
    \esp  It is a  total order on $E_{\alpha}$.
       
     \begin{maprop} \label{prop:9}
       the map $\Psi_{\alpha}$ is an order isomorphism from $(E_{\alpha} ,\leqslant_R)$ to $\inter{0}{q_{r+1}-1},\leqslant)$.
     $$ \Psi_{\alpha} : \begin{cases} E_{\alpha} \to \inter{0}{q_{r+1}-1} \\ d\to \sum\limits_{j=1}^rd_jq_{j-1} \end{cases}$$
      \end{maprop}
      
      \begin{mademo} see \textbf{[3]} 2.2. ( \nameref{biblio})   
      \end{mademo}

      \esp The following algorithm explains the inverse function of $\Psi_{\alpha}$. We will denote $m_k=q_k+q_{k-1}$ for any $k\in\inter{0}{r}$. So $m_r=q_{r+1}$.

  \begin{monalgo} \label{algo:1} let $n\in \inter{0}{m_r-1}$. \\
  With the following algorithm, we have  $d\in E_{\alpha}$ and $\Psi_{\alpha}(d)=n$.
  
   \shadowbox{\parbox{11cm}{\textbf{Input : } $n$ \esp \textbf{Output : } $(d_i)_{i\in \inter{1}{r}}$ \\ 
   for $k=r$ to $k=1$ with  step $-1$ : $\begin{cases} d_k=\max\left(0,\left\lfloor \frac{n-q_{k-2}}{q_{k-1}}\right\rfloor \right)\\ n \leftarrow n-d_kq_{k-1} \end{cases}$}}

   \end{monalgo}


  \textbf{Remark : } as a direct consequence : $E_{\alpha}$ has $q_{r+1}$ elements.\\

     $\bullet$ Now, we will deal with $\alpha$-numeration for elements of  $U_{\alpha}=  \{ \{ k\alpha\}, k\in\N \} $. Since, $\alpha=\frac{p_{r+1}}{q_{r+1}}$ and this fraction is reduced, we have $U_{\alpha}= \{ \frac{n}{q_{r+1}}, n\in \inter{0}{q_{r+1}-1}\}$. So, this set is very simple, but we will focus on the map $k\to \{ k\alpha \}$, with the order point of view :\\
     
     \esp We consider the \textbf{alternate lexicographic order ( ALO)} denoted $ \underset{A}{\leqslant}$ on $\N^r$ : 
     $$ d \underset{A}{\leqslant}d' \Leftrightarrow d=d' \text{ or } \exists j\in \inter{1}{r}, \begin{cases} (-1)^{j-1}d_j<(-1)^{j-1}d'_j \\
     \forall i\in \inter{1}{j-1}, d_i=d'_i \end{cases} $$
     \esp It is another total order on $E_{\alpha}$.\\
     
     \newpage

      We define also  : 
      $$\forall i\in\inter{-2}{r}, \esp \delta_i=(-1)^i( q_{i}\alpha - p_i)$$
      \esp We have, with $a_0=0$ here :
      $$ \delta_{-2}=\alpha \esp ; \esp \delta_{-1}=1 \esp ; \esp \delta_0=\{ \alpha \} =\alpha \esp ; \esp \forall i\in\inter{0}{r},  \esp \delta_i=-a_i\delta_{i-1}+\delta_{i-2}$$
      \esp Let $T$ be the Gauss map : $]0,1[\to [0,1[, x\to \{ 1/x \}$.\\
      \esp  By induction on $i$, with the fact that : $a_i=\left\lfloor \frac{1}{T^{i-1}(\alpha)} \right\rfloor$ if $i\leqslant r-1$, we obtain :
      $$ \forall i\in\inter{0}{r-1}, \frac{\delta_i}{\delta_{i-1}}=T^i(\alpha) $$
      
      \esp  Beware : for $i=r$, $T^{r-1}(\alpha)=[0,a_r,1]=\frac{1}{a_r+1}$, so :
      $$ \frac{\delta_r}{\delta_{r-1}}= \frac{\delta_{r-2}-a_r\delta_{r-1}}{\delta_{r-1}}= \frac{1}{T^{r-1}(\alpha)}-a_r=1 $$
      \esp So : $\delta_r=\delta_{r-1}$. We also have ( see [...])   $ \delta_r=\delta_{r-1}=\frac{1}{q_{r+1}}$.\\
      \esp To summarize this :
      $$ \forall i\in \inter{0}{r-1}, 0<\delta_i<\delta_{i-1} \esp ; \esp \delta_r=\delta_{r-1}=\frac{1}{q_{r+1}}$$

     \begin{maprop} \label{prop:10} .\\
      (i) the map $\Lambda_{\alpha} $ ( defined below) is an order isomorphism, with ALO on $E_{\alpha}$ : 
       $$ \Lambda_{\alpha} : \begin{cases} E_{\alpha} \to \left\{ \frac{n}{q_{r+1}}, n\in \inter{0}{q_{r+1}-1}\right\} \\ d\to \sum\limits_{j=1}^rd_j(-1)^{j-1}\delta_{j-1} \end{cases}$$
       
     (ii) we have : 
     $$ \forall n\in \inter{0}{q_{r+1}-1}, \{ n\alpha \}= \Lambda_{\alpha}(\Psi_{\alpha}^{-1}(n))$$

     \end{maprop}

  \begin{mademo} see \textbf{[3]} 2.2 ( \nameref{biblio}).   
 \end{mademo}

    \textbf{Remarks : } result (ii) means that the map $n\to \{n\alpha \}$ ( with $0\leqslant n < q_{r+1}$), is, from the order point of view, the " same thing"  as the identity $(E_{\alpha},RLO)\to (E_{\alpha},ALO)$. \\
     
     \esp We can sum up these formulae  : $ \forall n\in \inter{0}{q_{r+1}-1}$, with $d=\Psi_{\alpha}^{-1}(n)$ :
     $$n=\sum_{j=1}^r d_jq_{j-1} \esp ;\esp   \lfloor n\alpha\rfloor =\sum_{j=1}^r d_jp_{j-1}  \esp ; \esp \{ n\alpha \}=\sum_{j=1}^r (-1)^{j-1}d_j\delta_{j-1}$$

  \esp   The following algorithm expresses the inverse function of $\Lambda_{\alpha}$.

  \begin{monalgo} \label{algo:2} let $\beta\in \{ n/q_{r+1}, n\in\N \}$.\\ Applying the algorithm below ( with $\beta_0=\beta$), we have :\\
  (i) $b\in E_{\alpha}$.\\
  (ii)  $\beta = \Lambda_{\alpha}(b)$.
  
   \shadowbox{\parbox{11cm}{\textbf{Input : } $\beta$ \esp \textbf{Output : } $(b_i)_{i\in \N^*}$ \\ 
  for $k=1$ to $k=\infty$ with step $1$ : $\begin{cases} b_k=\min\left(a_k,\left\lceil \frac{\beta_{k-1}}{\delta_{k-1}}\right\rceil\right) \\ \beta_{k}=b_k\delta_{k-1}-\beta_{k-1} \end{cases} $}}
   \end{monalgo}
   
 

 \newpage

  

    \subsection{minimal points of $(\{ n\alpha - \beta \},n )_n$}

    \esp In this subsection, we prove several Lemma about minimal points of this sequence, in different cases, from the simplest to the more intricate. These results will be useful, when we will study in section 5, $\Irr$ and $\PF$ of $\frac{<a,b>}{d}$, since points of $L$ in a rectangle of vertice $d$ can be parametrized by such sequences ( see \ref{subsec:relations}).\\ 
  
  $\bullet$ First, we study  minimal points in $\R^2$, with the product order, of the sequence $(\{ n\alpha \},n)_{1\leqslant n \leqslant N }$, where $N$ is a positive integer. \\
  \esp With Proposition \ref{prop:10}, we can affirm that these points have the property that if $n=\Psi_{\alpha}(d)$, then for any other $d'\in E_{\alpha}$, we have $d<d'$ for ALO or for RLO ! \\
  \esp The following Lemma gives the result for the particular case $N=q-1$, where $q$ is the denominator of the reduced fraction of $\alpha$. We denote, as usual, $[a_0,a_1,\cdots,a_r,1]$ the CFE of $\alpha$. \\

  \begin{monlem} \label{lem:17} the minimal points in $\R^2$ of the set $\{ (\{ n\alpha \},n) , 1\leqslant n \leqslant q-1 \}$ are for the following integers $n$, given by three equivalent formulations : \\
  (i) with the $\alpha$-numeration of $n$ ( the $0$-tail is not written) :
  $$(1)\esp ; \esp (1,j) , j\in \inter{1}{a_2} \esp  ; \esp  (1,a_2,0,j), j\in\inter{1}{a_4} ; \cdots ; (1,a_2,0,\cdots,0,j) , j\in\inter{1}{a_{2\lfloor r/2 \rfloor}} $$
 we can sum up this with : $(1)$ and sequences $ (d_k)_{k\in\inter{1}{r}}$ such that  $d_1=1$ and :
 $$   \forall i\in\inter{1}{k-1}, d_{2i}=a_{2i}, d_{2i+1} =0, \esp d_{2k}\in \inter{1}{a_{2k}}, \forall i> 2k, d_i=0, \esp k\in\inter{1}{\lfloor r/2 \rfloor}$$ 
  (ii) with denominators of convergents of $\alpha$ :
  $$ q_0 \esp ; \esp   q_{2k-2}+jq_{2k-1}, j\in\inter{1}{a_{2k}}, k\in \inter{1}{\lfloor r/2 \rfloor} $$
  
  (iii) with diophantine approximation : these are the denominators of the reduced semiconvergents of $\alpha$, that are lower than $ \alpha$.
  \end{monlem}
 \begin{mademo}
   - let us verify that these points are minimal : for $(1)$, it is obvious. We denote $A_{p}=(a_2,0,a_4,0,\cdots,a_{2p-2},0)$, for any positive integer $p$ ( with $A_1=()$, the empty sequence).\\
   \esp Let $d=(1,A_p,j)$ for $j\in \inter{1}{a_{2p}}$ and $p\leqslant r/2$. If $d'<d$ for ALO and $d'\not = 0$, then $d'=(1,A_p,j')$ with $j'>j$, or $d'=(1,A_p,j,0,k)$, with $k\in\inter{1}{a_{2p+2}}$. In both cases, we have $d'>d$ for RLO. \\
   
   - now, let us verify that no other point is minimal : if $d'$ is not among these $\alpha$-numeration of integers $n$, then let $k$ be the least integer such that $d'_k$ is different from all the $d_k$ of (i). Then, we claim that $k$ is odd. Indeed, if $k$ were even, we would have $d'=(1,A_p,d'_k,...)$, with $k=2p+2$ and $d'_k=0$. So $d'=(1,A_p)$, according to our definition of $\alpha$-admissible sequences,  but this is in the list (i). \\
   \esp So, $k$ is odd and if we denote $d=(1,A_p,d_{k-1})$, with $k=2p+1$, then $d<d'$ for ALO and RLO.\\
   
  - (ii) and (iii) are direct consequences of (i) 
   \end{mademo}

  \newpage

 $\bullet$ We will also need a more general case : but it is sufficient to sort among the former values of $n$, those which are lower than $N$ :
  
   \begin{monlem}   \label{lem:18} let $N\in\inter{1}{q-1}$. We denote $(N_1,N_2,\cdots,N_s)$ the $\alpha$-numeration of $N$, with $N_s>0$. The minimal points in $\R^2$ of the set $\{ (\{ n\alpha \},n) , 1\leqslant n \leqslant N\}$ are for the following values of $n$ : \\
   
 (i) with the $\alpha$-numeration of $n$ ( the $0$-tail is not written) :
  $$(1)\esp ; \esp (1,j) , j\in \inter{1}{a_2} \esp  ; \esp  (1,a_2,0,j), j\in\inter{1}{a_4} ; \cdots ; (1,a_2,0,\cdots,a_{2\lfloor s/2 \rfloor-2},0,j) , j\in\inter{1}{\mu} $$
  where $\mu=N_s$ if $s$ is even, $\mu=a_{s-1}$ if $s$ is odd.\\
 we can sum up this with : $(1)$ and $(d_k)_{k\in\inter{1}{r}} $, such that $ d_1=1$ and :
 $$  \forall i\in\inter{1}{k-1}, d_{2i}=a_{2i}, d_{2i+1} =0 \esp , \esp d_{2k}\in \inter{1}{a'_{2k}}, \forall i> 2k, d_i=0, \esp k\in\inter{1}{\lfloor s/2 \rfloor}$$ 
 where $a'_{2k}=a_{2k}$ if $2k<s$, and $a'_s=N_s$.\\
 
  (ii)  with denominators of convergents of $\alpha$  :
  $$ q_0 \esp ; \esp   q_{2k-2}+jq_{2k-1}, j\in\inter{1}{a'_{2k}}, k\in \inter{1}{\lfloor s/2 \rfloor} $$
   where $a'_{2k}=a_{2k}$ if $2k<s$, and $a'_s=N_s$.

  
  \end{monlem}

  \begin{mademo}
  It is a consequence of the former Lemma, when we remember ( see Proposition \ref{prop:9})  that $\Psi_{\alpha}$ is an increasing map from $E_{\alpha}$ to $\inter{0}{q-1}$, with RLO. So, we just have to select, among the $d$ of the previous results those which are lower than $(N_1,N_2,\cdots,N_s)$ for RLO. Indeed, no $(\{ k\alpha \},k)$, with $k>N$ is lower than any $(\{ n\alpha \},n)$, with $n\leqslant N$...
 \end{mademo}

  $\bullet$ We generalize a bit more again , adding another condition :
  
    \begin{monlem} \label{lem:19} let $N\in\inter{1}{q-1}$ and  $(N_1,N_2,\cdots,N_s)$ its $\alpha$-numeration of $N$, with $N_s>0$.\\
    \esp Let $\beta \in \{ \{ n\alpha\},n\in\N\}$, and $(b_1,b_2,\cdots, b_{s'})$ its $\alpha$-numeration, with $b_{s'}>0$. We suppose $\beta\not = 0$.\\
     The minimal points in $\R^2$ of the set $\{ (\{ n\alpha \},n) \text{ such that }, 1\leqslant n \leqslant N\ \text{ and } \{ n\alpha \}\leqslant \beta \}$ are : 
     
 $\blacktriangleright$ Case 1 : if $b_1\geqslant 2$ or $b_2=0$, that is to say $\alpha\leqslant \beta$, these are the same as in former Lemma ! 
 

 $\blacktriangleright$ Case 2 : else ( $b_1=1$ and $b_2>0$). \\ 
 We denote $t=\min\{ i>0, b_{2i+1}>0 \}$ and $t=\lfloor s'/2 \rfloor$ if $\forall i>0, b_{2i+1}=0$ ( $s' $ is even in that case).\\ 
 (i) with the $\alpha$-numeration of $n$ :
 $$ (1,a_2,0,a_4,\cdots,a_{2t-2},0,j) \esp, \esp j\in \inter{b_{2t}}{a'_{2t}} \esp ; \esp  \text{ and } $$
 $$(1,a_2,0,a_4,\cdots,a_{2k-2},0,j)\esp ; \esp j \in \inter{1}{a'_{2k}} \esp ; \esp k\in\inter{t+1}{\lfloor s/2 \rfloor}$$ 
 
\esp  where $a'_{2k}=a_{2k}$ if $2k<s$ \esp ; \esp   $a'_s=N_s$.\\
 
  (ii)  with denominators of convergents of $\alpha$  :
  $$   q_{2k-2}+jq_{2k-1}, \esp  j\in\inter{c_{2k}}{a'_{2k}}, \esp  k\in \inter{t}{\lfloor s/2 \rfloor} $$
   with the same notations as in (i) and : $c_{2k}= \begin{cases} b_{2k} \text{ if } k=t \\
              1 \text{ else }\end{cases}$.
      
 \end{monlem}
 
 \textbf{Remark : } if $t>\lfloor s/2 \rfloor$, then there is no minimal element, so the set is empty : there is no integer $n$ such that $n\in \inter{1}{N}$ and $\{ n\alpha \} \leqslant \beta$. \\
 
 \begin{mademo}
 Let denote $E=\{ (\{ n\alpha \},n), n\in\inter{1}{N}\}$ and $E'=\{ (\{ n\alpha \},n)\in E, \{ n\alpha \}\leqslant \beta \}$. We claim that :
 $$ \min(E')=\min(E)\cap E'$$
 \esp Indeed, The second subset is obviously included in the first one and no element of $E\backslash E'$ is lower than any element of $E'$.\\
 \esp So, we just have to select among minimal elements $ (\{ n\alpha \},n)$ of $E$, whose $\alpha-$numeration, say $d$, is given by former Lemma, those that verify : $\{ n\alpha \}\leqslant \beta$, that is to say $d\leqslant_A b$.\\
 $\blacktriangleright$ Case 1 : if $b_1>1$ or $b_2=0$, that is to say  $\alpha\leqslant \beta$, then every $d$ of former Lemma verify the condition $d\leqslant_A b$.\\
  $\blacktriangleright$ Case 2 : else, we have $b_1=1$ and $b_2>0$. \\
  $\blacktriangleright \blacktriangleright$ subcase 1 :   suppose that $b_{2j+1}>0$ for some $j>0$. With notation of our Lemma, we have :
 $$ \forall i\in\inter{1}{t-1}, b_{2i}=a_{2i} \text{ and } b_{2i+1}=0$$
 \esp Let $d$ such that $ d_1=1$ and :
 $$  \forall i\in\inter{1}{k-1}, d_{2i}=a_{2i}, d_{2i+1} =0 \esp , \esp d_{2k}\in \inter{1}{a'_{2k}}, \forall i> 2k, d_i=0$$
 \esp  where $k\in\inter{1}{\lfloor s/2 \rfloor}$, with  $a'_{2k}=a_{2k}$ if $2k<s$, and $a'_s=N_s$ ( see former Lemma). \\
 \esp If $k<t$, then $d>_A b$ : indeed, $d_i=b_i$ for all $i<2k+2$, for $d_{2k}$ must be $a_{2k}=b_{2k}$, while $d_{2k+2}=0< b_{2k+2}$, because $b_{2k+1}=0$ and $b_{2t+1}>0$. \\
 \esp If $k=t$, then $d_i=b_i$ for all $i<2k$, so $d\leqslant_A b$ if and only if $d_{2k}\geqslant b_{2k}$, for $d_{2k+1}=0<b_{2k+1}$. \\
 \esp If $k>t$, then $d\leqslant_A b$, for $d_i=b_i$ for all $i<2t$ and $d_{2t}=a_{2t}\geqslant b_{2t} $ and $ d_{2t+1}=0< b_{2t+1}$. \\ 
   $\blacktriangleright \blacktriangleright$ subcase 2 : suppose that $b_{2i+1}=0$ for all $i>0$. Then, $s'$ can not be odd, so   :
  $$b=(1,a_2,0,a_4,\cdots,0,a_{s'-2},0,b_{s'})$$
  \esp The same arguments as above prove that a minimal element of former Lemma is in $E'$ if and only if $d_i=b_i$ for all $i<s'$ and $d_{s'}\geqslant b_{s'}$. The formula still rules.
 \end{mademo}

 \vspace*{0.6cm}
 
  $\bullet$ Now, we will need another generalization :   we study the minimal points in $\R^2$, with the product order, of the sequence $(\{ n\alpha -\beta \},n)_{n_0\leqslant n < q }$, where $n_0=0$ or $1$, $q$ is a positive integer and $\alpha=\frac{p}{q},\beta= \frac{p'}{q}$, with $p,p'\in \inter{1}{q-1}$ and $p,q$ coprime. \\
  
  \esp First, we remark that : 
  $$ \forall x\in \R, \{ x-\beta \} = \begin{cases} \{ x \} -\beta \in [0,1-\beta [ \text{ if } \{ x \} \geqslant \beta \\
              \{ x \} + 1- \beta \in [ 1- \beta , 1 [ \text{ if } \{ x \} < \beta \end{cases} $$
 
  \esp If the sequence begins with the index $i$ at $0$, the point $(1-\beta,0)$ obtained for $n=0$ is obviously a minimal point of the sequence. So, the other minimal points $(\{ n\alpha \},n)$ must verify $\{ n\alpha \}\geqslant \beta$. We have the same situation if the sequence begins at $i=1$ and if $\alpha \geqslant \beta$. \\ 
  
  \esp If the sequence begins at $i=1$ and if $\alpha < \beta$, then the " first" minimal points will be such that $\{ n\alpha\}<\beta$ and the following such that $\{ n\alpha \}\geqslant \beta$. So, the lowest integer $n$ such that $\{ n\alpha \}\geqslant \beta$ is important.\\
  
  \esp The case $\alpha=\beta$ is obvious ( only one or two minimal points for $n=0$ and $n=1$). \\

  \newpage

  \begin{monlem} \label{lem:20} let $q$ be  a positive integer and $\alpha=\frac{p}{q},\beta= \frac{p'}{q}$, with $p,p'\in \inter{1}{q-1}$ and $p,q$ coprime. We denote $\alpha=[0,a_1,\cdots,a_r,1]$ the CFE and $(b_1,b_2,\cdots,b_s)$ the $\alpha$-numeration of $\beta$, where $b_s>0$. \\
\esp We denote $t= \begin{cases} (s+1)/2 \text{  if }  \forall i>0,b_{2i}=0 \\ \min\{ i\in \N^*, b_{2i}\not = 0 \} \text{ else } \end{cases} $  ( $s$ is odd in the first case).\\

 \esp Minimal points of $(\{ n\alpha -\beta \},n)_{1\leqslant n < q }$ are obtained for following $n$, given by their $\alpha$-numeration : $(1)$ ( first ), $(b_i)_i$ ( last) and :
  $$ (1,a_2,0,\cdots,a_{2k-2},0,j) \esp ; \esp j\in\inter{1}{a_{2k}} \esp ; \esp  k\in\inter{1}{t-1} $$
  $$(b_1,b_2, \cdots,b_{2k-1}, j)\esp ; \esp j\in\inter{0}{b_{2k}-1} \esp ; \esp k\in\inter{t}{\lfloor s/2 \rfloor}$$

  \end{monlem}

  \vspace*{0.5cm}
  
  \textbf{Remark 1 : } if $b_2>0$, that is to say if $t=1$, then the first shape of  $ \alpha$-numeration of  $n$ is absent. In addition, if $b_1=1$ and $b_2>0$, then we count two times the first minimal point for $n=1$.\\
 
\textbf{Remark 2 : } if $b_{2i}=0$ for all $i$, then the second shape of  $ \alpha$-numeration of $n$ is absent, except for $(b_i)_i$.\\
 
 \textbf{Remark 3 : } if the sequence begins at $i=0$ instead of $i=1$, then the point $(1-\beta,0)$ obtained for $n=0$ is obviously a minimal point of the sequence. So, the other minimal points $(\{ n\alpha \},n)$ must verify $\{ n\alpha \}\geqslant \beta$. They are obtained for $n$ with $\alpha$-numeration :
    $$(b_1,b_2, \cdots,b_{2k-1}, j)\esp ; \esp j\in\inter{0}{b_{2k}-1} \esp ; \esp k\in\inter{t}{\lfloor s/2 \rfloor}$$
 
  \vspace*{0.5cm}
 
  \begin{mademo} 
  \esp As long as $\{ n\alpha \}<\beta$, minimal points of our sequence are obtained for the same values of $n$ than those of $(\{n\alpha \},n)_n$, because $\{ n\alpha-\beta \}=\{n\alpha \}+1-\beta$.\\
  \esp We remark that ( see definition of $E_{\alpha}$) :
  $$ \forall i\in\inter{1}{t-1}, b_{2i}=0 \esp ; \esp b_{2i-1}=a_{2i-1} $$
  \esp  Now, we need the $\alpha$-numeration, say $\nu$, of $n_1$, the least integer $n$ such that $\{ n\alpha \} \geqslant \beta$. If $ \alpha \geqslant\beta$ then $n_1=1$. Else, with Proposition \ref{prop:9}, $\nu$ is the minimum of elements $d=(d_1,\cdots,d_r)$ of $E_{\alpha}$ for RLO, such that $d \geqslant_A b$, where $b=(b_k)_{k\in \inter{1}{r}}$ and $\geqslant_A$ means for ALO.\\
 \esp  We claim that $\nu=(b_1,b_2,\cdots,b_{2t-1})$. Indeed, the condition $d \geqslant_A b$ implies that :
 $$ \forall i\in\inter{1}{t-1}, d_{2i}=0=b_{2i} \esp ; \esp d_{2i-1}=a_{2i-1}=b_{2i-1} \text{  and } d_{2t-1}\geqslant b_{2t-1}$$
 \esp But, the above $\delta$ satisfies $\nu\geqslant_A b$ and $\nu$ is minimal ( for RLO) among these one. Furthermore :
  $$ \forall n<n_1, \{ n\alpha - \beta \}\in [ 1- \beta,1[ \esp ; \esp \{ n_1\alpha - \beta \}\in [ 0,1- \beta[ $$
  \esp Note that this is also true if $b_{2i}=0$ for all $i>0$ and that, in that case, we have $\nu=b$ and so $\{ \nu\alpha \}=\beta$ and no minimal point after this one... \\
  \esp We return to the general case : we have seen just above that the lowest integer $n$ such that $\{ n\alpha \}\geqslant \beta$ is $n_1=(b_1,b_2,\cdots,b_{2t-1})_{\alpha}$. So, for $n<n_1$, we use values of $n$ of Lemma \ref{lem:18} : we obtain first part of our result ( that involves the $(a_i)$), for $b_{2t-1}>0$ : indeed $b_{2t}>0$ and $b_{2t-2}=0$.\\
  \esp Now, if $n\geqslant n_1$, let denote $d=(d_1,\cdots,d_r)$ its $\alpha$-numeration. Then, the minimality condition for $(\{ n\alpha - \beta \},n)$ is equivalent to : $d\geqslant_A b$ and $ d$ is minimal among these ( elements of $E_{\alpha}$ greater than $b$ for ALO) for the product of orders (ALO,RLO). So, $d$ must verify $(b_1,\cdots,b_s)\leqslant_A d <_A (b_1,\cdots,b_{2t-1})$. It gives the successive values of $d$ : $(b_1,\cdots,b_{2t-1},j)$ as $j\in \inter{0}{b_{2t}-1}$. Later, if $b_{2t+2}\not = 0$, the first ( for RLO) $d$ that verify $d\geqslant_A b$ and that is lower ( for ALO) than $(b_1,b_2,\cdots,b_{2t-1},b_{2t}-1)$ is $(b_1,b_2,\cdots,b_{2t-1},b_{2t},b_{2t+1})$. Then, we have $(b_1,b_2,\cdots,b_{2t-1},b_{2t},b_{2t+1},j)$ with $j$ taking successive values of $\inter{1}{b_{2t+2}-1}$...and so on.
  \end{mademo}

\newpage

   \subsection{number of integers $k$ such that $\{ k\alpha \}\leqslant \beta $ and $k\leqslant \nu$}

   \label{subsec:counting}
     
   \esp We consider two rationals $\alpha=\frac{p}{q}$ and $\beta=\frac{m}{q}$ in $[0,1[$ and  a positive integer $\nu$, such that :
   $$ \gcd(p,q)=1 \esp ; \esp 0<m<q \esp ; \esp 0< \nu < q $$
   \esp We denote $n=(n_k)_k$ and $b=(b_k)_k$ the respective $\alpha$-numeration of $\nu$ and $\beta$.  We denote $\sigma$ the usual shift on sequences.\\
   
   \esp We will also use the two total orders on finite sequences of reals : RLO, denoted $\leqslant_R$ and ALO, denoted $\leqslant_A$ ( see \ref{subsec:alphanum}). \\

 \esp We also denote :
   $$ C(\alpha,\beta,\nu)= \#\{ k\in \inter{0}{\nu}, \{ k\alpha \}\leqslant \beta \} $$ 
   
     \vspace*{0.5cm}

   \esp We will denote $[a_0,a_1,\cdots,a_r,1]$ the CFE of $\alpha$ ( with $a_0=0$).\\
   
   \esp We also use the following notations : 
    $$ \alpha_0=\alpha \esp ; \esp \forall k\in \inter{1}{r},\alpha_k=\left\{ \frac{1}{\alpha_{k-1}}\right\}$$
 $$ \nu_0=\nu \esp ; \esp \forall k\in \inter{1}{r-2}, \nu_k= \begin{cases} \lfloor \nu_{k-1}\alpha_{k-1}\rfloor  \text{ if } n_{k}\not = 0 \text{ or } n_{k+1}=0 \\ 
   \lfloor \nu_{k-1}\alpha_{k-1}\rfloor +1 \text{ else } \end{cases}  $$
 $$ \beta_0=\beta \esp ; \esp \forall k\in \inter{1}{r}, \beta_k= \frac{1}{\alpha_{k-1}}(b_k\alpha_{k-1}-\beta_{k-1})$$

  \vspace*{1cm} 
   
   \begin{maprop} \label{prop:11} we denote $n=(n_k)_k$ the $\alpha$-numeration of $\nu$ and $b=(b_k)_k$ the $\alpha$-numeration of $\beta$. We denote $s$ the minimum of the lengths of $n$ and $b$, when we drop the eventual infinite " $0$-tail". So, $n_s$ or $b_s$ is not null, but $\sigma^s(n)$ or $\sigma^s(b)$ is the null sequence.\\
   $$ C(\alpha,\beta,\nu)=D+ \sum_{i=1}^{s } (-1)^{i-1}[  b_i \nu_i +\tau_i+\epsilon_i-\epsilon'_i]$$
   \esp where :
 $$D = \textbf{1}_{ n\leqslant_A b}+ \textbf{1}_{ b\leqslant_R n}-\textbf{1}_{ n=b} \esp ; \esp
\tau_i= \begin{cases} 1 \text{ if } n_in_{i+1}=0 \text{ and }  \sigma^i(n)\not = (0) \\ \min(b_i,n_i) \text{ else} \end{cases} $$
 $$\epsilon'_i= \begin{cases} 1 \text{ if }   \sigma^i(b) <_R \sigma^i(n) \\ 0 \text{ else } \end{cases}\esp ;  \esp \epsilon_i= \begin{cases} 1 \text{ if }  b_i< n_i \text{ and } \sigma^i(b) <_A \sigma^i(n) \\ 0 \text{ else } \end{cases}$$

   \end{maprop}
  \begin{mademo} see \textbf{[3]} 4.4. ( \nameref{biblio}) \end{mademo}

 \newpage

\section{Invariants of $\frac{<a,b>}{d}$}  

 \label{sec:invariants}

      \esp Let $a$ and $b$ be two coprime positive integers. We want to study the invariants of $ \frac{<a,b>}{d}$.\\
      \esp  We remind that we can suppose,  without loss of generality, that $a,b,d$ are pairwise coprime. In addition, we can also suppose that $d\not\in <a,b>$, because $ \frac{<a,b>}{d}=\N$ if $d\in <a,b>$.

    \subsection{multiplicity}

      \esp The multiplicity of $ \frac{<a,b>}{d}$ is the smallest positive integer $k$ such that $kd\in <a,b>$. We could deduce a formula from the result concerning the irreducible elements, because the multiplicity is the smallest element among these one. But, the following theorem is much more interesting, for it gives a very simple expression of $m(  \frac{<a,b>}{d})$ in terms of " best rational in an interval"...  \\
      
     \begin{monlem} \label{lem:21} let $n$ be a positive integer and $a,b$ be two coprime positive integers. We suppose that $n=xa-yb=ab(\frac{x}{b}-\frac{y}{a})$, with $x,y$ integers. Then :
     $$ n\in <a,b>\backslash  ( a \inter{0}{b-1}) \Leftrightarrow \left\lfloor \frac{x}{b} \right\rfloor \not =  \left\lfloor \frac{y}{a} \right\rfloor  $$
     \end{monlem}
     \begin{mademo}
     - we suppose that $ n\in <a,b>\backslash (a\inter{0}{b-1})$, then we can find two integers $i,j$, such that $n=ai+bj$, $j$ is positive and $i$ is non negative. Indeed, if $j=0$, then $n=ai$ and $i\geqslant b$, so $n=ai'+bj'$, with $i'=i-b$ and $j'=a$. So : $(x-i)a=(y+j)b$. But, $a$ and $b$ are coprime, so it exists an integer $k$ such that :
     $$ x-i=kb \text{ and } y+j=ka $$
     So : $\frac{y}{a}<k\leqslant \frac{x}{b}$ and we deduce $\left\lfloor \frac{x}{b} \right\rfloor \not =  \left\lfloor \frac{y}{a} \right\rfloor$. \\
     - we suppose that $\left\lfloor \frac{x}{b} \right\rfloor \not =  \left\lfloor \frac{y}{a} \right\rfloor$. We have $\frac{n}{ab}= \frac{x}{b}- \frac{y}{a} $ and $n$ is positive, so, there exists an integer $k$ such that :  $\frac{y}{a}<k\leqslant \frac{x}{b}$. Now, we denote $i=x-kb$ and $j=ka-y$. Then, $i\in\N,j\in \N^*$ and $n=ax-by=ai+bj$. So, $ n\in <a,b>$. \\
     \esp Now, suppose that  $n\in  a\inter{0}{b-1}$. Then we can use our first argument with $j=0$ and $i\in \inter{0}{b-1}$. We obtain :  it exists an integer $k$ such that :
     $$ x-i=kb \text{ and } y=ka $$
     So, $\frac{y}{a}=k=\left\lfloor \frac{y}{a} \right\rfloor$. But $\frac{x}{b}=k+\frac{i}{b}\in [k,k+1[$, so $ \left\lfloor \frac{x}{b}\right\rfloor =  \left\lfloor \frac{y}{a} \right\rfloor$ : contradiction.    \end{mademo}

     \begin{montheo}  \label{theo:4} let $a,b,d$ be three pairwise coprime positive integers such that $d\not\in <a,b>$ and $d\geqslant 2$. \\
      We can find   $\alpha' \in \frac{1}{b}\Z$ and $\alpha\in \frac{1}{a}\Z$ such that $\alpha'-\alpha=\frac{d}{ab}$. Then : \\
     $$ m\left(\frac{<a,b>}{d}\right) = \text{ the denominator of the reduced  best rational in } [\alpha,\alpha']$$ 
       \end{montheo}

     \begin{mademo} We can write $\alpha'=\frac{x}{b}$ and $\alpha=\frac{y}{a}$, with $x$ and $y$ are integers such that  $d=ax-by$ ( basic arithmetic, since $a$ and $b$ are coprime). We remark that : since $d\not\in <a,b>$, then, with our previous Lemma, $\lfloor \alpha \rfloor= \lfloor \alpha' \rfloor $. We denote $ p/q$ the reduced best rational in $[\alpha,\alpha']$. We have two cases :  if $p/q=\alpha$, then $q=a$ for $y$ and $a$ are coprime ( $ d$ and $a$ are coprime). So $qd\in <a,b>$. Else, we have $qd=ab(q\alpha -a\alpha')$ and $\lfloor q\alpha \rfloor \not = \lfloor q\alpha' \rfloor$ ( see Proposition \ref{prop:8}), since $p/q$ is the best rational in $]\alpha,\alpha']$ and $\alpha < \alpha'$. So our previous Lemma gives : $qd \in <a,b>$. In both cases, we have $qd\in <a,b>$.\\
     \esp Now, for $k\in \inter{1}{q-1}$, $kd=akx-bky$ and Proposition \ref{prop:8} gives : $\lfloor k\alpha \rfloor  = \lfloor k\alpha' \rfloor$ ( there is no rational $i/k$ in $]\alpha,\alpha']$). So : $kd\not\in <a,b>$ or $kd\in a\N$. If $kd\in a\N$, then $k\in a\N$, for $a$ and $d$ are coprime. But $q\leqslant a$, for $ \alpha=\frac{y}{a}$, so $k<a$, contradiction.          \end{mademo}

 \newpage
     
      \subsection{minimal generators}
      
    $\bullet$ We begin with a result that describes precisely the minimal generators of $ \frac{<a,b>}{d}$, namely $\Irr(\frac{<a,b>}{d})$ and so its cardinality $e(\frac{<a,b>}{d})$. The case $d>\max(a,b)$ is the most intricate...
      
        \begin{montheo} \label{theo:5}
             let $a,b,d$ be three pairwise coprime positive integers such that $d\not\in <a,b>$. \\
         \esp     Let $m\in \inter{1}{d-1}$ such that $am+b = 0 [d]$. Let $\alpha=m/d$ and $[0,a_1,\cdots,a_r,1]$ its CFE.\\
         \esp Let $(p_k/q_k)_k$ its reduced convergents and for  $k\in \inter{-1}{r+1},\esp \delta_k=(-1)^k(q_k\alpha-p_k)$, as defined in \ref{subsec:alphanum}.\\
         \esp We denote $\tau=\frac{am+b}{d}$ and         $\mu_i=\tau q_i-ap_i$ for all $i\in \inter{0}{r}$.\\
     \esp     We denote $(N_1,\cdots,N_s)$ the $\alpha$-numeration of $a-1$ ( with $N_s>0$). \\
     \esp We denote $x_0=\begin{cases} d(\delta_{s-2}-N_s\delta_{s-1}) \text{ if } s \text{ is even}\\ d\delta_{s-1} \text{ else} \end{cases} \esp ; \esp  \esp  \nu= \begin{cases} N_s \text{ if } s \text{ is even} \\ 0 \text{ else }\end{cases}$  .\\
   \esp All sets below are written with pairwise distincts elements. \\

    $\blacktriangleright$ Case 1 : $d<a<b$.
             $$ e\left( \frac{<a,b>}{d}\right)=3+ \sum_{k=1}^{\lfloor r/2 \rfloor }a_{2k} $$
        $$ \Irr \left( \frac{<a,b>}{d}\right)=\left\{ a \esp ; \esp b \esp ; \esp  \tau \esp  ;  \esp   \mu_{2k-2}+j\mu_{2k-1}, \esp j\in\inter{1}{a_{2k}}, k\in \inter{1}{\lfloor r/2 \rfloor} \right\} $$      
     
   $\blacktriangleright$ Case 2 : $a<d<b$ or ($a<b<d$ and $x_0+b>d$).
  
   $$ e\left( \frac{<a,b>}{d}\right)=2+ \sum_{k=1}^{\lfloor \frac{s-1}{2} \rfloor }a_{2k}+ \nu $$

           $$\Irr \left( \frac{<a,b>}{d}\right)=\left\{ a \esp ; \esp \tau \esp ; \esp   \mu_{2k-2}+j\mu_{2k-1}, \esp \begin{cases} j\in\inter{1}{a_{2k}} \text{ if } k\in \inter{1}{\lfloor \frac{s-1}{2} \rfloor}  \\
   j\in\inter{1}{N_s} \text{ if } k=s/2 \end{cases}\right\} $$

     $\blacktriangleright$ Case 3 : $a<b<d$ and $m< x_0+b \leqslant d $. \\
     \esp Same results as in Case 2, except that $a\not\in \Irr \left( \frac{<a,b>}{d}\right)$ and so $ e\left( \frac{<a,b>}{d}\right)=1+ \sum\limits_{k=1}^{\lfloor \frac{s-1}{2} \rfloor }a_{2k}+ \nu$. \\
     
    $\blacktriangleright$ Case 4 : $a<b<d $ and $x_0+b\leqslant m$. \\
\esp We denote $\beta=\frac{x_0+b-1}{d}$ and $(b_1,b_2,\cdots, b_{s'})$ its $\alpha$-numeration, where $b_{s'}>0$.\\
\esp We also denote $t=\min\{ i>0, b_{2i+1}>0 \}$ and $t=s'/2$ if $\forall i>0, b_{2i+1}=0$ ( $s'$ is even in that case). \\
\esp Then , $a\not\in \Irr\left( \frac{<a,b>}{d}\right)$ and :
$$\Irr\left( \frac{<a,b>}{d}\right)=\left\{   \mu_{2k-2}+j\mu_{2k-1} \esp , \esp j\in\inter{c_{2k}}{a'_{2k}} \esp , \esp k\in \inter{t}{\lfloor s/2 \rfloor} \right\}$$
  \esp  where $c_{2k}=\begin{cases} b_{2k} \text{ if } k= t \\ 1 \text{ else } \end{cases}$ and $\esp a'_{2k}= \begin{cases} a_{2k} \text{ if } \esp 2k<s \\ N_s \text{ if } \esp  2k=s \end{cases}$.\\
     $$  e\left( \frac{<a,b>}{d}\right)=1+a_{2t}-b_{2t}+\sum_{k=t+1}^{\lfloor \frac{s-1}{2} \rfloor }a_{2k}+\nu$$

        \end{montheo}

   \textbf{Remark : }  the hypothesis of Case 4 is equivalent to $\beta<\alpha$, that is to say $b_1=1$ and $b_2>0$.    

\newpage


   
    \begin{mademo} as in section \ref{sec:representation}, we denote $S'=  \frac{<a,b>}{d}$ and $\psi : \begin{cases} \Z^2\to \frac{1}{d}\Z \\ (x,y)\to \frac{ax+by}{d} \end{cases}$ and $L=\psi^{-1}(\Z)$. \\
    \esp We can describe $(\inter{1}{d-1})^2\cap L$ with the finite sequence $(d\{n\alpha \},n)_{n\in \inter{1}{d-1}}$. Indeed : \\
    \esp For all $\forall n\in \inter{1}{d-1}, d\{n\alpha \}\in\inter{1}{d-1}$, for $\{ n\alpha \}\not = 0$ if $n\not = 0 [d]$, because $m$ and $d$ are coprime ( $b$ and $d$ are coprime). Moreover :
    $$\forall n\in \inter{1}{d-1},ad\{n\alpha \}+bn=(am+b)n-ad\lfloor n\alpha \rfloor =0[d]$$
    \esp In addition, each $d\times d$ square in $L$ contains a unique point on each line and on each row... \\
    \esp Now, let us prove that, if $a<d$, then : $x_0=\min\{ x(u), u\in L\cap ( \N\times \inter{1}{a-1})\}$, in order to use Proposition \ref{prop:2}. First, we can replace, in this formula, $\N$ by $\inter{0}{d-1}$, for $L$ is $(d,0)$-invariant. We have seen above that $L\cap ( \inter{0}{d-1}\times \inter{1}{a-1})$ is  enumerated by  $(d\{n\alpha \},n)_{n\in \inter{1}{a-1}}$, because $a<d$. So, if we denote $x'_0=\min(L\cap ( \inter{0}{d-1}\times \inter{1}{a-1}))$, then $x'_0=d\min ( \{n\alpha\}, n\in\inter{1}{a-1})$. Using our $\alpha$-numeration and Lemma \ref{lem:18}, we obtain that this minimum is obtained for the following value of $n$ :
    $$n=\begin{cases} q_{s-3}+a_{s-1}q_{s-2}=q_{s-1} \text{ if } s \text{ is odd } \\ q_{s-2}+N_sq_{s-1}  \text{ if } s \text{ is even } \end{cases}$$
    \esp This implies that $x_0=x'_0$ and our formula...\\
   
      $\blacktriangleright$  Case 1 : with Proposition \ref{prop:2} : 
    $$ \Irr(S')=\{ a,b \}\cup \psi(\min((\inter{1}{d-1})^2\cap L))$$
   \esp But, with remarks at the beginning of this proof and Lemma \ref{lem:17}, we have $ \min(  (\inter{1}{d-1})^2\cap L))$ enumerated by $(d\{ \alpha \},1)$ and the $(d\{n\alpha \},n)$  with :
    $$  n=(1,a_2,0,\cdots,a_{2k-2},0,j)_{\alpha}, \esp k\in\inter{1}{\lfloor r/2 \rfloor}, j\in\inter{1}{a_{2k}}$$
    \esp We obtain ( see Proposition \ref{prop:10} and formulae around) :
    $$ n=q_{2k-2}+jq_{2k-1} \esp ; \esp \{ n\alpha \}= \delta_{2k-2}-j\delta_{2k-1}$$
    \esp But $\psi$ is into on $(\inter{1}{d-1})^2$ since $d<a<b$, so we deduce the formula for $e\left(\frac{<a,b>}{d}\right)$. Moreover, we have $\psi(d\{ \alpha \},1)=a\{ \alpha \}+\frac{b}{d}=\tau$ and $ \forall k\in\inter{1}{r}, j\in\inter{1}{a_{2k}}$ : 
    $$  n=q_{2k-2}+jq_{2k-1} \Rightarrow
    \psi((d\{ n\alpha \},n))=a(\delta_{2k-2}-j\delta_{2k-1})+\frac{b}{d}(q_{2k-2}+jq_{2k-1})$$
    
    $\blacktriangleright$  Case 2 : if $a<d<b$, the arguments are essentially the same as in Case 1. The only difference is a limitation on the value of $n$ : with Proposition \ref{prop:2}, we have $n\in\inter{1}{a-1}$ instead of $n\in\inter{1}{d-1}$. With Lemma \ref{lem:18}, we obtain the result. \\
    \esp If $a<b<d$ and $x_0+b>d$ : with Proposition \ref{prop:2} and its notations,  we have $x_1=d$ and $a\in \Irr(S')$. In addition :
    $$   \Irr(S')=\{ a \}\cup \psi(\min((\inter{1}{d-1}\times \inter{1}{a-1})\cap L))$$
    \esp for $a$ is obtained with the point $(d,0)$ in $L$. So, the results are the same as in the case $a<d<b$. \\
    
    $\blacktriangleright$ Case 3 : if $a<b<d$ and $m< x_0+b \leqslant d $. With Proposition \ref{prop:2} again, we have $x_1=x_0+b-1$ and $a\not\in  \Irr(S')$. Likewise :
    $$  \Irr(S')= \psi(\min((\inter{1}{x_0+b-1}\times \inter{1}{a-1})\cap L))$$   
    \esp But $(\inter{1}{x_0+b-1}\times \inter{1}{a-1})\cap L$ is enumerated by the $(d\{n\alpha \},n)$, such that $\{n\alpha \}\leqslant \beta$ and $n\leqslant a-1$. Our hypothesis $ m< x_0+b$ means that $\alpha\leqslant \beta$, so the point $(d\{ \alpha \},1)$ is in $(\inter{1}{x_0+b-1}\times \inter{1}{a-1})\cap L$. We deduce that the condition "$\{n\alpha \}\leqslant \beta$" is dispensable for the minimal points of $(\inter{1}{x_0+b-1}\times \inter{1}{a-1})\cap L$. Thus, we are in the same condition as in previous Case, except for $a$...\\
    
      $\blacktriangleright$ Case 4 : if $a<b<d $ and $x_0+b\leqslant m$. Let us remark that $\beta<\alpha$, hence $b_1=1$ and $b_2>0$. We can take up the arguments of Case 3, except that the condition "$\{n\alpha \}\leqslant \beta$" is important. So, we can conclude with Lemma \ref{lem:19}.    \end{mademo}

\newpage
  
 $\bullet$ What are the extremal possible values of  $ e\left( \frac{<a,b>}{d}\right)$ and when does it happen ? \\
  
   \begin{maprop} [ Extremal cases]  \label{prop:12}
   
   \esp We keep the notations of previous theorem.\\
  
    $\blacktriangleright$ Case 1 : if $d<a<b$, then : 
    $$3\leqslant  e\left( \frac{<a,b>}{d}\right)\leqslant d+1$$ 
$ \triangleright  \esp   e\left( \frac{<a,b>}{d}\right)=3\Leftrightarrow d / (a+b) $. In that case, $\Irr\left(  \frac{<a,b>}{d}\right)= \{ a,b, \frac{a+b}{d} \} $.\\
     $\triangleright  \esp  e\left( \frac{<a,b>}{d}\right) = d+1 \Leftrightarrow d /( b-a )$. In that case, $\Irr\left( \frac{<a,b>}{d}\right)=\left\{ a+ \frac{n(b-a)}{d}, n\in\inter{0}{d} \right\} $. \\
 
    $\blacktriangleright$ Case 2 : if $a<d<b$ or ($a<b<d$ and $x_0+b>d$) , then :
    $$2\leqslant  e\left( \frac{<a,b>}{d}\right)\leqslant a$$
     $  \triangleright  \esp  e\left( \frac{<a,b>}{d}\right)=2\Leftrightarrow \frac{m}{d}<\frac{1}{a-1}$. In that case, $\Irr\left( \frac{<a,b>}{d}\right)=\{ a, \tau \}$.\\
     $\triangleright  \esp   e\left( \frac{<a,b>}{d}\right)=a \Leftrightarrow \frac{m}{d}>1-\frac{1}{a-1}$. In that case, $\Irr\left( \frac{<a,b>}{d}\right)=\left\{ a+n(\tau-a), n\in\inter{0}{a-1}\right\} $. \\
    
     $\blacktriangleright$ Case 3 : if  $a<b<d$ and $m< x_0+b \leqslant d $, then : 
     $$    e\left( \frac{<a,b>}{d}\right)\leqslant a-1$$
       $\triangleright  \esp   e\left( \frac{<a,b>}{d}\right)=a-1 \Leftrightarrow \frac{m}{d}>1-\frac{1}{a-1}$. In that case, $\Irr\left( \frac{<a,b>}{d}\right)=\left\{ a+n(\tau-a), n\in\inter{1}{a-1}\right\} $. \\
     
      $\blacktriangleright$ Case 4 : if $a<b<d $ and $x_0+b\leqslant m$, then :
      $$    e\left( \frac{<a,b>}{d}\right)\leqslant a-2$$
        $ \triangleright  \esp  e\left( \frac{<a,b>}{d}\right)=a-2 \Leftrightarrow \begin{cases} \{ 2\alpha \}\leqslant \beta \\
              \alpha \geqslant 1- \frac{\{ 2 \alpha \}}{a-3} \end{cases}$. In that case, $\Irr\left( \frac{<a,b>}{d}\right)=\left\{ a+n(\tau-a), n\in\inter{2}{a-1}\right\} $. \\
   \end{maprop}
   
 \vspace*{0.5cm}  
   
   \textbf{Remark 1 : } for all of these cases, $\Irr\left( \frac{<a,b>}{d}\right)$ has elements in arithmetic progression, with expression $a+n(\tau-a)$, but the range of $n$ is different for each fo these cases.\\
   
   \textbf{Remark 2 : } if $a<d<b$ and $d$ divides $a+b$, then $m=1$, so $\frac{m}{d}<\frac{1}{a-1}$ and $\Irr\left( \frac{<a,b>}{d}\right)=\{ a, \tau \}$, this is consistent with respect to  Case 1. \\
   \esp Similarly, if $a<d<b$ and $d$ divides $b-a$, then $m=d-1$, so  $ \frac{m}{d}>1-\frac{1}{a-1}$ and $e$ is maximal, as in Case 1. \\

  \textbf{Examples : } \\
   We take $a=151$ and $b=503$. Then :
   $$ \Irr\left( \frac{<a,b>}{6}\right)=\{ 109,151,503 \} \esp ; \esp  \Irr\left( \frac{<a,b>}{32}\right)= \{ 151, 162, 173, 184, 195, 206, 217, 228,... 503 \}$$
    $$ \Irr\left( \frac{<a,b>}{218}\right)=\{ 3, 151 \} \esp ; \esp  \Irr\left( \frac{<a,b>}{176}\right)= \{ 151, 153, 155, 157, 159, 161, 163, 165, ... 451 \}$$

\vspace*{1cm}

   \begin{mademo}
     $\blacktriangleright$ Case 1 :  the formula of Theorem \ref{theo:5} (i) proves that $3\leqslant  e\left( \frac{<a,b>}{d}\right)$. Proposition \ref{prop:2} gives the other inequality, as there are at most $d-1$ minimal elements in $\inter{0}{d-1}^2$.\\
   --- Now, $  e\left( \frac{<a,b>}{d}\right) =3$ if and only if $a_{2k}=0$ for all $k>0$, with notations of Theorem \ref{theo:5}. But, the $a_i$ are positive integers for positive indices $i$, so it means that $r=0$ or $1$. This is equivalent to the fact that $\alpha=[a_0,a_1,1]$ or $\alpha=[a_0,1]$. But, $\alpha\not\in \Z$, since $d$ does not divide $b$, so : 
   $$  e\left( \frac{<a,b>}{d}\right) =3 \Leftrightarrow \frac{1}{\{\alpha \}}\in \N^* \Leftrightarrow \frac{d}{m}\in\N^* $$
   \esp But, $d$ and $m$ are coprime, so $ e\left( \frac{<a,b>}{d}\right) =3$ is equivalent to $ m=1 [d]$, that is to say $ a+b=0[d] $. In that case, we have 3 minimal generators of $ \frac{<a,b>}{d}$ : $a,b$ and $\tau$. \\
   --- Now, $ e\left( \frac{<a,b>}{d}\right)$ is maximal if and only if all points $(d\{ n\alpha \},n)_{n\in \inter{1}{d-1}}$ are minimal in $L\cap \inter{1}{d-1}^2$ ( see proof of Theorem \ref{theo:5}). This means that the sequence $(\{ n \alpha \})_{n\in\inter{1}{d-1}}$ is decreasing ( with notations of Theorem \ref{theo:2}). This is equivalent to $d\{ \alpha \}=d-1$, since all these points have different $x$ and different $y$ in $\inter{1}{d-1}$. In that case, $d \{ n\alpha\}=d-n$ for all $n\in \inter{1}{d-1}$. Hence, all the points above are pairwise not comparable. We can conclude : 
   $$  e\left( \frac{<a,b>}{d}\right)=d+1 \Leftrightarrow m=d-1 [d] \Leftrightarrow d / (b-a) $$
   \esp In that case, the minimal generators are : $a,b$ and the $\frac{1}{d}\left(a(d-n)+bn\right)$ where $n\in\inter{1}{d-1}$.  \\
   
   $\blacktriangleright$ Case 2 : if $a<d<b$  or ($a<b<d$ and $x_0+b>d$), again, with Theorem \ref{theo:5} and Proposition \ref{prop:2}, we obtain both inequalities. \\
   --- we have with same arguments as in Case 1 :
   $$e(S')=2 \Leftrightarrow s=1 \Leftrightarrow a-1\leqslant a_1 \Leftrightarrow a-1 \leqslant \lfloor 1/\alpha \rfloor \Leftrightarrow a-1 \leqslant d/m$$
   \esp   In that case, according to Theorem \ref{theo:5}, $\Irr(S')=\{ a , \tau \}$.\\
   --- we have $e(S')=a$ if and only if all points $(d\{ n\alpha \},n)_{n\in \inter{1}{a-1}}$ are minimal in $L\cap \inter{1}{d-1}\times \inter{1}{a-1}$ ( see proof of Theorem \ref{theo:5}).  This means that the sequence $(\{ n \alpha \})_{n\in\inter{1}{a-1}}$ is decreasing, that is to say $u_{n+1}-u_n=\alpha-1$ if we denote $u_0=1$ and $u_n=\{ n\alpha \}$. This is equivalent to , $(a-1)(1-\alpha)< 1$ and  gives our condition. In that case, we obtain : 
   $$ \forall n\in\inter{1}{a-1}, \{n\alpha \}=  \alpha - (n-1)(1-\alpha)= n\alpha - (n-1)$$ 
   \esp We deduce easily the expression of irreducible elements of $S'$. \\
   
   $\blacktriangleright$ Case 3 :   if  $a<b<d$ and $m< x_0+b \leqslant d $... see Theorem \ref{theo:5} and Case 2. \\
   
    $\blacktriangleright$ Case 4 : if $a<b<d $ and $x_0+b\leqslant m$. We use Proposition \ref{prop:2}, with the fact that there is no point in $L\cap \inter{1}{x_1}\times \inter{0}{a-1}$, on the lines $y=0$ and $y=1$. So it remains $a-2$ lines for possible minimal points... \\
    --- we have $e(S')=a-2$ if and only if each of the $a-2$ remaining lines contains one minimal point. This is equivalent to the fact that the sequence $(\{ n\alpha \})_{2\leqslant n\leqslant a-1}$ is decreasing and $\leqslant \beta$. This gives our conditions, after similar considerations as in Case 2. The result on $\Irr(S')$ in that case is obtained as in Case 2.   
   \end{mademo}

 \vspace*{0.5cm}

     \newpage

      \subsection{pseudo-Frobenius numbers}
      
     \esp We will distinguish several cases to make the statements more readable and in our proofs, we will use different " half-turn" in $\Z^2$ to compute maximal points of a lattice in a rectangle ( see \ref{subsec:minimal}) : \\
     
     $\bullet$ Let us begin with special cases, when $\PF$ is reduced to a single element or is obtained from $\Irr$ via a  symmetry :

      \begin{montheo} [ Easy cases] \label{theo:6}.\\ 
        let $a,b,d$ be 3 pairwise coprime positive integers such that $d\not\in <a,b>$ and $a<b$. \\
      $\blacktriangleright$ Case 1 : if $(a-1)(b-1)=1[d]$, then  $\frac{<a,b>}{d}$  is symmetric ( its type is one) and :
 $$  f\left(\frac{<a,b>}{d}\right)= \frac{f(<a,b>)}{d}=\frac{ab-a-b}{d}$$
   $\blacktriangleright$ Case 2 : if $a=1[d]$, then :
    $$ \PF\left(\frac{<a,b>}{d}\right)= \frac{(a-1)b}{d}-\left(\Irr\left(\frac{<a,b>}{d}\right)\backslash\{ b\} \right)$$
  
    $\blacktriangleright$ Case 3 : if $b=1[d]$, then :
    $$ \PF\left(\frac{<a,b>}{d}\right)= \frac{(b-1)a}{d}-\left(\Irr\left(\frac{<a,b>}{d}\right)\backslash\{ a\} \right)$$
    \esp For cases 2 and 3, we obtain : $t=e-1$ for $\frac{<a,b>}{d}$.
    \end{montheo}       
    \textbf{Remark : } Case 1 is not compatible with Case 2 or Case 3. But we can have $a=b=1[d]$ : in that case, we have $d<\min(a,b)$ and $b-a=0[d]$, so $\Irr(<a,b>/d)=\left\{ a+ \frac{n(b-a)}{d}, n\in\inter{0}{d} \right\} $. We obtain, if we denote $f=ab-a-b$ : 
    $$\PF\left(\frac{<a,b>}{d}\right)=\left\{ \frac{1}{d}[f-nb-(d-1-n)a], \esp n\in\inter{0}{d-1}\right\}$$    

    \begin{mademo}  we will denote $S'=\frac{<a,b>}{d}$ and, as in Lemma \ref{lem:7} :   $\psi$ $\begin{cases} \Z^2\to (1/d) \Z \\ (x,y)\to\frac{ax+by}{d}\end{cases}$ and : 
    $$ L=\{ (x,y)\in\Z^2, ax+by=0[d]\} \esp ; \esp R= \inter{b-d}{b-1}\times \inter{-d}{-1} $$ 
     
       $\blacktriangleright$ Case 1 : if $(a-1)(b-1)=1[d]$, then $(b-1,-1)$ is an element of $L$. So, it is the maximum of $L\cap R$. We conclude that $\frac{<a,b>}{d}$ has only one $\PF$-number, which is $\psi(b-1,-1)$... \\
        $\blacktriangleright$ Case 2 : if $a=1[d]$, then $(b,-1)$ is in $L$, but not in $R$. Let us consider $\sigma : (x,y)\to (b-x,-1-y)$, that maps $L\cap R$ onto $L\cap R'$, where $R'=\inter{1}{d}\times \inter{0}{d-1}$. So, $\max(L\cap R)= \sigma(\min(L\cap R'))$.\\
        \esp But $L\cap R'=(L\cap R")\backslash \{ (0,d) \}$, if we denote $R"=\inter{0}{d}^2\backslash \{ (0,0) \}$. This implies that $\min(L\cap R')=\min(L\cap R")\backslash \{ (0,d) \}$, for $(0,d)$ is a minimal point of $L\cap R"$. In addition, $d<a<b$, so $\Irr(S')=\psi(\min(L\cap R"))=\psi(\min(L\cap R"))\cup \{ b\}$ ( see Proposition \ref{prop:2}). Now, we obtain, with Lemma \ref{lem:7}  :
        $$ \PF(S')=\psi(\max(L\cap R))= \frac{ab-b}{d}-(\Irr(S')\backslash \{ b\})$$
          $\blacktriangleright$ Case 3 : if $b=1[d]$, then $(b-1,0)$ is in $L$, but not in $R$. Let us consider $\sigma' : (x,y)\to (b-1-x,-y)$. We have $d<b$.  If $d<a$, we argue as in Case 2. \\
          \esp Else :   $\PF(S')=\psi(\max(L\cap A))$, where $A=\inter{b-d}{b-1}\times \inter{-(a-1)}{-1}$ ( see Proposition \ref{prop:3}). Yet, $\sigma'(L\cap A)=L\cap A'$, where $ A'=\inter{0}{d-1}\times \inter{1}{a-1}$. But,  $\Irr(S')=\psi(\min(L\cap A'_2))$, where $A'_2= \inter{1}{d}\times \inter{0}{a-1}$ ( see Proposition \ref{prop:2}). Now $(d,0)$ is a minimal point of $L\cap A'_2$ and there is no point in $L\cap A'$ with $x=0$, so : $\min(L\cap A')=\min(L\cap A'_2)\backslash \{ (d,0)\}$ and we conclude as above...  
      \end{mademo}

   \newpage

     $\bullet$ In the previous proof, we have used the half-turn $u\to v-u$ in $\Z^2$, with $v=(b,-1)$ or $v=(b-1,0)$. Now, we will use it with $v=(b-1,-1)$.
      
     \begin{montheo}[ The  case $d<\min(a,b)$ ].  \label{theo:7} \\
      let $a,b,d$ be three pairwise coprime positive integers, such that $d<\min(a,b)$ and $d$ does not divide $ab-a-b$. Let $m\in \inter{1}{d-1}$ such that  $ am+b=0[d]$.  \\
      \esp    We denote $\alpha = \frac{m}{d}$, $[0,a_1,a_2,\cdots,a_r,1]$ its  CFE  and $(p_k/q_k)_k$ its convergents.\\
      \esp We denote $\tau=\frac{am+b}{d}$, $\beta=1-\{  \frac{b-1+m}{d}\}$ and  $(b_1,b_2,\cdots,b_{r})$ the   $\alpha$-numeration of $\beta$.\\
      \esp  We denote $B=\{  i\in \inter{1}{\lfloor r/2\rfloor }, b_{2i}\not = 0 \}$ ,  $\mu_i=\tau q_i-ap_i$ for all $i\in \inter{0}{r}$ and  :
      $$f=ab-a-b=f(<a,b>) \esp ;\esp f_1 = \frac{f}{d}-a(1-\beta) \esp ; \esp f_2=  \frac{f}{d}-\frac{b}{d}\sum_{i=1}^rb_iq_{i-1}$$
   \esp Then : 

   $$ \PF\left(\frac{<a,b>}{d}\right)= \{ f_1,f_2 \}\cup  \left\{ \frac{f}{d}+a\beta-j\mu_{2k-1}- \sum_{i=1}^{2k-1}b_i\mu_{i-1},\esp j\in\inter{0}{b_{2k}-1}, k\in B\right\}$$

     $$ t\left(\frac{<a,b>}{d}\right)
    =  2+  \sum\limits_{k=1}^{\lfloor r/2 \rfloor}b_{2k} $$

  \end{montheo}

  \vspace*{0.5cm}

  \textbf{Remark 1 : } if $d$ does not divide $f$ and $d<a<b$, then the type of $\frac{<a,b>}{d}$ is at least equal to 2, so $\frac{<a,b>}{d} $ is not symmetric.\\

    \begin{mademo} let $S'=\frac{<a,b>}{d}$.   
      We use Lemma \ref{lem:7} ( or Proposition \ref{prop:3}) : 
     $$\PF(S')= \psi(\max(R\cap L ))$$
     \esp where  $\psi : \begin{cases} \Z^2\to \frac{1}{d}\Z \\ (x,y)\to \frac{ax+by}{d} \end{cases} \esp ; \esp L=\psi^{-1}(\Z)$ and  $R= \inter{b-d}{b-1}\times \inter{-d}{-1} $.  \\
     \esp But, we can describe $R\cap L$ with the finite sequence $(b-1-g(n),-n-1)_{n\in\inter{0}{d-1}}$, where $g(n)$ is the integer of $\inter{0}{d-1}$ such that $a(b-1-g(n))-b(n+1)=0[d]$. We obtain : $\forall n\in\inter{0}{d-1},$
     $$a(b-1-g(n))-b(n+1)=0[d] \Leftrightarrow ag(n)=a(b-1+m+mn)[d] \Leftrightarrow g(n)=b-1+m+mn [d]$$
     \esp Indeed, $d$ and $a$ are coprime. So : 
     $$   \forall n\in\inter{0}{d-1} , \esp  g(n)= d \left\{ \frac{b-1+m+mn}{d} \right\}$$
     \esp We remark that :
     $$ g(0)=0 \Leftrightarrow b-1+m=0[d] \Leftrightarrow a(b-1+m)=0 [d] \Leftrightarrow ab-a-b=0[d]$$
  \esp Since $f\not = 0 [d]$, then $g(0)\not = 0$ and :  
     $$\max((b-1-g(n),-n-1)_{n\in\inter{0}{d-1}})= \sigma ( \min((g(n),n)_{n\in\inter{0}{d-1}}))$$
     \esp Where $\sigma : \begin{cases} \Z^2\to \Z^2 \\ (x,y)\to (b-1-x,-1-y) \end{cases}$.\\
           \esp We use now Lemma \ref{lem:20} ( and Remark 3 below it), that gives minimal points of $(\{ n\alpha -   \beta \},n)_{n\in\inter{0}{d-1}}$, where $\alpha= \frac{m}{d}$ and $\beta=\{ - \frac{b-1+m}{d}\}=1-\{  \frac{b-1+m}{d}\} \in ]0,1[$, since $d$ does not divide $b-1+m$ ( see above). \\
    \esp For a general $n$ ( corresponding to a minimal point of $(\{ n\alpha -   \beta \},n)_{n\in\inter{0}{d-1}}$), we obtain the PF-number of $S'$ : 
    $$\frac{f-bn}{d}-a\{n\alpha-\beta \}$$        
    \esp For $n=0$, we have $g(0)=d(1-\beta)$ ( this corresponds to the trivial minimal point $(1-\beta,0)$). This gives the PF-number for $S'$ : $f_1=\frac{f}{d}-a(1-\beta)$.\\
    \esp For $n$, whose $\alpha$-numeration is $(b_k)_k$ ( this means that $n$ is the integer in $\inter{0}{d-1}$, such that $\{ n\alpha \}=\beta$), we obtain $\{ n\alpha-\beta\}=0$ and so this gives the PF-number for $S'$ : $f_2=\frac{f-nb}{d}$, where $n=\sum\limits_{i=1}^rb_iq_{i-1}$.\\
    \esp For the general case $n$, whose $\alpha$-numeration is $(b_1,\cdots,b_{2k-1},j)$, with $j\in\inter{0}{b_{2k}-1}$ and $k\in B$, we obtain : 
    $$ \frac{f}{d}-\frac{b}{d}\left( \sum_{i=1}^{2k-1}b_iq_{i-1}+jq_{2k-1}\right)-a\left( \sum_{i=1}^{2k-1}b_i(q_{i-1}\alpha-p_{i-1})+j(q_{2k-1}\alpha-p_{2k-1})\right)+a\beta$$
    \esp Indeed, $\{n\alpha - \beta \}=\{ n\alpha \}-\beta$, becasue $\{ n\alpha \}<\beta$ for these values of $n$. ( see Lemma \ref{lem:20}). This gives our formula...\\
    \esp For the type $t(S')$, we remark that $\psi$ is injective over $R$, if $d<\min(a,b)$ ( see Proposition \ref{prop:3}).
     \end{mademo}

\vspace*{0.5cm}

$\bullet$ To be complete, we end with the case $d>\min(a,b)$, that is somwhat different :  in the proof we use an half-turn $u\to v-u$, with $v=(b,0)$.

  \begin{montheo} [ The case $d>\min(a,b)$]. \label{theo:8}\\ 
    let $a,b,d$ be 3 pairwise coprime integers, such that $1<\min(a,b)<d$ and $d$ does not divide $ab-a-b$.\\
    \esp Let $m\in \inter{1}{d-1}$ such that  $ am+b=0[d]$.      We denote $\alpha = \frac{m}{d}$, $[0,a_1,a_2,\cdots,a_r,1]$ its  CFE  and $(p_k/q_k)_k$ its convergents. \\
    \esp We denote $\tau=\frac{am+b}{d}$ ( a positive integer), $\beta=\{ a\alpha \}=1-  \left\{\frac{b}{d}\right\},(b_1,b_2,\cdots,b_s)$ the $\alpha$-numeration of $\beta$, where $b_s>0$ and $\mu_i=\tau q_i-ap_i$ for all $i\in \inter{0}{r}$.\\
    \esp We also denote $B=\{  i\in \inter{1}{\lfloor s/2\rfloor }, b_{2i}\not = 0 \}$ and $k_0= \begin{cases} (s+1)/2 \text{  if } B=\emptyset \\ \min(B) \text{ else } \end{cases} $  ( $s$ is odd in the first case).  Last notations :

  $$\F_0 = \left\{ a\left\lfloor \frac{b}{d}\right \rfloor - \tau \right\} \text{ if } \alpha < \beta \esp ; \esp \F_0=\emptyset \text{ else }$$
  $$ \F_1= \left\{ a\left\lfloor \frac{b}{d}\right \rfloor - \mu_{2k-2}-j\mu_{2k-1} ,  j\in\inter{1}{a_{2k}}, k\in \inter{1}{k_0-1}\right\}$$
  
   $$ \F_2= \left\{ a\left\lfloor \frac{b}{d}\right \rfloor +  a -j\mu_{2k-1}- \sum_{i=1}^{2k-1}b_i\mu_{i-1}, j\in\inter{0}{b_{2k}-1}, k\in B \right\}$$

  $\blacktriangleright $ Case 1 : if $ a<d<b $
      $$ \PF\left(\frac{<a,b>}{d}\right)= \F_0\sqcup \F_1\sqcup \F_2$$
     $$   t\left(\frac{<a,b>}{d}\right)    =   \sum_{k=1}^{k_0-1}a_{2k}+\sum\limits_{k=k_0}^{\lfloor s/2 \rfloor}b_{2k}+ \textbf{1}_{\alpha<\beta} $$
     
   $\blacktriangleright $ Case 2 : if $ d>\max(a,b) $
      $$ \PF\left(\frac{<a,b>}{d}\right)=  \F_2 \esp ; \esp  t\left(\frac{<a,b>}{d}\right)  =\sum\limits_{k=k_0}^{\lfloor s/2 \rfloor}b_{2k}$$

  \end{montheo}

  \textbf{Remark : } $(b_k)_k$ is also the $\alpha$-numeration of $a$. Moreover, if $d>\max(a,b)$, we can not have $b_{2i}=0$ for all $i$, because $\PF(<a,b>/d)$ can not be empty ! \\
 \esp If $a<d<b$ and $\alpha>\beta$, that is if $b_1=1$ and $b_2>0$, then $k_0=1$ and $\F_1$ is empty.\\

  \begin{mademo}
   \esp We suppose that $a<b$, so $a<d$. This proof will be similar to the proof of Theorem \ref{theo:7}. We mention the notations used in this one : $S'=\frac{<a,b>}{d}$ and with Proposition \ref{prop:3} : 
     $$\PF(S')= \psi(\max(R\cap L )) \esp \text{  where } \esp \psi : \begin{cases} \Z^2\to \frac{1}{d}\Z \\ (x,y)\to \frac{ax+by}{d} \end{cases}$$
     \esp and $ L=\psi^{-1}(\Z) \esp ; \esp R= \inter{b-\delta}{b-1}\times \inter{-(a-1)}{-1} $,  where $\delta=\min(b-1,d)$.  \\
     
     \esp The difference with the proof of Theorem \ref{theo:7} lies in $R$ and in the parametrization of $R\cap L$ : we will use the sequence $(b-h(n),-n)_{n\in\inter{1}{a-1}}$, where $h(n)$ is the integer of $\inter{1}{d}$ such that $a(b-h(n))-bn=0[d]$. But, in the case $d>\max(a,b)$, we will ask that $h(n)<b$ !! \\
  $\blacktriangleright $ Case 1 : $a<d<b$.    We obtain : 
     $$ \forall n\in\inter{1}{a-1}, \esp ah(n)=ab-nb=(n-a)am[d]\Leftrightarrow h(n)=nm-am[d]$$
     \esp Yet $a<d$ and $(b,-a)$ is a point of $L$, so there is no point of $L$ among the $(b,-n)_{n\in\inter{1}{a-1}}$. Hence, $h(n)\not=0[d]$ for $n\in\inter{1}{a-1}$. So :
     $$   \forall n\in\inter{1}{a-1},\esp h(n)=d\left\{ (n-a)\alpha  \right\}=d\left\{ n\alpha-\beta  \right\} $$
     \esp So, $L\cap R$ is the image of $U=(h(n),n)_{n\in\inter{1}{a-1}}$ by the map $\sigma : (x,y)\to (b-x,-y)$ and :
     $$\max(L\cap R)= (b,0)-\min(U)$$
     \esp If $(h(n),n)$ is a minimal point of $U$, then we obtain the following PF number of $S'$  :
     $$ \psi(b-h(n),-n)= \frac{b(a-n)}{d}- a\{ n\alpha - \beta \}$$
     \esp Now, we use Lemma \ref{lem:20} and first remark that $n=a$, that is $n=(b_1,\cdots,b_s)_{\alpha}$, is excluded. Secondly, we     distinguish two kinds of minimal point of $(\{ n\alpha - \beta \},n)_{n\in \inter{1}{a-1}}$ : \\
     \esp - those that satisfy $\{ n\alpha \}<\beta$,  which corresponds to ( see Proof of Theorem \ref{theo:7}) : 
     $$n=q_{2k-2}+ j q_{2k-1}  \esp ; \esp \{ n\alpha \}= \delta_{2k-2}-j\delta_{2k-1}, \esp  j\in\inter{1}{a_{2k}}, k\in \inter{1}{k_0-1}$$
     \esp For these minimal points, we have $\{ n\alpha - \beta \}=\{ n\alpha \}+1-\beta= \{ n\alpha \}+ \frac{b}{d}-\lfloor b/d \rfloor$. We deduce the formula for the elements of $\F_1$. \\
      \esp - those that satisfy $\{ n\alpha \}\geqslant \beta$,  which corresponds to :
     $$ n=\sum_{i=1}^{2k-1}b_iq_{i-1}  \esp ; \esp \{ n\alpha \}= \sum_{i=1}^{2k-1}b_i(q_{i-1}\alpha - p_{i-1}), \esp j\in\inter{0}{b_{2k}-1}, k\in \inter{k_0}{\lfloor s/2 \rfloor}$$
       \esp For these minimal points, we have $\{ n\alpha - \beta \}=\{ n\alpha \}-\beta= \{ n\alpha \}+ \frac{b}{d}-\lfloor b/d \rfloor-1$. We deduce the formula for the elements of $\F_2$. \\
    \esp Now, if we want to conclude for PF-numbers of $S'$, we just have to remark that Lemma \ref{lem:20} also gives a minimal point of $U$ for $n=1$, which is $(d\{ \alpha-\beta\},1)$. If $\alpha > \beta$, that point is in $\F_2$ and already counted. Else, it gives a new PF-number : $ a\lfloor b/d \rfloor - \tau$.\\
    \esp   For the type, we just have to count elements of $\PF(S')$, keeping in mind that the sets $\F_0,\F_1,\F_2$ are pairwise disjoints and that the parametrization of their elements is injective.\\
     $\blacktriangleright $ Case 2 : if $d>\max(a,b)$. We keep above arguments with the additional condition that $h(n)<b$ : this gives the condition $\{ n\alpha - \beta \}< 1-\beta$, that is $\{n\alpha \}\geqslant \beta$. We conclude with Lemma \ref{lem:20}.    
    \end{mademo}




\textbf{Remark   : }  unlike in the case $d<\min(a,b)$, we can find $a,b,d$ such that $a<d<b$, $d$ does not divide $f$ and $\frac{<a,b>}{d}$ is symmetric : for example, with $a=11, b=89$ and $d=20$, we obtain $\frac{<a,b>}{d}=<5,11>$. \\

 \subsection{the Frobenius number}

  $\bullet$ Once we can compute PF-numbers of a numerical semigroup, we know that  the Frobenius number is the greatest of them. In the expressions of PF-numbers that appear in Theorems \ref{theo:7} and \ref{theo:8}, we have two parameters $j$ and $k$ that have  values in a finite interval of integers. Moreover, we remark that, if $k$ is fixed, then these expressions on $j\to H_k(j)$ are affine, so monotonous. Hence, the maximum is obtained for an extremal value of $j$ : that is to say $0$ or $1$ or $a_{2k}$ or $b_{2k}-1$ ( if $b_{2k}\not = 0$ !). Remark that $(b_k)_k$ is not the same for Theorem \ref{theo:7} and \ref{theo:8}. \\
 \esp   We denote $B=\{ i\in \inter{1}{\lfloor r/2 \rfloor}, b_{2i}\not = 0 \}$ and obtain  : 
  
  $\blacktriangleright $ if $d<\min(a,b)$
    $$ f\left(\frac{<a,b>}{d}\right)=\max\left( f_1,f_2,  \frac{f}{d}+a\beta-j\mu_{2k-1}- \sum_{i=1}^{2k-1}b_i\mu_{i-1}, \esp  j=0, j =b_{2k}-1, k\in B\right)$$
    
   $\blacktriangleright $ if $a<d<b$ 
   $$ f\left(\frac{<a,b>}{d}\right)=\max(  a\left\lfloor \frac{b}{d}\right \rfloor - \tau , a\left\lfloor \frac{b}{d}\right \rfloor - \mu_{2k-2}-j\mu_{2k-1} , \esp j=1,j=a_{2k}, k\in \inter{1}{k_0-1},  $$
$$ ,  a\left\lfloor \frac{b}{d}\right \rfloor +  a -j\mu_{2k-1}- \sum_{i=1}^{2k-1}b_i\mu_{i-1},\esp j=0, j=b_{2k}-1, k\in B )$$
   
   $\blacktriangleright $ if $d>\max(a,b)$ : 
   $$ f\left(\frac{<a,b>}{d}\right)=\max\left(  a\left\lfloor \frac{b}{d}\right \rfloor +a-j\mu_{2k-1}- \sum_{i=1}^{2k-1}b_i\mu_{i-1},\esp j=0, j =b_{2k}-1, k\in B\right)$$
   
  \esp Now, for which $k$ do we obtain a maximum ? To answer this question, we think of our representation ( $\psi : (x,y)\to (ax+by)/d$) of $S=\frac{<a,b>}{d}$ in $\Z^2$ and remember that  the points representing $\PF(S)$ ( see Lemma \ref{lem:8} and Proposition \ref{prop:3}) can be parametrized by an $x$-decreasing sequence $(u_k)_{k\in\inter{1}{t}}$  that is order-concave in $\Z^2$ ( with the partial product order), which means that : $(u_k-u_{k-1})_k$ is decreasing. But Theorem \ref{theo:7} and \ref{theo:8} give the PF-numbers of $S$, say $(\psi(u_i))_i$, respecting the previous order ( see proof of these theorems). Consequence : if we denote $(h_i)_i$ that sequence of integers, then $(h_i-h_{i-1})_i$ is decreasing, for $\psi$ is increasing ( $a$ and $b$ are positive). So, $h_i$ is  the Frobenius number of $S$ for the unique $i$ such that $h_i-h_{i-1}>0$ and $h_{i+1}-h_i<0$. \\
  
  \esp Yet, the slope of our affine functions $H_k$ is $-\mu_{2k-1}$, see notations in Theorem \ref{theo:7} and \ref{theo:8}. So, $(\mu_{2k-1})_k$ is increasing and  we will discuss the sign of $\mu_{2k-1}$. We recall that :
  $$\forall i\in\inter{0}{r}, \mu_i=\tau q_i-ap_i=aq_i\left(\alpha+\frac{b}{ad}-\frac{p_i}{q_i}\right)$$
  \esp where $a,b,d$ are three positive pairwise coprime integers, $\alpha=\frac{m}{d}$, $m\in\inter{1}{d-1}$ coprime with $d$, such that $am+b=0[d]$  and $(p_i/q_i)_{i\in\inter{0}{r}}$ are the usual convergents of $\alpha$. So, for all integer $k\in\inter{1}{\lfloor (r+1)/2 \rfloor}$ :
  $$\mu_{2k-1}>0  \esp \Leftrightarrow \esp \frac{p_{2k-1}}{q_{2k-1}}-\alpha < \frac{b}{ad}  \esp \Leftarrow \esp q_{2k-1}^2> \frac{ad}{b}$$
   \esp But, we will need more precise notations...
   
   \newpage


  \begin{montheo}  \label{theo:9}  let $a,b,d$ be three pairwise coprime integers.\\
    \esp Let $m\in \inter{1}{d-1}$ such that  $ am+b=0[d]$.      We denote $\alpha = \frac{m}{d}$, $[0,a_1,a_2,\cdots,a_r,1]$ its  CFE  and $(p_k/q_k)_k$ its convergents.  We denote $\tau=\frac{am+b}{d}, \mu_i=\tau q_i-ap_i$ for $i\in \inter{-1}{r}$ , $f=f(<a,b>)=ab-a-b$.
    
    $$ A_-=\left\{ i, b_{2i}\not = 0 \text{ and } \mu_{2i-1}<0 \right\} \esp ; \esp  A_+=\left\{ i, b_{2i}\not = 0 \text{ and } \mu_{2i-1}>0 \right \}$$
     $$k_1=\max(A_-) \esp ; \esp k_2=\min(A_+)$$

    $\blacktriangleright$ Case 1 : if $d$ divides $f$ then $f\left(\frac{<a,b>}{d}\right)=\frac{f}{d}$. \\
    $\blacktriangleright$ Case 2 : if $d<\min(a,b)$, then we denote   $\beta=1-\left\{\frac{b+m-1}{d}\right\}=(b_1,\cdots,b_r)_{\alpha}$. 
    
    $$ \text{ if } A_-=\emptyset \esp , \esp  f\left(\frac{<a,b>}{d}\right)=\frac{f}{d}+a\beta-a$$
    $$ \text{ if } A_+=\emptyset \esp , \esp  f\left(\frac{<a,b>}{d}\right)=\frac{f}{d}-\frac{b}{d}\sum_{i=1}^rb_iq_{i-1}$$
    $$ \text{ else } \esp , \esp  f\left(\frac{<a,b>}{d}\right)=\frac{f}{d}+a\beta-\min\left(-\mu_{2k_1-1}+\sum_{i=1}^{2k_1}b_i\mu_{i-1},\sum_{i=1}^{2k_2-1}b_i\mu_{i-1}\right)$$

    $\blacktriangleright$ Case 3 : if $a<d<b$, we denote  $1-\left\{\frac{b}{d}\right\}=(b_1,\cdots,b_r)_{\alpha}$.
     $$ \gamma = \min \left( i\in\inter{1}{\lfloor (r+1)/2 \rfloor },+\mu_{2i-1}>0 \right) \esp ; \esp k_0= \begin{cases} k_0=\lfloor r/2\rfloor \text{ if } b_{2i}=0 \text{  for all integer } i \\ \min(i,b_{2i}\not = 0) \text{ else} \end{cases}$$

    $$ \text{ if } \gamma < k_0 \esp , \esp  f\left(\frac{<a,b>}{d}\right)=a\left\lfloor \frac{b}{d}\right\rfloor-\mu_{2\gamma-2}$$
    
       $$  \text{ if } \gamma > k_0 \esp , \esp  f\left(\frac{<a,b>}{d}\right)=a\left\lfloor \frac{b}{d}\right\rfloor+a-\min\left(-\mu_{2k_1-1}+\sum_{i=1}^{2k_1}b_i\mu_{i-1},\sum_{i=1}^{2k_2-1}b_i\mu_{i-1}\right)$$
    
       $$ \text{ if } \gamma = k_0 \esp , \esp  f\left(\frac{<a,b>}{d}\right)=a\left\lfloor \frac{b}{d}\right\rfloor-\min\left(\mu_{2k_0-2},-a+\sum_{i=1}^{2k_0-1}b_i\mu_{i-1}\right)$$
     
  $\blacktriangleright$ Case 4 : if $d>\max(a,b)$, , we denote  $1-\left\{\frac{b}{d}\right\}=(b_1,\cdots,b_r)_{\alpha}$. Then :
    $$ \text{ else } \esp , \esp  f\left(\frac{<a,b>}{d}\right)=a\left\lfloor \frac{b}{d}\right\rfloor+a-\min\left(-\mu_{2k_1-1}+\sum_{i=1}^{2k_1}b_i\mu_{i-1},\sum_{i=1}^{2k_2-1}b_i\mu_{i-1}\right)$$

 \end{montheo}
  
 \vspace*{1cm}

  \textbf{Question : } could we find a simpler formula for $f(\frac{<a,b>}{d})$ ? If we refer to our proof and notations of Theorem \ref{theo:4}, we can claim that : 
   $$f\left(\frac{<a,b>}{d}\right)=\max\{ k\in\N\backslash (a \inter{0}{b-1}), \lfloor k\alpha \rfloor =  \lfloor k\alpha' \rfloor \} $$
   \esp But, it is not so easy to compute this max...\\

 \newpage
 
  \begin{mademo}
    $\blacktriangleright$ Case 1 : has already been shown at Theorem \ref{theo:6}.\\
    
     $\blacktriangleright$ Case 2 : we use and refer to Theorem \ref{theo:7}. We remind the result of this Theorem : if we denote $B=\{  i\in \inter{1}{\lfloor r/2\rfloor }, b_{2i}\not = 0 \}, \esp f_1 = \frac{f}{d}-a(1-\beta) \esp ; \esp f_2=  \frac{f}{d}-\frac{b}{d}\sum\limits_{i=1}^rb_iq_{i-1}$, then : 

   $$ \PF\left(\frac{<a,b>}{d}\right)= \{ f_1,f_2 \}\cup  \left\{ \frac{f}{d}+a\beta-j\mu_{2k-1}- \sum_{i=1}^{2k-1}b_i\mu_{i-1},\esp j\in\inter{0}{b_{2k}-1}, k\in B\right\}$$
   
     $\blacktriangleright\blacktriangleright$ subcase 1 : if $A_-=\emptyset$, then $\mu_{2k-1}>0$ for all $k$ in $B$, with notations of Theorem \ref{theo:7}. So, the Frobenius is obtained for the "first" pseudo-Frobenius in the list, namely $f_1$.\\
     $\blacktriangleright\blacktriangleright$ subcase 2 : if $A_+=\emptyset$, then $\mu_{2k-1}<0$ for all $k$ in $B$, with notations of Theorem \ref{theo:7}. So, the Frobenius is obtained for the "last" pseudo-Frobenius in the list, namely $f_2$. \\ 
      $\blacktriangleright\blacktriangleright$ subcase 3 : if both $A_-$ and $A_+$ are non empty, then the Frobenius is obtained for one of the two following pseudo-Frobenius, listed in Theorem \ref{theo:7} :
      $$ f_-=\frac{f}{d}+a\beta- ( b_{2k_1}-1)\mu_{2k_1-1}-\sum_{i=1}^{2k_1-1}b_i\mu_{i-1} \esp ; \esp f_+=\sum_{i=1}^{2k_2-1}b_i\mu_{i-1} $$
      \esp $f_-$ being the greater of PF-numbers parametrized by $k$ such that $\mu_{2k-1}<0$ ( $k=k_1$ and $j=b_{2k-1}-1$)  and $f_+$  being the greater of PF-number parametrized by $k$ such that $\mu_{2k-1}>0$ ( $k=k_2$ and $j=0$). Our formula is a direct consequence. \\
      
      $\blacktriangleright$ Case 3 : we use and refer to Theorem \ref{theo:8}. We remind the principle notations and result in that case :  
       
  $$\F_0 = \left\{ a\left\lfloor \frac{b}{d}\right \rfloor - \tau \right\} \text{ if } \alpha < \beta \esp ; \esp \F_0=\emptyset \text{ else }$$
  $$ \F_1= \left\{ a\left\lfloor \frac{b}{d}\right \rfloor - \mu_{2k-2}-j\mu_{2k-1} ,  j\in\inter{1}{a_{2k}}, k\in \inter{1}{k_0-1}\right\}$$
  
   $$ \F_2= \left\{ a\left\lfloor \frac{b}{d}\right \rfloor +  a -j\mu_{2k-1}- \sum_{i=1}^{2k-1}b_i\mu_{i-1}, j\in\inter{0}{b_{2k}-1}, k\in \inter{k_0}{\lfloor s/2 \rfloor}\right\}$$

      $$ \PF\left(\frac{<a,b>}{d}\right)= \F_0\sqcup \F_1\sqcup \F_2$$

     $\blacktriangleright\blacktriangleright$ subcase 1 : if $\gamma < k_0$, then the Frobenius of $S$ is among the elements of $\F_1$, for $\mu_{2k-1}>0$ for $k\geqslant k_0-1$. It is obtained for one of the following pseudo-Frobenius listed in $\F_1$ : 
     $$ f_-= a\left\lfloor \frac{b}{d}\right\rfloor- \mu_{2\gamma -4}-a_{2\gamma -2}\mu_{2\gamma -3}  \esp ; \esp f_+=  a\left\lfloor \frac{b}{d}\right\rfloor- 0- \mu_{2\gamma-2}$$
     \esp $f_-$ being the greater of PF-numbers parametrized by $k$ such that $\mu_{2k-1}<0$ ( $k=\gamma -1$ and $j=a_{2\gamma -2}$)  and $f_+$  being the greater of PF-number parametrized by $k$ such that $\mu_{2k-1}>0$ ( $k=\gamma $ and $j=0$). But, $f_+=f_-$ and we deduce our formula. \\
      $\blacktriangleright\blacktriangleright$ subcase 2 : if $\gamma > k_0$, then the Frobenius of $S$ is among the elements of $\F_2$, becasue $\mu_{2k-1}<0$ for $k= k_0$. Now, the arguments are the same as in Case 1 subcase 3...\\
       $\blacktriangleright\blacktriangleright$ subcase 3 : if $\gamma = k_0$, then $\mu_{2k_0-3}<0$ and $\mu_{2k_0-1}>0$. So, the Frobenius of $S$ is one of the following integers :
       $$ f_-= a\left\lfloor \frac{b}{d}\right\rfloor-\mu_{2k_0-2} \esp ; \esp f_+= a\left\lfloor \frac{b}{d}\right\rfloor-a- \sum_{i=1}^{2k_0-1}b_i\mu_{i-1}$$
        \esp $f_-$ being the greater of PF-numbers parametrized by $k$ such that $\mu_{2k-1}<0$ ( $k=k_0 -1$ and $j=a_{2k_0 -2}$)  and $f_+$  being the greater of PF-number parametrized by $k$ such that $\mu_{2k-1}>0$ ( $k=k_0 $ and $j=0$).\\
        
       $\blacktriangleright$ Case 4 : the arguments are the same as in Case 3 subcase 2, for we only consider $\F_2$...     
     \end{mademo}

  \newpage

      \subsection{genus}
      
      \begin{montheo} \label{theo:10}
        Let $a,b,d$ be  three pairwise coprime positive integers and $m$ an integer such that $b=am$ mod $d$.
        $$ g\left(\frac{<a,b>}{d}\right)
       =   \frac{g(<a,b>)}{d}+ \frac{1}{2}\left[C(\alpha,\beta,\nu)-1- \beta\nu\right]$$
      
      \esp where $\alpha = \left\{ \frac{m}{d} \right \}, \beta =  \left\{ \frac{b-1}{d} \right \},\nu= d \left\{ \frac{a-1}{d} \right \}$ and $C=\#\{ k\in\inter{0}{\nu}, \{ k\alpha \}\leqslant \beta \}$.\\
       ( see   Proposition \ref{prop:11} in \ref{subsec:counting} for an expression of $C$).\\
      \esp In particular, we obtain : 
      $$ a=1  \esp [d] \text{ or } b =1  \esp [d]  \Rightarrow   g\left(\frac{<a,b>}{d}\right)
       =   \frac{g(<a,b>)}{d}= \frac{(a-1)(b-1)}{2d}$$
      
      \end{montheo}
      
      \begin{mademo}
     \esp According to Proposition \ref{prop:1}, we have to count points in $R\cap L$, where :
      $$ R=\inter{1}{b-1}\times \inter{-(a-1)}{-1} \esp ; \esp L=\{ (x,y)\in\Z^2, ax+by=0 \text{  mod } m \}$$
     \esp But, we know that $L$ contains exactly one point one every " horizontal" or " vertical" segment of length $d$. So, if we denote $N=\# ( R\cap L)$, we have :
      $$ N= (b-1)\left\lfloor \frac{a-1}{d} \right\rfloor + d\left\{ \frac{a-1}{d} \right\}\left\lfloor \frac{b-1}{d} \right\rfloor + N' $$
    \esp  where the first term counts points of $L$ in the sub-rectangle $\inter{1}{b-1}\times \inter{-d\left\lfloor \frac{a-1}{d} \right\rfloor}{-1}$. The second term counts points of $L$ in the sub-rectangle $\inter{1}{d\left\lfloor \frac{b-1}{d} \right\rfloor }\times \inter{-(a-1)}{-d\left\lfloor \frac{a-1}{d} \right\rfloor-1}$. And finally, $N'$ counts points of $L$ in the subrectangle : $\inter{d\left\lfloor \frac{b-1}{d} \right\rfloor +1}{b-1}\times \inter{-(a-1)}{-d\left\lfloor \frac{a-1}{d} \right\rfloor-1}$. But, since $L$ is invariant via translation of $(d,0)$ and $(0,d)$, $N'$ is the number of points of $L$ in the rectangle  
      $\inter{1}{(b-1)\%d}\times \inter{ -(a-1)\%d}{-1}$, where $\%$ denotes the remainder of the euclidean division.\\
      \esp To express $N'$, we use the following parametrization of points of $L$ in a rectangle  :
      $$L\cap\inter{0}{d-1}\times \inter{-(a-1)\% d}{-1}=\{ (x_k,-k) , k\in\inter{1}{(a-1)\% d}, \text{ such that } ax_k-bk=0 \text{ mod } d \}$$
      \esp but :
      $$ ax_k-bk=0 \text{ mod }d \Leftrightarrow 
      ax_k=amk \text{ mod } d \Leftrightarrow  x_k=mk \text{ mod } d $$
      \esp indeed $a$ and $d$ are coprime. But, $x_k\in \inter{0}{d-1}$, so $x_k = d \left\{ \frac{km}{d}\right\}$. Now, we count points with $x_k\in\inter{1}{(b-1)\%d}$, so the condition on $x_k$ is expressed by $0< \{ \frac{km}{d} \} \leqslant \{ \frac{b-1}{d}\}$. It gives :
      $$ N'= C(\alpha,\beta,\nu)-1 $$
      \esp ( $-1$ because we omit $k=0$), with 
      $\alpha = \left\{ \frac{m}{d} \right \}, \beta =  \left\{ \frac{b-1}{d} \right \},\nu= d \left\{ \frac{a-1}{d} \right \}$ and $C$ defined as above.\\
      \esp If $d$ divides $a-1$ or $b-1$, then $\beta$ or $\nu$ is null, so $C(\alpha,\beta,\nu)=1$ and $\beta\nu=0$. 
      \end{mademo}

\vspace*{0.5cm} 
      
      \textbf{Remark 1 : } that expression of $g$ is not symmetric in terms of $a$ and $b$, even if the product $\beta\nu$ is. But it proves that :
      $$ C(\alpha,\beta,\nu)=C(\alpha',\beta',\nu')$$
      where $\alpha'=\{\frac{m'}{d}\}$ and $m'$ is the inverse of $m$ mod $d$, $\beta'= \frac{\nu}{d}, \nu'=d\beta$.\\
      
      \textbf{Remark 2 : } $C(\alpha,\beta,\nu)-1$ is the number of integer $k$ such that $k\in\inter{1}{\nu}$ and $\{ k\alpha \}\leqslant \beta$. But $\beta$ is the probability that a real in $[0,1[$ lies in $[0,\beta]$. So, if we choose randomly $\nu$ reals in $[0,1[$, then $\beta\nu$ is the expectation of the number of these reals that lie in $[0,\beta]$.\\
      \esp Hence, the term $C(\alpha,\beta,\nu)-1-\beta\nu$ can be interpreted as a "signed deviation from the random mean value"... \\

      \textbf{Remark 3 : } this expression of $g$ and especially that of $f$ is too intricate to help us about Wilf's property...

     \newpage

    \subsection{comparison of $t$ and $e$, Wilf's property}

$\bullet$ We have already mentionned that minimal generators and pseudo-Frobenius numbers of a numerical semigroup $S$ have similar properties : they are respectively minimal and maximal sets of $S\backslash \{ 0 \}$ and $\Z\backslash S$ for the order induced by $S$. That kind of symmetry could imply some relation between their finite cardinality. Unfortunately, this is not the case : we have $t(S)<e(S)$, if $e(S)\leqslant 3$, but for $e(S)\geqslant 4, t(S)$ can be as large as we want ( \textbf{[7]} \nameref{biblio}). Nonetheless, we know that $t(S)<m(S)$, considering rests modulo $m(S)$. \\
\esp So, this is remarkable that $t(S)<e(S)$ when $S$ is a  proportionally modular semigroup, that is to say $S=\frac{<a,b>}{d}$, for some  three pairwise coprime integers $a,b,d$ :

     \begin{montheo} \label{theo:11} if $a,b,d$ are 3 pairwise coprime positive integers such that $d\not\in <a,b>$, then : 
     
     $$t\left(\frac{<a,b>}{d}\right)<  e\left(\frac{<a,b>}{d}\right)$$
     \esp So, Wilf's property holds for these numerical semigroups !    
      
  \end{montheo}

 \vspace*{1cm}
  
  \textbf{Remark : } before presenting a proof, we want to underline that, depending on the respective positions of $a,b,d$, this proof will be either obvious or laborious ( never difficult nonetheless). We are however convinced that this inequality, between $e$ and $t$, could be deduced from a deeper result about minimal points in a lattice of $\Z^2$. This will be perhaps the subject of a next article...\\
  \esp To make this 3-pages proof readable, here is a summary : since we use Theorem \ref{theo:5}, \ref{theo:6}, \ref{theo:7}  and  \ref{theo:8}, we must distinguish several cases.\\
  
  $\blacktriangleright$ Case 1 : if $d$ divides $a-1, b-1$ or $ab-a-b$ ( the frame of Theorem  \ref{theo:6}). Obvious case because $\Irr(S)$ and $\PF(S)$ are symmetric ( except for one element ), if $S=\frac{<a,b>}{d}$.\\
  
    $\blacktriangleright$ Case 2 : if $d<\min(a,b)$ ( the frame of Theorem  \ref{theo:7}), it is a consequence of definition of $\alpha$-numeration ( $b_i\leqslant a_i$ with notations of this Theorem).\\
    
     $\blacktriangleright$ Case 3 : if $a<d<b$, the proof is more technical and we will study two subcases  depending on $\alpha$-numeration of $a-1$ and $a$ :  a general subcase, where these numerations have the same length and an exceptional subcase. But, we keep some room to manoeuvre.\\
     
       $\blacktriangleright$ Case 4,5,6 : if $d>\max(a,b)$. Proofs are becoming longer, more intricate and we show that some exceptional subcases are impossible. Case 6 fills almost half  of the proof...\\

    \newpage

  \begin{mademo} \esp We denote $S=\frac{<a,b>}{d}$. We will use some previous Theorems and their notations, as well as our $\alpha$-numeration. We  refer the reader to the statements of these theorems and to section \ref{subsec:alphanum}.\\
  $\blacktriangleright$ Case 1 : we suppose that we are in one of the three " easy cases" of Theorem  \ref{theo:6} : if $d$ divides $ab-a-b$, then $t(S)=1$ and $e(S)\geqslant 2$. If, $d$ divides $a-1$ or $b-1$, then Theorem  \ref{theo:6} proves that $t(S)=e(S)-1$. \\
   
  $\blacktriangleright$ Case 2 : we  suppose that $d<\min(a,b)$. \\
  \esp   With Theorems  \ref{theo:5} and  \ref{theo:7}, we obtain :
  $$   e\left( \frac{<a,b>}{d}\right)=3+ \sum_{k=1}^{\lfloor r/2 \rfloor }a_{2k} \esp ; \esp  t\left(\frac{<a,b>}{d}\right)=2+\sum\limits_{k=1}^{\lfloor r/2 \rfloor}b_{2k}  $$
  \esp But, $\forall k\in \inter{1}{\lfloor r/2 \rfloor}, b_{2k}\leqslant a_{2k}$, so the inequality is true. \\
  \esp We notice, that $e(S)=t(S)+1$ if and only if : $ \forall k\in \inter{1}{\lfloor r/2 \rfloor}, b_{2k}= a_{2k}$. \\
  
  $\blacktriangleright$ Case 3 : we suppose that $a<d<b$.\\
  \esp We use Theorems  \ref{theo:5} and  \ref{theo:8} with the common notations $\alpha, (a_k)_k,k_0$. We also denote  $(b_k)_k$, the $\alpha$-numeration of $a$,  $\esp (N_k)_k$ the $\alpha$-numeration of $a-1$ and $s=\max(i,N_i\not = 0)$.\\
  $\blacktriangleright\blacktriangleright$ Subcase 1 : if $b_k=0$ for all $k>s$. 
   $$ e(S)-t(S)= 2+\nu- \chi_{\alpha < \beta}- c_s +\sum_{k=k_0}^{\lfloor (s-1)/2 \rfloor}(a_{2k}-b_{2k})$$
  
 \esp Where $c_s=b_s$ if $s$ is even and $0$ else.  The sum above is non negative ( see Case 2) and, since $(b_k)_k$ is the successor of $(N_k)_k$ in $E_{\alpha}$ for RLO, we have $N_s=b_s$ or $b_s-1$. Hence,  $\nu-c_s=0$, except if $s$ is even and $N_s=b_s-1$ ( see below). So, in the general case, we have $e(S)-t(S)\geqslant 1$. \\
  \esp Now, what does happen if $s$ is even and $N_s=b_s-1$ ? This  is verified if and only if $(b_k)_k$ is the lowest  element in $(E_{\alpha},RLO)$ such that $b_s$ is fixed and $b_k=0$ for $k>s$, that is to say $b_{s-2k}=a_{s-2k}$  and $b_{s-2k+1}=0$ for $k>1$ ( except $b_1=1$). In that case, we obtain $c_{s}=N_s+1=\nu+1$ and the sum is null so : $e(S)-t(S)=1-\chi_{\alpha<\beta}=1$, for $\beta<\alpha$ ( $b_1=1$ and $b_2>0$). So, we have $e(S)-t(S)=1$. \\

 $\blacktriangleright\blacktriangleright$ Subcase 2 : if   $(N_1,N_2,\cdots,N_{s})$  is the greatest element in $(E_{\alpha},RLO)$ ( see \ref{subsec:alphanum}) of that form, that is to say if $(N_k)_k=(a_1,a_2,a_3,\cdots,a_{s})$. Then $b_{s+1}=1$ and $b_{s-2k}=0$ for all $k\geqslant 0$ ( except $b_1=1$ if $s$ is odd) and $b_{s-2k+1}=a_{s-2k+1}$ if $k\in \inter{1}{\lfloor s/2 \rfloor}$ . Then 
        $$ e(S)-t(S)= 2+\nu- \chi_{\alpha < \beta}- c_{s+1}+\sum_{k=k_0}^{\lfloor (s-1)/2 \rfloor}(a_{2k}-b_{2k}) $$
 \esp  where $c_{s+1}=1$ if $s$ is odd and $c_{s+1}=0$ else.
 We make a distinction according to the parity of $s$ : \\
  $\blacktriangleright \blacktriangleright\blacktriangleright $ if $s$ is even : then $b_{2k}=0$ for all $k$, $c_{s+1}=0$ and $\nu=N_{s}=a_{s}$, so $e(S)-t(S)\geqslant 2$. \\  
   $\blacktriangleright \blacktriangleright\blacktriangleright $ if $s$ is odd : then $b_{2k}=a_{2k}$ for all $2k< s$, $c_{s+1}=1$ and $\nu=0$. In addition, $\alpha >\beta$, for $b_1=1$ and $b_2>0$ ( we remark that if $s=1$, then $b_2=1$). So $e(S)-t(S)=1$.\\
  
 $\blacktriangleright$ Case 4 : if $d>\max(a,b)$ and $x_0+b>d$.\\
 \esp We use again Theorem  \ref{theo:5} and  \ref{theo:8} with their notations. For the notation $(b_k)_k$, we will use it for the $\alpha$-numeration of $a$ ( that is also the $\alpha$-numeration of $1-\frac{b}{d}$). Moreover $s=\max(i,N_i>0)$.
 $$ t(S)=\sum_{k=k_0}^{\lfloor r/2\rfloor}b_{2k} \esp ; \esp e(S)=2+\nu+\sum_{k=1}^{\lfloor (s-1)/2\rfloor}a_{2k}$$
 \esp As in the previous case, we distinguish same two subcases and similar arguments :\\
  $\blacktriangleright\blacktriangleright$ Subcase 1 : if $b_k=0$ for all $k>s$. 
   $$ e(S)-t(S)= 2+\nu- c_s +\sum_{k=1}^{k_0-1}a_{2k}+\sum_{k=k_0}^{\lfloor (s-1)/2 \rfloor}(a_{2k}-b_{2k})$$
   \esp Where $c_s=b_s$ if $s$ is even and $0$ else. If $s$ is odd, then $\nu=0=c_s$. If $s$ is even, then $\nu=N_s$ and $c_s=b_s \in \{ N_s,N_s+1 \}$. In both cases, we have $\nu-c_s\geqslant -1$, so $e(S)-t(S)\geqslant 1$. \\
    $\blacktriangleright\blacktriangleright$ Subcase 2 : see Case 3 subcase 2. \\
    
    $\blacktriangleright$ Case 5 : if $d>\max(a,b)$ and $m<x_0+b\leqslant d$.\\
    \esp We are in the same situation as in Case 4, except that $e(S)$ is one unity less. So, we have to look carefully at this case. \\
    $\blacktriangleright\blacktriangleright$ Subcase 1 : if $b_k=0$ for all $k>s$. 
   $$ e(S)-t(S)= 1+\nu- c_s +\sum_{k=1}^{k_0-1}a_{2k}+\sum_{k=k_0}^{\lfloor (s-1)/2 \rfloor}(a_{2k}-b_{2k})$$
   \esp Where $c_s=b_s$ if $s$ is even and $0$ else. We remark that if $k_0>1$, then we can easily conclude, thanks to the term $a_2\geqslant 1$ and the Case 4. \\
   \esp If $s$ is odd, then $\nu=0=c_s$ and $e(S)-t(S)\geqslant 1$. \\
   \esp If $s$ is even : if $b_s=N_s$, we also have $e(S)-t(S)\geqslant 1$. If $b_s=N_s+1$, then $b_{2k}=a_{2k}$  and $b_{2k+1}=0$ for all $k<s/2$, except $b_1=1$. We will show that this  is not possible. Indeed, we have ( see \ref{subsec:alphanum} ) :
   $$ 1-\frac{b}{d}= \delta_0-\sum_{i=1}^{\frac{s}{2}-1}a_{2i}\delta_{2i-1}-(N_s+1)\delta_{s-1}= \cdots = \delta_{s-2}-(N_s+1)\delta_{s-1}$$
   \esp yet :
   $$x_0=d(\delta_{s-2}-N_s\delta_{s-1})$$
  \esp  so :
   $$ x_0+b=d(1+\delta_{s-1})>d$$
   \esp which contradicts our hypothesis. \\
   $\blacktriangleright\blacktriangleright$ Subcase 2 : again we refer to Case 3 Subcase 2. If $s$ is even, then with same arguments, $e(S)-t(S)\geqslant 1$. If $s$ is odd, we will show that it is not possible : indeed, we would have
   $$ (b_k)_k=(1,a_2,0,a_4,0,\cdots,0,a_{s-1},0,1)$$
   \esp so :
   $$ 1-\frac{b}{d}= \delta_{s-1}-\delta_s $$
  \esp yet, $x_0=d\delta_{s-1}$, so : $x_0+b=d+d\delta_s>d$ : contradiction.\\

  $\blacktriangleright $ Case 6 :  if $d>\max(a,b)$ and $x_0+b \leqslant m$. \\
  \esp Again, we use Theorem \ref{theo:5} and  \ref{theo:8}. Since some common notations name different things in each theorem, we will precise this : we use $k_1$ instead of $t$ in Theorem \ref{theo:5} and we do not use $s$ in Theorem \ref{theo:8}, since it names $s=\max(i,N_i>0)$ in Theorem \ref{theo:5}. We use $(b'_k)_k$ instead of $(b_k)_k$ in Theorem \ref{theo:5}, for $(b_k)_k$ names the $\alpha$-numeration of $a$ in Theorem \ref{theo:8}.\\
  \esp So $\alpha=m/d$, $\beta=\frac{x_0+b-1}{d}$ and $(b'_k)_k$ is the $\alpha$-numeration of  $\beta$. Again, $(b_k)_k$ is the $\alpha$-numeration of $a$.\\
  \esp  \textbf{N.B. : } the hypothesis $x_0+b \leqslant m$ is equivalent to $\beta<\alpha$, that is to say $b'_1=1$ and $b'_2>0$. \\ 
  \esp Moreover :
  $$ \sigma = \max(k,b_k>0) \esp ; \esp s'=\max(k, b'_k>0)$$
    \esp We mention that, since $(b_k)_k$ is the successor of $(N_k)_k$ for ALO in $E_{\alpha}$, then $\sigma=s$ in the general case or $\sigma=s+1$ in one exceptional case ( see Case 3 subcase 2).\\ 
  \esp We have with Theorems \ref{theo:5} and  \ref{theo:8} :
  $$ e(S)= 1+a_{2k_1}-b'_{2k_1}+\nu+\sum_{k=k_1+1}^{\lfloor (s-1)/2\rfloor}a_{2k} \esp ; \esp t(S)=\sum_{k=k_0}^{\lfloor r/2\rfloor}b_{2k} $$ 
  
  \esp We underline the definitions of $k_0$ and $k_1$ :
  $$ k_0=\begin{cases} \frac{\sigma+1}{2} \text{ if } b_{2i}=0 \text{ for all } i \\ \min(i>0,b_{2i}>0) \text{ else } \end{cases}  \esp ; \esp k_1   =\begin{cases} \frac{s'}{2} \text{ if } b'_{2i+1}=0 \text{ for all } i>0 \\ \min(i>0,b'_{2i+1}>0) \text{ else } \end{cases} $$

  \esp \textbf{Step 1} : we will prove that $k_0> k_1$.\\
  \esp   First, we remark that $b'_{2i+1}=0$ if $i\in\inter{1}{k_1-1}$ and $b'_{2i}=a_{2i}$ if $i<k_1$. In addition, $b'_{2k_1}>0$ if $k_1>1$. Indeed : if $b'_{2i+1}$ for all $i>0$, then $2k_1=s'$. Else, $b'_{2k_1+1}>0$ and $b'_{2k_1-1}=0$. Yet : $(b'_k)_k=(1,a_2,0,a_4,0,\cdots,0,a_{2k_1-2},0,b'_{2k_1},...)$ and so :
  $$ k_1>1 \Rightarrow (b'_k)_k <_A (1,a_2,0,a_4,0,\cdots,0,a_{2k_1-2})$$
 \esp  Then :
  $$ \frac{x_0+b-1}{d} <\delta_0-\sum_{i=1}^{k_1-1}a_{2i}\delta_{2i-1}= \delta_{2k_1-2}$$
  \esp But, this inequality is also true if $k_1=1$, since $\delta_0=\alpha$ ( see N.B. above).\\
  \esp On the other hand, using similar arguments for $(b_k)_k$, we  obtain :
  $$(b_k)_k \leqslant_A(a_1,0,a_3,\cdots,0,a_{2k_0-3},0,b_{2k_0-1})$$
  \esp Indeed : $b_{2i}=0$ and $b_{2i-1}=a_{2i-1}$  for $i\in\inter{1}{k_0-1}$. In addition, if $b_{2i}>0$ for some $i$, then $b_{2k_0}>0$ and we have a strict inequality. Else, $2k_0-1=\sigma$ and we have an equality !  \\
  \esp  So : 
   $$ 1-\frac{b}{d}\leqslant  b_{2k_0-1}\delta_{2k_0-2}+\sum_{i=0}^{k_0-2}a_{2i+1}\delta_{2i}=1-\delta_{2k_0-3} +b_{2k_0-1}\delta_{2k_0-2}$$
   \esp Hence :
   $$   1-x_0+d \delta_{2k_1-2} > b \geqslant  d(\delta_{2k_0-3}-b_{2k_0-1}\delta_{2k_0-2})$$
   \esp But we have  $x_0-1\geqslant 0$ and 
   $$\delta_{2k_0-3}-b_{2k_0-1}\delta_{2k_0-2}=(a_{2k_0-1}-b_{2k_0-1})\delta_{2k_0-2}+\delta_{2k_0-1}\geqslant \delta_{2k_0-1}$$
   \esp So :
   $$ \delta_{2k_1-2}> \delta_{2k_0-1}$$
   \esp Now, $(\delta_i)_i$ is non increasing, so $k_0 > k_1$. \\
   

\esp \textbf{Step 2 : } \\ 
    $\blacktriangleright\blacktriangleright$ Subcase 1 : if $b_k=0$ for all $k>s$. That is to say, if $\sigma=s$.\\
   $$ e(S)-t(S)= 1+\nu- c_s+a_{2k_1}-b'_{2k_1} +\sum_{i=k_1+1}^{k_0-1}a_{2i}+\sum_{i=k_0}^{\lfloor (s-1)/2 \rfloor}(a_{2i}-b_{2i})$$
   \esp where $c_s=b_s$ if $s$ is even and $0$ else. If $s$ is odd, then $\nu=0=c_s$ and $e(S)-t(S)\geqslant 1$. If $s$ is even, then $\nu=N_s$ and $c_s=b_s \in\{ N_s,N_s+1 \}$. If, $b_s=N_s$, same conclusion. If $b_s=N_s+1$, then, as in Case 5, subcase 1, we obtain $x_0+b>d$, which contradicts our hypothesis.\\

   $\blacktriangleright\blacktriangleright$ Subcase 2 : if  $(N_k)_k=(a_1,a_2,a_3,\cdots,a_{s})$. That is to say, if $\sigma=s+1$. Then $b_{s+1}=1$ and $b_{s-2k}=0$ for all $k\geqslant 0$  ( except $b_1=1$ if $s$ is odd) and $b_{s-2k+1}=a_{s-2k+1}$ if $k\in\inter{1}{\lfloor s/2 \rfloor}$. \\
     $$ e(S)-t(S)= 1+\nu- c_{s+1}+a_{2k_1}-b'_{2k_1} +\sum_{i=k_1+1}^{k_0-1}a_{2i}+\sum_{i=k_0}^{\lfloor (s-1)/2 \rfloor}(a_{2i}-b_{2i})$$  
 \esp  where $c_{s+1}=1$ if $s$ is odd and $c_{s+1}=0$ else. If, $s$ is even, then $\nu-c_{s+1}=N_s>0$, so $e(S)-t(S)\geqslant 1$. If $s$ is odd, this is not possible ( see Case 5, subcase 2). 
 \end{mademo}

\newpage


\normalsize 
\esp \textbf{All suggestions, corrections, critics, remarks, questions... are welcome and can be made at the following email adress : }  emmanuel.cabanillas@ac-montpellier.fr

 \end{document}